\documentclass[11pt]{amsart}
\usepackage{anyfontsize}
\usepackage{tikz-cd}
\usepackage[headings]{fullpage}
\usepackage{amsmath}
\usepackage{amsfonts}
\usepackage{latexsym}
\usepackage{amssymb}
\usepackage{placeins}

\newtheorem{proposition}{Proposition}[section]
\newtheorem{theorem}{Theorem}[section]

\newtheorem{lemma}{Lemma}[section]
\newtheorem{definition}{Definition}[section]
\newtheorem{corollary}{Corollary}[section]

\def\fin{\sf fin}
\def\vol{{\sf vol}}
\def\d{{\sf D}}
\def\grad{\mathop{\sf grad}}
\def\div{\mathop{\sf div}}
\def\homo{\mathop{\sf Hom}}
\def\ker{\mathop{\sf Ker}}

\def\im{\mathop{\sf Im}}

\def\curl{\mathop{\sf curl}}
\def\rk{\mathop{\sf rank}}
\def\det{\mathop{\sf det}}
\def\ann{\mathop{\sf ann}}
\def\dep{\mathop{\sf depth}}
\def\dim{\mathop{\sf dim}}
\def\codim{\mathop{\sf codim}}

\def\ass{\mathop{\sf Ass}}
\def\supp{\mathop{\sf supp}}
\def\pt{{\sf PT}}
\def\t{{\sf T}}

\def\dbydt{\sf \frac{d}{dt}}
\def\dbydts{\sf\frac{d^2}{dt^2}}
\def\ext{\sf Ext}

\def\Cinf{\mathcal{C}^\infty}
\def\Dr{\mathcal{D}'}

\def\Mellk{\mathcal{M}_{\ell,k}}
\def\Sellk{\mathcal{S}_{\ell,k}}
\def\Bellk{\mathcal{B}_{\ell,k}}
\def\Sc{\mathcal{S}}
\def\Sr{\mathcal{S}'}
\def\Er{\mathcal{E}'}
\def\D{\mathcal{D}}
\def\A{\mathcal{A}}
\def\C{\mathbb{C}}
\def\M{\mathcal{M}}

\newcommand{\del}{\partial}

\newcommand{\R}{\ensuremath{\mathbb{R}} }

\newcommand{\V}{\ensuremath{\mathcal{V}} }

\newcommand{\Rn}{\R^{n}}

\newcommand{\pDe}{[D_1, \ldots , D_n]}

\newcommand{\F}{\mathcal{F}}
\begin{document}

\title[Controllability and Vector Potential]{Controllability and Vector Potential}

\maketitle

\begin{center}
\author{Shiva Shankar
\footnote{\tt shunyashankar@gmail.com}
}
\end{center}
\begin{abstract} Kalman's fundamental notion of a controllable state space system \cite{k} has been generalised to higher order systems by Willems \cite{w}, and further to distributed systems defined by partial differential equations \cite{ps}. It turns out, that for systems defined in several important spaces of distributions,  controllability is now identical to the notion of vector potential in physics, or of vanishing homology in mathematics. 
These notes will explain this relationship, and a  few of its consequences. It will also pose an important question: does a controllable system, in any space of distributions, always admit a vector potential?  In other words, is Kalman's notion of a controllable system, suitably generalised, nothing more - nor less - than the possibility of describing the dynamics of the system by means of a vector potential? 

Furthermore, it also turns out that the category of distributed systems bears many formal similarities to the category of affine algebraic sets. This raises a second important question: what is the category for which these distributed systems are `local models', just as affine algebraic sets are local models for the category of algebraic varieties? It would then be possible to extend the theory of control described in these notes to this larger category of systems. 
\end{abstract}

\vspace{.8cm}

\tableofcontents

%\vspace{1cm}
\newpage

%\section{The Controllability Question}
\section{The Solvability Question for PDE, and its Dual}

\vspace{.2cm}

To bring into perspective the nature of the question that shall interest us here, we begin with a question that was of central importance in the 1950s and 60s.
\vspace{2mm}

\begin{center} The solvability question for systems of partial differential equations \end{center} 

\vspace{1mm}

Let \index{$A$, ring of differential operators}$A = \mathbb{C}\pDe$,
$D_j=\frac{1}{\imath}\partial_j$, $j = 1, \ldots, n$, be the $\mathbb{C}$-algebra of constant coefficient partial differential operators, and \index{$\mathcal{D}'$}$\Dr$ the space of distributions, on $\Rn$. We denote elements of $A$ by lower case letters $a, p, q, \ldots$, or by $a(D), p(D), q(D), $ for emphasis. (When we consider equations from Physics, such as the Laplace or the Maxwell equations, we take $A$ to be generated by $\partial_1, \ldots, \partial
_n$, and then we denote elements by $a(\partial), p(\partial)$ and so on. We also write $\C[\dbydt]$ for the ring of ordinary differential operators.) The ring $A$ acts on $\Dr$ by differentiation, and gives it the structure of a (topological) $A$-module. Let $\F$ be an $A$-submodule of $\Dr$. The attributes of the dynamical systems that we study will assume values in $\F$, and we call it the \index{signal space}{\it space of signals}. While many of the formal properties of systems described in these notes hold true for an arbitrary $A$-submodule $\F$, our focus will be mainly on $\Dr$, the space $\Cinf$ of smooth functions, the space $\Sr$ of temperate distributions, the Schwartz space $\Sc$ of rapidly decreasing functions, and the spaces $\Er$ and $\D$ of compactly supported distributions and smooth functions respectively on $\Rn$. We refer to them collectively as the \index{classical spaces}{\em classical signal spaces}.
\vspace{1mm}

Let $P(D) =  (p_{ij}(D))$ be an $\ell \times k$ matrix with entries in $A$; let its rows be denoted $r_1, \ldots, r_\ell$. It defines the  differential operator 
\begin{equation} 
\begin{array}{cccc}
P(D): & \F^k  & \longrightarrow & \F^\ell\\
& f=(f_1, \ldots, f_k) & \mapsto &  (r_1f, \ldots, r_\ell f) ,
 \end{array}
\end{equation}
where $r_if =  \sum_{j=1}^k p_{ij}(D)f_j, ~i = 1, \ldots , \ell$.
The solvability question for $P(D)$ in $\F$ asks: given $g \in \F^\ell$, is there an $f \in \F^k$ such that $P(D)f = g$? 
\vspace{1mm}

Consider the set $Q$ of all relations between the rows of $P(D)$, i.e. the set of all $(a_1, \ldots, a_\ell) \in A^\ell$ such that $a_1r_1 + \cdots + a_\ell r_\ell = 0$. $Q$ is an $A$-submodule of $A^\ell$, and as $A$ is Noetherian,  is finitely generated, say by $\ell_1$ elements. Let $Q(D)$ be the $\ell_1 \times \ell$ matrix whose rows are these $\ell_1$ generators. Then, the sequence
\begin{equation}
A^{\ell_1} \stackrel{Q^\t(D)}{\longrightarrow} A^\ell \stackrel{P^\t(D)}{\longrightarrow} A^k
\end{equation}
is exact by construction (the superscript {\tiny$\t$} denotes transpose). Applying the (left exact) functor\footnote{We refer to \cite{ma} and \cite{hor} for definitions in algebra and analysis  respectively.} 
 $\homo_A(-,~\F)$ to the above sequence, gives the sequence
\begin{equation}
\F^k \stackrel{P(D)}{\longrightarrow} \F^\ell \stackrel{Q(D)}{\longrightarrow} F^{\ell_1} .
\end{equation}
This sequence is a complex, i.e. the image of the map $P(D)$ is contained in the kernel of the differential operator $Q(D)$. It provides a necessary condition for the solvability question: to solve $P(D)f = g$, it is necessary that $Q(D)g = 0$.

Further, suppose this sequence is also exact, i.e. suppose that the image of the map $P(D)$ is {\em equal} to the kernel of $Q(D)$. Then, a  necessary {\em and} sufficient condition for the solvability of $P(D)f = g$ would be that $Q(D)g = 0$.
\vspace{1mm}
 
 {\it Thus, the solvability question for $P(D): \F^k \rightarrow \F^\ell$ admits an answer precisely when \phantom{Xx}its image  equals a kernel (namely that of the map $Q(D): \F^\ell \rightarrow \F^{\ell_1}$).} \\

 The problem of determining when images of partial differential operators are equal to kernels, was solved by Ehrenpreis, Malgrange, H\"ormander, Palamodov and others. Their solution is the description of the $A$-module structure of $\F$, but before we summarise it (for the classical spaces), we pose the dual question.
 \vspace{1mm}

Given $P(D)$ as above, now let $R$ be the set of all relations between its columns. $R$ is an $A$-submodule of $A^k$, and suppose that it is generated by $k_1$ elements, say. Let $R(D)$ be the $k \times k_1$ matrix whose columns are these generators. Then the sequence 
\begin{equation} A^\ell \stackrel{P^\t(D)}{\longrightarrow} A^k \stackrel{R^\t(D)}{\longrightarrow} A^{k_1}\end{equation}
is a complex (which is not necessarily exact). Again, applying the functor $\homo_A(-,~\F)$ to this sequence gives the complex
\begin{equation}\F^{k_1} \stackrel{R(D)}{\longrightarrow} \F^k \stackrel{P(D)}{\longrightarrow} \F^\ell . \end{equation}
Thus, the kernel of the map $P(D)$ contains the image of $R(D)$. We now ask:
\vspace{2mm}

\noindent ($\ast$) \hspace{.02cm} {\it What is the question about $P(D): \F^k \rightarrow \F^\ell$ for which the answer is `precisely \phantom{XXX}when its kernel equals an image (namely that of the map $R(D): \F^{k_1} \rightarrow \F^k$)'?} 
\vspace{1.5mm}

These notes are about this question, and the answer to which, for the classical spaces, is `the kernel of $P(D)$ is {\it controllable} precisely when it is an image' (Theorem 7.1). They are organised as follows: the next section explains the solution to the solvability question in terms of the $A$-module structure of the classical spaces, and introduces Willems' notion of controllability \cite{w}, a vast generalisation of Kalman's \cite{k}. The section following it provides necessary and sufficient conditions for controllability. The fourth section shows that a distributed system admits a decomposition into its controllable and uncontrollable parts. The fifth section constructs a coarse topology on infinite dimensional affine subspaces of the space of distributed systems. With respect to this topology, the Zariski topology, controllable systems are generic within the class of all underdetermined systems, whereas uncontrollable systems are generic within the class of all overdetermined systems. The sixth section points out the difficulties involved when the signal space is neither an injective nor a flat $A$-module. An important example here is the direct and inverse limits of the Sobolev spaces, the natural stage for questions related to energy and dissipation. The final section studies the Nullstellensatz question for PDE in the classical spaces, and introduces spaces of periodic functions on $\R^n$ for which the question has been answered.

Willems \cite{w} is a summary of the foundational work on which these notes rest, and \cite{bo, quad} shorter summaries. Pommaret \cite{pom} is a comprehensive monograph. There are also other interpretaions of the notion of a  controllable system which are not considered here, for instance the notion of minimality in \cite{l}. While the primary focus in these references and in these notes are questions in control theory, they have a wider applicability, as for instance in the construction of special `wave like' solutions from a vector potential in \cite{hhs}. The notion of a vector potential is fundamental in classical physics, but its utility in mathematical questions is yet to be fully explored.\\

%\newpage

\section{Controllable Dynamical Systems from Kalman to Willems}

\vspace{.2cm}

In 1960, Kalman \cite{k} introduced the notion of a \index{controllable state space system, Kalman}controllable state space system.
\begin{definition} A linear finite dimensional system whose evolution is described by
\begin{equation} {\dbydt} x = Xx + Uu, \end{equation}
where the state (or phase) $x:\R \rightarrow \mathbb{C}^\ell$, and \index{input}input $u: \R \rightarrow \mathbb{C}^m$, are smooth functions of time, and $X:\mathbb{C}^\ell \rightarrow \mathbb{C}^\ell, ~U:\mathbb{C}^m \rightarrow \mathbb{C}^\ell$ are linear maps, is said to be controllable if the following is satisfied: given $x_1, x_2$ in $\mathbb{C}^\ell$, and $t_1, t_2$ in $\R$ with $t_1 < t_2$, there is an input $u$ such that the solution to the above equation with this input satisfies $x(t_1) = x_1, ~x(t_2) = x_2$.
\end{definition}
This notion  was quickly recognised to be fundamental, and on it was built the superstructure of post-war control theory.
Nonetheless, there were foundational issues with the state space model itself. For one, the system's evolution had to be described by first order differential equations; another, more serious problem, was the need to specify, ab initio, a causal structure that declared some signals to be inputs. These, and other problems with the Kalman paradigm, were overcome in a spectacular fashion by J.C.Willems \cite{w}, who proposed a far larger class of models for linear dynamical systems, and a definition of controllability which faithfully generalised  Kalman's definition.  To motivate this development, we first make a few comments regarding the nature of models.
\vspace{1.5mm}

A model is a picture of reality, and the closer it is to reality, the better the picture it will be, and
more accurate will be the theory that describes this model. A model seeks to represent a
certain phenomenon whose attributes are certain qualities that vary with space,
time etc. The closest we can get to this reality, this phenomenon, is to take all possible variations
of the attributes of the phenomenon, itself, as the model. This collection, considered all together
in our minds, is in engineering parlance, the {\it  behaviour} of the phenomenon. 
Indeed, paraphrasing Wittgenstein, we could have said: 
%\vspace{.3mm}

\begin{center} {\it The behaviour is all that is the case.} \end{center}

\vspace{1.5mm}
A priori, any variation of these attributes could perhaps have occured, but the laws governing the
phenomenon limit the variations to those that can, and do, occur. We shall consider phenomena described by local
laws, i.e. laws expressed by differential equations. Here by local we mean that the variation of the
attributes of the phenomenon at a point depends on the values of the attributes in arbitrarily small
neighbourhoods of the point, and not on points far away.

This is a familiar situation in physics, mathematics and engineering. For instance:
\vspace{2mm}

1. Planetary motion. A priori, the earth could perhaps have traversed any trajectory around the
sun, but Kepler's laws restrict it to travel an elliptic orbit, sweeping equal areas in equal times, and
such that the square of the period of revolution is proportional to the cube of the major semi-axis.
Kepler's laws are of course local, namely Newton's equations of motion.
\vspace{2mm}

2. Magnetic fields. A priori, any function $B$ could perhaps have been a magnetic field, but in fact
must satisfy the law that there do not exist magnetic monopoles in the universe. By Gauss, this law is also
local:
$\div B = 0 $.
\vspace{2mm}

3. Complex analysis. Complex analysis studies those functions which admit a convergent power
series about each point. This is the law that defines the subject.

This law is again local, namely the Cauchy-Riemann equations.\\

We now further confine ourselves to phenomena which are described by local laws that
are also linear and shift invariant. The behaviour of such phenomena can be modelled by kernels of linear maps defined by constant coefficient partial differential operators.

Thus, for example, suppose that the attributes of a phenomenon are described by some $k$-tuple of complex numbers, which depend smoothly on some $n$ independent real valued parameters. Let these attributes be denoted by $f = (f_1, \ldots , f_k)$, where each $f_i: \R^n \rightarrow \mathbb{C}$ is smooth. Then, a law the phenomenon obeys is a partial differential operator
$p(D) = (p_1(D), \ldots , p_k(D))$ such that $p(D)f =\sum p_i(D)f_i = 0$. In other words a law is an operator $p(D): (\mathcal{C}^\infty)^k \rightarrow \mathcal{C}^\infty$, and to say that $f$ obeys this law is to say that $p(D)f = 0$, or that $f$ is a zero
of $p(D)$. The behaviour $B$ of the phenomenon is thus 
\[  B ~ = \bigcap_{{\rm all ~laws} ~p(D)} \{f ~| ~p(D)f = 0 \},
\]
the common zeros of all the laws governing the phenomenon. Shift invariance implies that the entries $p_1(D), \ldots,p_k(D)$ of $p(D)$ are all in $A = \mathbb{C}[D_1, \ldots ,D_n]$. As the sum of two laws, and the scaling of a law by an element of $A$, are also laws, it follows that the set $P$ of all laws governing the phenomenon is an $A$-submodule of $A^k$. If $P$ is generated by $r_1, \ldots ,r_\ell$, then writing these $\ell$ elements of $A^k$ as rows of a matrix, say $P(D)$, implies that $B$ is the kernel of
\[ P(D): (\mathcal{C}^\infty)^k \rightarrow (\mathcal{C}^\infty)^\ell. \] %which is Equation (1) for $\F = \mathcal{C}^\infty$. \\

Thus, the phenomena we study are exactly those whose behaviours can be modelled by kernels of constant coefficient partial differential operators (as in Equation (1)).

These kernels define the class of \index{distributed system}{\it distributed systems} (or \index{distributed behaviour}{\it distributed behaviours)}. When $n = 1$, ~i.e. when $A = \mathbb{C}[\dbydt]$, they define the class of  \index{lumped systems}{\it lumped systems}.\\

\noindent Remark 2.1: Control of such a system is the process of restricting its behaviour to a suitable subset by the imposition of additional laws, i.e. by enlarging the rows of $P(D)$. This formulation of control does not involve notions of inputs and outputs, and hence avoids any a priori notion of causality.\\   
%These notes will not however address this issue any further.

\noindent Example 2.1: Let $A = \mathbb{C}[\dbydt]$. The collection of all signals $x$ and $u$ in Kalman's state space system of Definition 2.1, is the kernel of the operator
\begin{equation}
\begin{array}{lccc}
P({\dbydt}) = \left ({\dbydt} I_\ell - X, ~-U \right ): & (\mathcal{C}^\infty)^{\ell + m} & \longrightarrow & (\mathcal{C}^\infty)^\ell \\
  & (x,u)\phantom{X} & \mapsto & P(\dbydt)\left [ \begin{array}{l} x\\u \end{array}   \right ] ~,
 \end{array}
\end{equation}
but now, the difference between state and input has been obliterated. In other words, we do not impose any causal structure on the system. \\

\noindent Remark 2.2: The kernel $\ker_\F(P(D))$ of the operator of Equation (1) depends upon the submodule $P \subset A^k$ generated by the rows of the matrix $P(D)$, and not on a specific choice of generators of $P$ such as the rows of $P(D)$. Indeed, the maps
\[ 
\begin{array}{ccc}
\homo_A(A^k/P,~ \F)  & \longrightarrow & \ker_\F (P(D))\\
 \phi & \mapsto & (\phi([e_1]), \ldots , \phi([e_k]))  
 \end{array}
\]
where, $[e_1], \ldots ,[e_k]$ denote the images of the standard basis $e_1, \ldots ,e_k$ of $A^k$ in $A^k/P$, and
\[
\begin{array}{ccc}
\ker_\F(P(D))  & \longrightarrow & \homo_A(A^k/P, ~\F)  \\ 
 f=(f_1, \ldots, f_k) & \mapsto & \phi_f   
 \end{array}
\]
where $\phi_f([e_i]) = f_i$, $1 \leqslant i \leqslant k$, are inverses of each other.
Hence, we denote this kernel also by $\ker_\F(P)$, and call it the kernel of $P$ in $\F$.\\

The following statements are immediate.
\begin{lemma}
(i) If $P_1 \subset P_2$ are submodules of  $A^k$, then $\ker_\F(P_2) \subset \ker_\F(P_1)$ in $\F^k$.  ($\ker_\F(P_2)$ is said to be a \index{sub-system}sub-system (or \index{sub-behaviour}sub-behaviour) of $\ker_\F(P_1)$.)

\vspace{1mm}
\noindent(ii) If $\{P_i\}$ is  a collection of submodules of $A^k$,  then $\ker_\F(\sum_i P_i) = \bigcap_i \ker_\F(P_i)$, and $\sum_i\ker_\F(P_i) \subset \ker_\F(\bigcap P_i)$.
\end{lemma}

\vspace{2mm}

In our more general setting, controllability cannot any longer be the ability to move from one state to another in finite time  (as in Definition 2.1), because there is now no notion of state. We need another definition altogether, which specialises to the Kalman definition for state space systems of Example 2.1. This new definition (2.3 below) is the insight of J.C.~Willems, and the genesis of his `Behavioural Theory'. \\

\begin{definition} Let $\F$ be an $A$-submodule of $\mathcal{D}'$. 
\vspace{1.5mm}

\noindent (i) A disjoint pair of subsets $U_1, U_2$ of $\mathbb{R}^n$ whose closures do not intersect, is said to be \index{admissible pair}admissible with respect to $\F$, if for any two signals $f_1, f_2$ in $\F$, there is an $f$ in $\F$ such that $f = f_1$ on some neighbourhood of $U_1$ and $f = f_2$ on some neighbourhood of $U_2$.
 
Such an $f$ is said to patch $f_1$ on $U_1$ with $f_2$ on $U_2$.
\vspace{1.5mm} 

\noindent (ii) A nested pair of subsets $U \subset V$ of $\Rn$ such that the closure of $U$ is contained in the interior of $V$, is admissible with respect to $\F$, if for any $f \in \F$ there is an $f_c \in \F$ such that $f_c = f$ on some neighbourhood of $U$ and equal to 0 on a neighbourhood of the complement $V^c$ of $V$. 

Such an $f_c$ is said to be a cutoff of $f$ with respect to $U \subset V$ in $\F$.
\end{definition}

The disjoint pair $U_1, U_2$ in (i) is admissible if the nested pair $U_1 \subset U_2^c$ is admissible as in (ii). Similarly, the nested pair $U \subset V$ in (ii) is admissible if the disjoint pair $U, V^c$ is admissible as in (i). The two definitions are thus equivalent.\\

\noindent Example 2.2: If $U = \{(x_1,x_2) \in \mathbb{R}^2 |~ x_2 < 0\}$ and $V = \{(x_1,x_2) \in \mathbb{R}^2 | ~x_2 < e^{-x_1^2}\}$, then the closure of $U$ is contained in the interior of $V$. Let $f$ be any function in the Schwartz space $\Sc$ which is nonzero on $\mathbb{R}^2$ and which decays as $e^{-|x|}$ at infinity. Then no cutoff of $f$ with respect to this $U \subset V$ can be in $\mathcal{S}$ (the derivatives of any such cutoff would not be rapidly decreasing). Similarly, no cutoff of the constant function $1$ in $\Sr$ with respect to the above $U \subset V$ can be tempered.

The nested pair $U \subset V$ is thus not admissible with respect to either $\mathcal{S}$ or $\Sr$.

\begin{proposition} (i) Any  $U \subset V$, with the closure of $U$ contained in the interior of $V$, is admissible with respect to $\Dr$, $\Cinf$, $\Er$ or $\D$.

\vspace{1mm}
\noindent (ii) Suppose that the distance between $U$ and the complement $V^c$ of $V$ is bounded away from 0, i.e. $\|x-y\| \geqslant \epsilon > 0$ for all $x \in U, y \in V^c$ (for instance, suppose the boundary of $U$ is compact). Then such a pair is admissible with respect to $\Sc$ or $\mathcal{S}'$.
\end{proposition}
\noindent Proof: (i) Let $\rho$ be any smooth function which equals 1 on some neighborhood $U'$ of $U$ contained in the interior of $V$, and 0 on the complement of $V$. Then for an element $f$ in $\mathcal{C}^\infty$ or in $\D$, $\rho f$ is a cutoff of $f$ with respect to $U \subset V$.

Now let $g$ be in $\Dr$ or in $\Er$. As $\rho g(f) = g(\rho f)$ for every $f \in \D$ or $\Cinf$ respectively, it follows that $\rho g$ is the required cutoff.

\vspace{2mm}
(ii) Let $U'$ be an open subset containing $U$ and contained in the interior of $V$, such that the distance between $U'$ and the complement of $V$ is $\frac{\epsilon}{2}$. Let $\chi$ be the characteristic function of $U'$, and let $\kappa$ be a smooth `bump' function supported in the ball of radius $\frac{\epsilon}{4}$ centered at the origin, i.e. $\kappa$ is identically 1 in some smaller ball about the origin. Then the convolution $\rho = \kappa \star \chi$ is a smooth function that is identically 1 on $U'$, 0 on the complement of $V$, and such that all its derivatives are bounded on $\mathbb{R}^n$. Thus for any $f$ in $\Sc$ or $\Sr$, $\rho f$ is the required cutoff. 
\hspace*{\fill}$\square$\\

\begin{definition} (J.C.Willems \cite{w}) A system  $\ker_\F(P) \subset \F^k$, where $P$ is a submodule of $A^k$ and $\F$ a submodule of $\Dr$, is \index{controllable system, Willems}controllable if given any two subsets $U_1$ and $U_2$ of $\R^n$ admissible with respect to $\F$, and any two elements $f_1$ and $f_2$ in $\ker_\F(P)$, then there is an element $f$ in $\ker_\F(P)$ such that $f = f_1$ on some neighbourhood of $U_1$ and $f = f_2$ on some neighbourhood of $U_2$. 

Equivalently, $\ker_\F(P)$ is controllable if and only if given any $f$ in it, and any nested $U \subset V$ admissible with respect to $\F$, then there is an $f_c$ in $\ker_\F(P)$ such that $f_c = f$ on some neighbourhood of $U$, and equal to 0 on the complement of $V$.
\end{definition}

%\vspace{1.5mm}
Thus, controllability of a system asserts that any two signals in it can be patched within the system itself, or equivalently, that a signal in it admits a cutoff within it. For emphasis, we sometimes refer to this notion as \index{behavioural controllability}{\it behavioural controllability}. 

%The following propositions are elementary.

\begin{proposition} Let $\F$ be a classical signal space (with its standard topology). Let $\ker_\F(P)$ be a controllable system as in the above definition. Then, its compactly supported elements  are dense in it. 
\end{proposition}
\noindent Proof: Let $V_1  \subset \cdots \subset V_i \subset V_{i+1} \subset \cdots$ be an exhaustion of $\R^n$ by compact sets such that $V_i$ is contained in the interior of $V_{i+1}$. Let $f$ be any element in $\ker_\F(P)$, and let $f_i \in \ker_\F(P)$ be a cutoff of $f$ with respect to $V_i \subset V_{i+1}$. Then the sequence $\{f_i\}$ converges to $f$.  
\hspace*{\fill}$\square$\\

\begin{proposition}  Consider the system given by the kernel of the operator $P(D):\F^k \rightarrow \F^\ell$ in Equation (1). If it is equal to the image of an operator $R(D):\F^{k_1} \rightarrow \F^k$ (i.e. if sequence (5) is exact), then it is controllable.
\end{proposition}
\noindent Proof: Given data $U \subset V \subset \R^n$ as in Definition 2.2 (ii) and $f \in \ker_\F(P(D))$, let $g \in \F^{k_1}$ be any element such that $R(D)g = f$. Let $g_c$ be any cutoff (in $\F^{k_1}$) with respect to $U \subset V$. Then $R(D)g_c$ is a cutoff of $f$ with respect to $U \subset V$ in $\ker_\F(P(D))$. \hspace*{\fill}$\square$\\

\noindent Example 2.3: The collection of all static magnetic fields is the system defined by Gauss' Law, and is controllable because it admits the vector potential \index{curl}$\curl$. In other words, the following sequence 
\[ (\Cinf)^3 \stackrel{\curl}{\longrightarrow} (\Cinf)^3 \stackrel{\div}{\longrightarrow} \Cinf \]

\vspace{1mm}
is exact. Here,  \index{divergence}$\div = (\del_1, \del_2, \del_3)$,   and  $\curl$ = $\left(\begin{array}{ccc}0 & -\del_3 & \partial_2 \\\partial_3 & 0 & -\partial_1 \\-\partial_2 & \partial_1 & 0\end{array}\right)$. 

\vspace{5mm}
This example prompts the following definition \cite{pom}.

\begin{definition} A system given by the kernel of a differential operator is said to admit a \index{vector potential}vector potential if it equals the image of some differential operator. 
\end{definition}
 
The primary goal of these notes is to prove that the converse of the above proposition is also true for a system defined in a classical space (Theorem 7.1 below). 
\vspace{1.5mm}

\noindent {\it Thus, in the classical spaces, the notion of a controllable system provides the answer to question ($\ast$) posed at the end of the previous section.}
\vspace{1.5mm}

The argument depends upon the answers to the solvability question that we started our discussion with.  We first summarise  these classical results due to Malgrange \cite{m1,m2} and Palamodov \cite{p}; they describe the $A$-module structure  of the classical spaces.
\vspace{1.5mm}

It remains an important open question whether a distributed system in any signal space $\F$ is controllable if and only if it admits a vector potential.

\vspace{4mm}

\begin{center} The Algebraic Structure of the Classical Spaces of Distributions
\end{center}
\vspace{1mm}

\begin{definition} An $A$-module $\F$ is \index{injective module}injective if $\homo_A (- ,~\F)$ is an exact functor. It is \index{flat module}flat if $~ - \otimes_A \F$ is an exact functor.
\end{definition}
Thus, if $\F$ is an injective $A$-submodule of $\Dr$, then the sequence (3) above is always exact, and we can answer the solvability question for every $P(D)$. 

\begin{theorem} The spaces $\Dr$, $\Cinf$, and \index{$\mathcal{S}'$}$\Sr$ are injective $A$-modules.
\end{theorem}
This is the celebrated \index{Fundamental Principle}Fundamental Principle of Malgrange and Palamodov.
\vspace{1.5mm}

\noindent Remark 2.3: The elementary case of the above theorem, namely that these spaces are divisible $A$-modules, are the classical existence theorems of Ehrenpreis-Malgrange for $\Dr$ and $\Cinf$, and of  H\"ormander for $\Sr$. They assert that in these spaces, $p(D): \F \rightarrow \F$ is surjective for any nonzero $p(D)$. In particular, any nonzero $p(D)$ has a tempered fundamental solution. 
\vspace{1.5mm}

Their duals, the spaces \index{$\mathcal{D}$}$\D$, \index{$\mathcal{E}'$}$\Er$ and \index{$\mathcal{S}$}$\Sc$ are not injective, not even \index{divisible module}divisible $A$-modules. For instance, the map $\dbydt: \D(\R) \rightarrow \D(\R)$ is not surjective,  as $f \in \D$ is in the image if and only if its integral over $\R$ equals 0. It is also not surjective when $\D$ is replaced by $\Er$ or $\Sc$.

\begin{theorem} \cite{m1} The spaces $\D$, $\Er$, and $\Sc$ are flat $A$-modules.
\end{theorem}

\noindent Remark 2.4: Sometimes we will think of the ring $A$ as the polynomial ring $\mathbb{C}[x_1, \ldots ,x_n]$. We will then identify the differential operator $p(D)$ with the corresponding polynomial $p(x)$, and  consider its affine variety $\V(p(x))$ in $\mathbb{C}^n$. We call it the affine variety of $p(D)$, or simply, the affine variety of $p$. A point $\xi \in \mathbb{C}^n$ is in $\V(p(x))$ if and only if the  exponential $e^{\imath<\xi, x>}$ is a solution of $p(D)$. More generally, the \index{exponential solution}exponential solutions of $p(D)$ are solutions of the form $q(x)e^{\imath{<\xi, x>}}$, where $q(x)$ is a polynomial, and $\xi$ a point in $\V(p(x))$. A theorem of Malgrange states that the set of solutions of $p(D)$ in $\Dr$ is the closed convex hull of its exponential solutions.

We also identify a matrix $P(D)$ with the corresponding matrix $P(x)$ by identifying each entry $p_{ij}(D)$ with the polynomial $p_{ij}(x)$.

\begin{definition} An injective $A$-module $\F$ is a \index{injective cogenerator}cogenerator if for every nonzero $A$-module $M$, $\homo_A(M,~ \F)$ is nonzero. A flat $A$-module $\F$ is \index{faithfully flat}faithfully flat if for every nonzero $A$-module $M$, $M \otimes_A \F$ is nonzero. 
\end{definition}

\begin{proposition} The $A$-modules $\Dr$ and $\Cinf$ are cogenerators.
\end{proposition}
\noindent Proof: Let $M$ be a nonzero $A$-module. Let $m(D)$ be any nonzero element of $M$ and
let $(m)$ be the cyclic submodule generated by it. Thus $(m)$ is isomorphic to
$A/i$ for some proper ideal $i$ of $A$.
As $\Cinf$ is an injective $A$-module, $\homo_A(M, ~\Cinf )$ surjects onto $\homo_A(A/i, ~\Cinf) \simeq \ker_{\Cinf}(i)$.
Therefore to prove the proposition it suffices to show that $\ker_{\Cinf}(i)$ is nonzero. But this is elementary, for as
$i$ is a proper ideal of $A$, its variety $\V(i)$ (in $\mathbb{C}^n$) is nonempty. Let $\xi$ be any point in it.
Then $m(D)e^{\imath{<\xi,x>}} = 0$, and so $\ker_{\Cinf}(i) \neq 0$.
\hspace*{\fill}$\square$\\

\noindent Remark 2.5: The injective module $\Sr$ is not a cogenerator. For example, let $A=\mathbb{C}[\dbydt]$. Consider the ideal generated by $\dbydt - 1$; then $\homo_A(A/(\dbydt-1), ~\Sr) = 0$ ($e^t$ is not tempered).

\begin{proposition} $\D$ and $\Er$ are faithfully flat $A$-modules.
\end{proposition}
\noindent Proof: A flat $A$-module $M$ is faithfully flat if and only if $mM \neq M$, for every maximal ideal $m$ of $A$. 

The Fourier transforms of the distributions
in $m\Er$, which extend to functions holomorphic on $\C^n$ by the \index{Paley-Wiener Theorem}Paley-Weiner Theorem, all vanish at the zero of $m$.
Thus $m\Er \neq \Er$ for every maximal ideal $m$, hence $\Er$ is faithfully flat. Similarly, $\D$ is also faithfully flat.
\hspace*{\fill}$\square$\\

\noindent Remark 2.6:
Let $p(x) = 1 + x_1^2 + \cdots + x_n^2$; then $p(x)^{-1} f$ is in $\Sc$, for every $f \in \Sc$, This implies (by Fourier transformation) that
$p(D): \Sc \rightarrow \Sc$ is surjective, so that $m\Sc = \Sc$ for any maximal ideal $m$ of $A$ that contains $p(D)$. Thus, the flat $A$-module $\Sc$ is not faithfully flat. 
\vspace{3mm}

This concludes the description of the classical spaces that these notes rest on. We observe a few immediate consequences.

\begin{proposition} Suppose that the $A$-submodule $\F \subset \Dr$ is an injective cogenerator. Then there is an inclusion reversing  bijection between $A$-submodules $P \subset A^k$ and distributed systems $\ker_\F(P) \subset \F^k$.
\end{proposition}
\noindent Proof: We have already observed in Lemma 2.1 that the corespondence between submodules and systems is inclusion reversing. Now suppose that $M \subsetneq N$, then $0 \rightarrow N/M \rightarrow A^k/M \rightarrow A^k/N \rightarrow 0$ is exact. Applying the functor $\homo_A(-, ~\F)$ to this sequence yields the exact sequence  $0 \rightarrow \ker_\F(N) \rightarrow \ker_\F(M) \rightarrow \homo_A(N/M, ~\F) \rightarrow 0$. As $\F$ is a cogenerator, $\homo_A(N/M, ~\F) \neq 0$, hence $\ker_\F(N) \subsetneq \ker_\F(M)$.

Next suppose that $M \not\subset N$; then $N \subsetneq M + N$, hence $\ker_\F(M) \cap \ker_\F(N) = \ker_\F(M+N) \subsetneq \ker_\F(N)$ by the above paragraph, hence $\ker_\F(M) \neq \ker_\F(N)$.  This establishes the proposition.
\hspace*{\fill}$\square$\\

\noindent Remark 2.7: For $\F$ any $A$-submodule of $\Dr$, and any distributed system $\ker_\F(P) \subset \F^k$, let $\M(\ker_\F(P))$ $
\subset A^k$ be the submodule of all the elements in $A^k$ that map every element in $\ker_\F(P)$ to 0. Clearly $P \subset \M(\ker_\F(P))$; indeed it is the largest submodule of $A^k$ whose kernel equals $\ker_\F(P)$.
The determination of $\M(\ker_\F(P))$ is the {\it Nullstellensatz problem for systems of  PDE} in $\F$ \cite{sn}. The above proposition avers that if $\F$ is an injective cogenerator, then $P = \M(\ker_\F(P))$, for every submodule $P \subset A^k$.

More generally, for any subset $X$ of $\F^k$ we set $\M(X)$ to be the submodule of all elements $p(D) \in A^k$ such that the kernel of $p(D): \F^k \rightarrow \F$ contains $X$.
\vspace{1mm}

The assignment $\M$ is also inclusion reversing, i.e. if $B_1 \subset B_2$ are two distributed systems in $\F^k$, then $\M(B_2) \subset \M(B_1)$. Thus we have two inclusion reversing assignments, $\ker_\F$ and $\M$, which define a Galois connection between the partially ordered sets of submodules of $A^k$ and systems in $\F^k$.
\vspace{1.5mm}

\noindent Example 2.4: If $\F$ is not an injective cogenerator, then $\M(\ker_\F(P))$ may not be equal to $P$. For instance, both the ideals $(1)$ and $(\dbydt - 1)$ of $\mathbb{C}[\dbydt]$ define the 0 system in $\Sr$, hence $\M(\ker_{\Sr}(\dbydt-1)) = (1)$. 
\vspace{1.5mm}

An argument, similar to the one in the above proposition, proves the following:
\begin{proposition} Suppose that $\F$ is an injective cogenerator. Then the complex $M \stackrel{f}{\rightarrow} N \stackrel{g}{\rightarrow} P$ is exact if and only if $\homo_A(P, ~\F) \stackrel{-\circ g}{\longrightarrow} \homo_A(N, ~\F) \stackrel{-\circ f}{\longrightarrow} \homo_A(M, ~\F)$ is exact. Hence, $0 \rightarrow M \rightarrow N \rightarrow P \rightarrow 0$ is split exact if and only if $0 \rightarrow \homo_A(P, ~\F) \rightarrow \homo_A(N, ~\F) \rightarrow \homo_A(M, ~\F) \rightarrow 0$ is split exact.
\end{proposition}
\noindent Proof: Suppose $\im (f) \subsetneq \ker (g) \subset N$. Then, by injectivity of $\F$, ~ $\homo_A(N/\im(f), ~\F)$ maps onto $\homo_A(\ker(g)/\im(f), ~\F)$, where the second $A$-module is nonzero because $\F$ is also a cogenerator.  Thus,  there exists an $A$-module map $\phi: N \rightarrow \F$ which restricts to 0 on $\im (f)$, but which is nonzero on $\ker (g)$. This $\phi$ maps to 0 in $\homo_A(M, ~\F)$, as $\phi \circ f = 0$. However, it is not  in the image of $-\circ g$, as $\phi$ is nonzero on $\ker (g)$.

The second assertion now follows immediately. 
\hspace*{\fill}$\square$\\

\noindent Example 2.5: The above proposition is not true for an injective $A$-module which is not a cogenerator. For instance $0 \rightarrow (\dbydt-1) \rightarrow \mathbb{C}[\dbydt] \rightarrow \mathbb{C}[\dbydt]/(\dbydt - 1) \rightarrow 0$ is not split exact, yet $ 0 \rightarrow \ker_{\Sr}(\dbydt - 1) \rightarrow \Sr {\rightarrow}$ $\homo_A(\dbydt - 1, ~\Sr) \rightarrow 0$ splits because $\ker_{\Sr}(\dbydt -1) = 0$.
\vspace{3mm}

\noindent Remark 2.8: We have seen that $\D$, $\Er$ and $\Sc$ are flat $A$-modules.  By the \index{equational criterion for flatness}equational criterion for flatness \cite{ma}, the kernel of every $P(D): \F^k \rightarrow \F^\ell$ is an image, i.e. sequence (5) is always exact if $\F$ is flat. Thus, every distributed system in these spaces admits a vector potential and is controllable.
\vspace{1mm}

We shall therefore not study systems in these spaces in any detail, confining ourselves only to the construction of counter examples to statements that are true for systems in injective modules, but not in general. (We refer to \cite{stop} for a discussion on the relationship between the structure of a topological $A$-module and its dual.) \\

\newpage

\section{Necessary and Sufficient Conditions for Controllability}

\vspace{.2cm}

Our goal is to answer question ($\ast$) of Section 1.

\vspace{1mm}
We make the following standing assumption about the $A$-module $\F$~:

\noindent  Let $V$ be a bounded subset of $\R^n$. Then there exists a subset $U$ whose closure is contained in the interior of $V$ such that the pair $U \subset V$ is admissible with respect to $\F$ (as in Definition 2.2(ii)).

Thus, every $f \in \F$ admits a cutoff $f_c \in \F$ with respect to such a pair.

\vspace{1mm}
By Proposition 2.1, the classical signal spaces $\Dr, ~\Cinf, ~\Sr, ~\Sc, ~\Er$, and $\D$,  satisfy this assumption.

\vspace{1mm}
Elementary consequences of this assumption are the following corollaries.

\begin{corollary} An $\F$ satisfying the above assumption is a faithful $A$-module. In other words, if $p(D)$ is a nonzero element of $A$, then the kernel of $p(D): \F \rightarrow \F$ is strictly contained in $\F$. \end{corollary} 
\noindent Proof: The above assumption implies that there are compactly supported elements in $\F$. No such element can be a homogenous solution of $p(D)$ by the Paley-Wiener Theorem. \hspace*{\fill}$\square$
\vspace{1mm}

\begin{corollary} The only sub-systems of $\F$ that are controllable are $0$ and $\F$.
\end{corollary}
\noindent Proof: The sub-system 0, which is the kernel of the map $1:\F \rightarrow \F$, is trivially controllable. So is $\F$, the kernel of $0:\F \rightarrow \F$, as it admits the vector potential $1:\F \rightarrow \F$. 
 
Suppose $p(D) \neq 0$ is such that the kernel of $p(D): \F \rightarrow  \F$ is  nonzero. By the standing assumption above, there is a relatively compact pair $U \subset V$ admissible with respect to $\F$. No nonzero $f \in \ker_\F(p(D))$ can admit a cutoff  with respect to this pair in $\ker_\F(p(D))$.
\hspace*{\fill}$\square$\\

\begin{proposition} Sequence (4) is exact if and only if $A^k/P$ is \index{torsion free}torsion free. 
\end{proposition}
\noindent Proof: To recollect notation: $P$ is the submodule of $A^k$ generated by the rows of the $\ell \times k$ matrix $P(D)$ of Equation (1), $R \subset A^k$ is the submodule of relations between the columns of $P(D)$ and it is generated by the $k_1$ columns of $R(D)$.
Let the kernel of $R^\t(D): A^k \rightarrow A^{k_1}$ be $P_1 \subset A^k$; clearly $\im(P^\t(D)) = P \subset P_1$. Thus (4) is a complex, and is exact if and only if $P = P_1$. 

Let $K$ be the field of fractions of the doman $A$. Tensoring complex (4) by $K$ gives the complex $K^\ell \stackrel{P^\t(D)}{\longrightarrow} K^k \stackrel{R^\t(D)}{\longrightarrow} K^{k_1}$ of finite dimensional $K$-vector spaces. 
%as the $A$-module maps $P^\t(D)$ and $R^{\t}(D)$ are also $K$-linear. 
As localisation is a flat functor, the image of $P^{\t}(D)$ now is $K \otimes_A P$, and the kernel of $R^\t(D)$ is now  $K \otimes_A P_1$. If the rank of $P^\t(D)$ equals $r$, then the rank of $R^\t(D)$ equals $k-r$ and its kernel has dimension $r$. Thus $K \otimes_A P = K \otimes_A P_1$, hence $K \otimes_A P_1/P = 0$ and so $P_1/P$ is a torsion module.

We next claim that the torsion submodule of $A^k/P$ equals $P_1/P$. For suppose $x \in A^k$ is such that $ax \in P$ for some nonzero $a \in A$. Then $1 \otimes ax  \in K \otimes_A P = K \otimes_A P_1$, hence $1 \otimes ax = \alpha_1 \otimes p_1 + \cdots + \alpha_r \otimes p_r$, for some $\alpha_i \in K$ and $p_i \in P_1$, for all $i$. Clearing the denominators of the $\alpha_i$, we have $1 \otimes bax = a_1 \otimes p_1 + \cdots + a_r \otimes p_r$, where the $b, a_1, \ldots ,a_r$ are all in $A$. As $P_1$ has no torsion, this equality implies that $bax = a_1p_1 + \cdots + a_rp_r \in P_1$. Thus, $R^{\t}(D)(bax)=0$, and as $ba \neq 0$, $R^{\t}(D)x=0$, which is to say that $x$ is itself in $P_1$. 
It follows that (4) is exact, i.e. $P = P_1$, if and only if $A^k/P$ is torsion free.
\hspace*{\fill}$\square$\\

\noindent Example 3.1: \noindent  Let $A= \mathbb{C}[\del_1, \del_2, \del_3]$, and $m$ its maximal ideal $(\del_1, \del_2, \del_3)$. Consider the following sequence
\begin{equation}0 \rightarrow A \stackrel{\div^\t}{\longrightarrow} A^3 \stackrel{\curl^\t}{\longrightarrow} A^3 \stackrel{\grad^\t}{\longrightarrow} A \stackrel{\pi}{\longrightarrow} A/m \rightarrow 0\end{equation}
where $\pi$ is the canonical projection, and \index{gradient}$\grad = \div^{\t}$ (Example 2.3). As $A$ is an integral domain, the map $\div^\t $ is an injection. The module of  relations between the columns of the matrix $\curl$ is a cyclic module generated by the one column of $\grad$, and the module of relations between the columns of $\div$ is generated by the columns of $\curl$.  Let $C_r$ and $D_r$  be the submodules of $A^3$ generated by the rows of $\curl$ and $\div$ respectively. As $A^3/C_r$ and $A^3/D_r$ are both torsion free, it follows  from the above proposition that the above sequence is exact, and hence a resolution of $A/m$.
\vspace{1.5mm}

We can now begin the programme to answer question ($\ast$) of Section 1.
\begin{theorem}(\cite{ps}) Let $\F \subset \Dr$ be an $A$-submodule, and $P$  a submodule of $A^k$.

\noindent (i) If $\F$ is an injective $A$-module, then $\ker_\F(P)$ admits a vector potential if $A^k/P$ is torsion free, and is hence controllable. 

\noindent (ii) If $\F$ is an injective cogenerator, then $\ker_\F(P)$ is controllable if and only if $A^k/P$ is torsion free. Thus $\ker_\F(P)$ is controllable if and only if it admits a vector potential.
\end{theorem}
\noindent Proof: (i) If $A^k/P$ is torsion free, then sequence (4) is exact by the above proposition, and this  in turn implies sequence  (5) is exact because $\F$ is an injective $A$-module. Thus $\ker_\F(P)$ admits a vector potential, and is hence controllable (by Proposition 2.3).

\vspace{1mm}
(ii) Now suppose that $\F$ is an injective cogenerator, and suppose that $A^k/P$ is not torsion free. Let $q(D) \in A^k \setminus P$ be such that $a(D)q(D) \in P$, for some nonzero $a(D) \in A$. Denote by $q$ and $aq$ the cyclic submodules of $A^k$ generated by $q(D)$ and $a(D)q(D)$ respectively. Then $P \subsetneq P + q$, hence $\ker_\F(P) \not\subset \ker_\F(q)$, by Lemma 2.1 and Proposition 2.6. However, $\ker_\F(P) \subset \ker_\F(aq)$, as $aq \subset P$.

Let $f  \in \ker_\F(P) \setminus \ker_\F(q)$. Let $V$ be a bounded open subset of $\R^n$ such that $q(D)f \neq 0$ on $V$. Let $U$ be a closed subset of $V$ such that the pair $U \subset V$ is admissible with respect to $\F$, namely the standing assumption above. Suppose $\ker_\F(P)$ were controllable, then some cutoff $f_c$ of $f$ with respect to $U \subset V$ would be in $\ker_\F(P)$. Clearly $q(D)f_c \neq 0$, and it has compact support.

As $a(D)q(D) \in P$, it follows that $a(D)(q(D)f_c) = 0$. This is a contradiction, as $a(D)$ cannot have
 a compactly supported homogeneous solution (by Paley-Wiener).
\hspace*{\fill}$\square$\\

\noindent Remark 3.1: If an injective $A$-submodule $\F$ of $\Dr$ is not a cogenerator, then $\ker_\F(P)$  could be an image, and hence  controllable, even though $A^k/P$  has torsion. For instance  $A/(\dbydt - 1)$ is a torsion module,  yet $ \Sr \stackrel{0}{\longrightarrow} \Sr \stackrel{\dbydt-1}{\longrightarrow} \Sr$ is exact. 
\vspace{2mm}

\noindent Example 3.2  (\index{deRham complex}deRham complex on $\R^3$): (i)  Suppose that $\F$ is $\Dr$, $\Cinf$ or $\Sr$. As these are injective $A$-modules, applying the functor $\homo(-, ~\F)$ to resolution (8) of Example 3.1 yields the exact sequence
\[0 \rightarrow \mathbb{C} \stackrel{\pi^\t}{\longrightarrow} \F \stackrel{\grad}{\longrightarrow} \F^3 \stackrel{\curl}{\longrightarrow} \F^3 \stackrel{\div}{\longrightarrow} \F \rightarrow 0,\]
where $\pi^\t$ is the injection of $\homo_A(A/m, ~\F) \simeq \mathbb{C}$, the space of constant functions, into $\F$.  Manifestly $\ker_\F(\div)$ and $\ker_\F(\curl)$ are controllable. 

The ideal of $A$ generated by the rows of $\grad$ is the maximal ideal $m = (\del_1, \del_2, \del_3)$. As $A/m$ is a torsion module, $\ker_\F(\grad)$ is not controllable when $\F$ is $\Dr$ or $\Cinf$, by Theorem 3.1 (ii). (Clearly we cannot patch together distinct constant functions within the space of constant functions.) This is also true when $\F = \Sr$.
\vspace{1mm}

\noindent (ii) Suppose now that $\F$ is $\D$, $\Er$ or $\Sc$. These are flat $A$-modules, hence tensoring  (8) with $\F$ yields the exact sequence \[ 0 \rightarrow \F \stackrel{\grad}{\longrightarrow} \F^3  \stackrel{\curl}{\longrightarrow} \F^3 \stackrel{\div}{\longrightarrow} \F \stackrel{\pi}{\longrightarrow} \F/m\F \rightarrow 0 \]
Now $\ker_\F(\grad)$ is also controllable because it equals 0 (namely Remark 2.8).
\vspace{1.5mm}

We now have two definitions of a  controllable  state space system, namely Kalman's Definition 2.1, and Willems' Definition 2.3 when the system (6) is rewritten as equation (7) of Example 2.1. We show below that the two definitions coincide; thus Willems' definition is a perfect generalisation of Kalman's.

We use the \index{PBH Test}Popov-Belevitch-Hautus (PBH) test for Kalman controllability of state space systems \cite{pw}.
\vspace{2mm}

\noindent PBH Test: Let $P(x) = (x I_\ell - X, -U)$, where $P(\dbydt) \in \mathbb{C}[\dbydt]$ is the operator of equation (7).  Then the state space system (6) of Definition 2.1 is controllable if and only if $P(\xi)$ is of full rank for every $\xi \in \mathbb{C}$. 
\vspace{1.5mm}

Thus, to determine whether (6) is controllable, it suffices to calculate all the maximal minors of $P(x)$, there are ${\ell + m \choose \ell}$
of them, and then to check that these minors do not have a common zero in $\mathbb{C}$. In other words, is the $\ell$-th determinantal ideal $\frak{i}_\ell(P(x))$, generated by the determinants of the $\ell \times \ell$ submatrices of $P(x)$, equal to $(1)$?
\vspace{1.5mm}

We generalise the above statement on the $\ell$-th determinantal ideal to the case of $A = \mathbb{C}[D_1,\ldots ,D_n]$. Henceforth, we adopt the following convention: given a submodule $P \subset A^k$, let $\ell \geqslant 0$ be such that $P$ can be generated by $\ell$ elements, but not by any set of $\ell - 1$ elements.  Then, any set of $\ell$ generators for $P$ is said to be {\it minimum}. We will always represent $P$ by an $\ell \times k$ matrix $P(D)$, where the $\ell$ rows is a \index{minimum set of generators}minimum set of generators. We call this integer $\ell$ the \index{minimum number of generators}minimum number of generators for $P$.

Given such a $P(D)$, consider the ideal of $\mathbb{C}[x_1, \ldots, x_n]$ generated by all the $\ell \times \ell$ minors of $P(x)$. If $\ell > k$, this ideal is equal to 0, otherwise it is the $(k-\ell)$-th Fitting ideal of $A^k/P$, and  is therefore independent of the choice of the matrix $P(D)$ whose $\ell$ rows is a minimum set of generators for $P$. It is the \index{cancellation ideal}{\it cancellation} ideal of $P$, or of $P(D)$ or $\ker_\F(P)$, and we denote it by $\frak{i}_\ell(P)$.

\vspace{1mm}
Similarly, the ideal generated by all the $k \times k$ minors of $P(x)$, where $P(D)$ is now any matrix whose rows generate $P$, is the 0-th Fitting ideal of $A^k/P$. It is the \index{characteristic ideal}{\it characteristic} ideal of $P$, or of $\ker_\F(P)$, and is denoted $\frak{i}_k(P)$.
\vspace{1mm}

The affine varieties $\mathcal{V}(\frak{i}_\ell(P))$ and $\mathcal{V}(\frak{i}_k(P))$ in $\mathbb{C}^n$ are the \index{cancellation variety}cancellation and \index{characteristic variety}characteristic varieties of $\ker_\F(P)$, respectively.

\begin{definition}Suppose that the submodule $P \subset A^k$  can be generated by not more than $k$ elements, i.e. $\ell \leqslant k$ above, then $\ker_\F(P)$ is \index{underdetermined system}underdetermined. It is \index{strictly underdetermined system}strictly underdetermined if $\ell < k$. If $\ell > k$, $\ker_\F(P)$ is an \index{overdetermined system}overdetermined system.
\end{definition}
Thus, the characteristic ideal $\frak{i}_k(P)$ equals 0 if $\ker_\F(P)$ is strictly underdetermined, and the cancellation ideal $\frak{i}_\ell(P)$ equals 0 if it is overdetermined. 
\vspace{1mm}

When $A=\mathbb{C}[\dbydt]$, a principal ideal domain, every submodule of $A^k$ is free (being finitely generated and torsion free) and has a basis not larger in number than $k$. Then, every matrix $P(\dbydt)$ representing $P$ as above, i.e. whose rows is now a basis for $P$, will be of full (row) rank. 
In this case, the PBH test generalizes perfectly to arbitrary lumped systems.

\begin{proposition} (\cite{w}) Let $\F$ be an injective cogenerator, and $P$ a submodule of $(\mathbb{C}[{\dbydt]})^k$ of rank $\ell$. Then $A^k/P$ is free (of rank $k - \ell$) if and only if the PBH condition $\frak{i}_\ell(P) = (1)$ is satisfied.  Thus, the PBH condition is necessary and sufficient for $\ker_\F(P)$ to be controllable. 
\end{proposition}
\noindent Proof: The structure theory for modules over a PID implies that there is a basis $e_1, \ldots , e_k$ for $A^k$,  and elements $a_1({\dbydt}), \ldots ,  a_\ell({\dbydt})$ of $A$, such that $a_1({\dbydt})e_1, \ldots, \linebreak a_{\ell} ({\dbydt}) e_\ell$ is a basis for $P$, namely the Smith form of any matrix $P(\dbydt)$ whose $\ell$ rows generate $P$. Thus, for $\xi \in \mathbb{C}$, the rank of $P(\xi)$ equals $\ell $ if and only $a_1(\xi), \dots , a_\ell(\xi)$ are all nonzero. This latter condition holds for every $\xi$ in $\mathbb{C}$ if and only if they are all nonzero constants. In other words, $P(x)$ has full row rank for every $\xi \in \mathbb{C}$ if and only if a basis for $P$ extends to a basis for $A^k$.

This last condition is in turn equivalent to saying that the exact sequence
\[
0 \rightarrow P \longrightarrow A^k \longrightarrow A^k/P \rightarrow 0
\]
splits. This implies that $A^k/P$ is a submodule of $A^k$, hence free. By Theorem 3.1 (ii), this is equivalent to the controllability of $\ker_\F(P)$.
\hspace*{\fill}$\square$\\

Specializing the proposition to the case of Example 2.1 shows that Willems' definition of behavioural controllability is a faithful generalization of Kalman's state controllability. 

\begin{corollary} (\cite{w}) A state space system is controllable in the sense of Definition 2.1 if and only if it is controllable in the sense of Definition 2.3. 
\hspace*{\fill}$\square$ \end{corollary}

We return to the case of $A = \mathbb{C}[D_1, \ldots ,D_n]$, and ask the following questions:

\noindent Let the rows of the matrix $P(D)$ be a minimum set of generators for $P \subset A^k$. What does it mean to say that $P(x)$ has full row rank for each $x \in \mathbb{C}^n$? Is this still a necessary and sufficient condition for controllability of $\ker_\F(P)$? If not, what is now the analogue of the PBH test?
\begin{proposition} Let $\F$ be an injective $A$-module, $P$ a submodule of $A^k$ and let $\ell \geqslant 0$ be the cardinality of a minimum set of generators for $P$. Then $A^k/P$ is free (of rank $k - \ell$) if and only if $\frak{i}_\ell(P) = (1)$.  Thus, the PBH condition is sufficient for $\ker_\F(P)$ to be controllable. 
\end{proposition}
\noindent Proof: Clearly, both the implications in the above statement imply that $\ell \leqslant k$; for suppose that $A^k/P$ is free, then $A^k \rightarrow A^k/P \rightarrow 0$ splits, hence
$P$ is a direct summand of $A^k$, hence a projective $A$-module.  Therefore, by Quillen-Suslin,\footnote{The theorem of Quillen-Suslin asserts that a finitely generated projective module over the polynomial ring is free.} $P$ is free, of rank $\ell \leqslant k$. The other statement, $\frak{i}_\ell(P) = 1$, manifestly implies the inequality. We may also assume that $\ell > 0$, because otherwise $P$ is the 0 submodule of $A^k$, and then $\frak{i}_0(P)$ is the $k$-th Fitting ideal of $A^k$, which by definition equals $(1)$.

Let $P(D) = (p_{ij})$ be any $\ell \times k$ matrix whose $\ell$ rows is a minimum set of generators for $P$, and $P(x)$ the corresponding matrix with polynomial entries. Now, $\rk(P(\xi)) = \ell$ for $\xi = (\xi_1, \ldots ,\xi_n) \in \mathbb{C}^n$ if and only if there is a minor, say the determinant $\det(M(\xi))$ of an $\ell \times \ell$ submatrix $M(\xi)$ of $P(\xi)$, which is not 0.
As the value $p_{ij}(\xi)$ of $p_{ij}(x)$ at $\xi$  is the image of $p_{ij}(D)$ under the morphism $A \rightarrow A/m_\xi$,
$m_\xi = (D _1 - \xi _1, \ldots , D_n - \xi _n)$, it follows that $\det(M(D))$ does not belong to $m_\xi$, and is therefore a unit in the localization $A_{m_\xi}$. Hence
there is a basis such that the matrix of the localization $P_{m_\xi} (D): A^k_{m_\xi} \rightarrow A^\ell_{m_\xi}$ has an $\ell \times \ell$ submatrix equal to the identity $I_\ell$. By row and column operations, all 
other entries of $P_{m_\xi}(D)$ can be made zero. This implies that
\[
0 \rightarrow P_{m_\xi} \longrightarrow A^k_{m_\xi} \longrightarrow A^k_{m_\xi}/P_{m_\xi} \rightarrow 0
\]
splits. Thus $A^k_{m_\xi}/P_{m_\xi}\simeq (A^k/P)_{m_\xi}$ is projective, and as this is true for every maximal ideal $m_\xi$, it follows that $A^k/P$ is projective. Again by Quillen-Suslin, $A^k/P$ is free (of rank $k - \ell$).

The second statement now follows from Theorem 3.1(i).
\hspace*{\fill}$\square$
\vspace{1.5mm}

\noindent Remark 3.2: The quotient $A^k/P$ is free if and only if it, and thus also $P$, are direct summands of $A^k$. Then $A^k \simeq P \oplus A^k/P$, hence $\ker_\F(P) \simeq \homo_A(A^k/P,~ \F)$ is a direct summand of $\F^k \simeq \homo_A(A^k, ~\F)$, $\F$ any $A$-submodule of $\Dr$.  Such a system is said to be \index{strong controllability}{\it strongly controllable}. Thus PBH is a sufficient condition for the strong controllability of $\ker_\F (P)$, for $\F$ any $A$-submodule of $\Dr$.  

In particular, a lumped controllable system is always strongly controllable (namely Proposition 3.2).

Now suppose that $\F$ is an injective cogenerator. Then $\ker_\F(P)$ is a direct summand of $\F^k$ if and only if $P$ is a direct summand of $A^k$. Indeed, by Proposition 2.7,   splittings of $A^k \rightarrow A^k/P \rightarrow 0$ are in bijective correspondence with splittings of $ 0 \rightarrow \ker_\F(P) \rightarrow \F^k$. Thus the PBH condition is both necessary and sufficient for strong controllability when $\F = \Dr$, or $\Cinf$.
\vspace{1mm}

\noindent Remark 3.3: In Example 3.1, $\frak{i}_\ell(D_r)$ is the maximal ideal $(\del_1, \del_2, \del_3)$, hence $D_r$ does not satisfy the PBH condition. Thus while $\ker_{\Dr}(\div)$ is controllable, it is not strongly controllable.  
\vspace{1mm}

The above proposition is a statement about a projection of the system $\ker_\F(0) \simeq \F^k$. Below, we study projections of a general system. It is the PDE analogue of the Elimination Problem of classical algebraic geometry.
\vspace{1mm}

The next result is the generalisation of the PBH test to $A = \mathbb{C}[D_1, \ldots ,D_n]$, under the condition that the cancellation ideal be nonzero. (This condition is always satisfied by state space systems of Example 2.1.) Later, we show that this condition is generically satisfied in the space of underdetermined distributed systems.

\begin{theorem} Let $P$ be a submodule of $A^k$, and let $P(D)$ be any $\ell \times k$ matrix whose $\ell$ rows generate $P$. Suppose that the cancellation ideal $\frak{i}_\ell(P) \neq 0$  (so that in particular $\ell \leqslant k$). Then $A^k/P$ is torsion free if and only if $P(x)$ has full row rank for all $x$ in the complement of an algebraic variety in $\mathbb{C}^n$ of dimension $\leqslant n-2$ (or in other words, that the Krull dimension of the ring $A/\frak{i}_\ell(P)$ be less than or equal to $n-2$). Thus,

\noindent (i) if $\F$ is an injective $A$-module, then $\ker_\F(P)$ is controllable if the dimension of the cancellation variety, $\mathcal{V}(\frak{i}_\ell(P))$, is less than or equal to $n-2$.

\noindent (ii) If $\F$ is an injective cogenerator, then $\dim(\mathcal{V}(\frak{i}_\ell(P))) \leqslant n-2$ is also necessary for $\ker_\F(P)$ to be controllable. \end{theorem}
\noindent Proof: Suppose to the contrary that $\rk(P(x))$ is less than $\ell$ for $x$ on an algebraic variety of dimension
$n-1$ in $\mathbb{C}^n$ - it cannot be less than $\ell$ on all of $\mathbb{C}^n$ because $\frak{i}_\ell(P)$ is nonzero, by assumption. Let $V$ be an irreducible component, and let it be the zero locus of the irreducible polynomial $p(x)$ (a codimension 1 prime ideal in a unique factorization domain is a principal ideal, generated by a prime element). This means that each of the $k\choose \ell$ generators of $\frak{i}_\ell(P)$ is divisible by $p(D)$. Let $p$ be the prime ideal $(p(D))$, and let $P_p(D)$ be the image of the matrix $P(D)$ in the localization $A_p$ of $A$ at $p$. Suppose that every entry of some row of $P_p(D)$, say the first row, is divisible by $p(D)$. As divisibility by $p(D)$ in $A$ and $A_p$ are equivalent, the corresponding row of $P(D)$ is also divisible by $p(D)$, and then clearly $A^k/P$ has a torsion element. Otherwise, at least one element of this first row of $P_p(D)$, say the first element, is not divisible by $p(D)$, and is therefore a unit in $A_p$. By row and column operations, every other element of the first column, and of the first row, can be made zero. Let the matrix so obtained after these row and column operations, with a unit in the (1,1) entry and the other entries of the first row and column equal to 0, be denoted $P^1_p(D)$, and let $P^1_p$ be the submodule of $A^k_p$ generated by its rows. 

The above argument for $P_p(D)$ can be repeated now for $P^1_p(D)$, and it follows that either $p(D)$ divides every element of some row, say the second, or that some element in that row, say  the (2,2) entry (the (2,1) entry is zero) is a unit in $A_p$. In the first case, $A^k_p/P^1_p$ has a torsion element, hence so does $A^k/P$; otherwise all other elements in the second row and second column can be made zero by column and row operations. Eventually, after at most $\ell-1$ such steps, the resultant matrix has units in the $(j,j)$ entries, $1\leqslant j \leqslant \ell-1$, and the other entries zero except in positions $(\ell,j), \ell \leqslant j \leqslant k$. If now $p(D)$ does not divide every entry of the $\ell$-th row, then it also does not divide the generators of the cancellation ideal $\frak{i}_\ell(P)$, which is a contradiction.

\vspace{1mm}
Conversely, suppose that $A^k/P$ has a torsion element. This implies that there is an element $x(D) = (a_1(D), \dots ,a_k(D))$ in $A^k \setminus P$ such that $r(D) = p(D)x(D) = (p(D)a_1(D), \ldots ,p(D)a_k(D))$ is in $P$, for some nonzero irreducible $p(D) \in A$. As $r(D)$ is in $P$, it is an $A$-linear combination of the $\ell$ rows $r_1(D), r_2(D), \ldots ,r_\ell(D)$ of $P(D)$, say $r(D) = b_1(D)r_1(D) + \cdots + b_\ell(D) r_\ell(D)$. Clearly the $b_j(D)$ are not all zero, and are also not all divisible by $p(D)$, because otherwise it would imply that $x(D)$ belongs to $P$, contrary to its choice. Thus, without loss of generality, let $b_1(D)$ be nonzero and not divisible by $p(D)$. Now let $B(D)$ be the followin $\ell \times \ell$ matrix: its first row is $(b_1(D), \ldots ,b_\ell(D))$, it has 1 in the $(j,j)$ entries, $2 \leqslant j \leqslant \ell$, and all other entries are 0. Then the product $B(D)P(D)$ is an $\ell \times k$ matrix whose first row is $r(D)$ and whose other rows are the rows $r_2(D), \ldots ,r_\ell(D)$ of $P(D)$. By construction of $B(D)$, its determinant $b_1(D)$ is not divisible by $p(D)$, and as every generator of the $\ell$-th determinantal ideal of $B(D)P(D)$ is divisible by $p(D)$, it follows that every generator of the cancellation ideal $\frak{i}_\ell(P)$ is also divisible by $p(D)$. This implies that $\rk(P(x)) < \ell$ at points $x$ in $\mathbb{C}^n$ where $p(x) = 0$, an algebraic variety of dimension $n-1$.

Statements (i) and (ii) now follow from Theorem 3.1.
\hspace*{\fill}$\square$
\vspace{1mm}

We have seen in Proposition 3.3 that the cancellation ideal $\frak{i}_\ell(P)$ of $P \subset A^k$ equals $(1)$ if and only if $P$ is a direct summand of $A^k$. The condition that $\frak{i}_\ell(P)$ be nonzero in the above theorem also admits an elementary interpretation.

\begin{proposition} Let $P$ be a submodule of $A^k$. Then its cancellation ideal $\frak{i}_\ell(P)$ is nonzero if and only if $P$ is free (of rank $\ell$).
\end{proposition}
\noindent Proof: Let $P(D)$ be an $\ell \times k$ matrix whose $\ell$ rows is a minimum set of generators for $P$; then $\frak{i}_\ell(P)$ is generated by the $\ell \times \ell$ minors of $P(x)$. If $\ell > k$, then there is nothing to be done, and so we assume $\ell \leqslant k$.
\vspace{.5mm}

Now if $P$ is not free, then there is a nontrivial relation between the rows of $P(x)$, hence all the maximal minors equal 0. 
\vspace{.5mm}

Conversely, suppose that $\frak{i}_\ell(P) = 0$. This means that after localizing at the 0 ideal, i.e over the function field $K = \mathbb{C}(x_1, \ldots , x_n)$, the $\ell$ rows of $P(x)$ span a subspace $L$ of $K^k$, whose projections to $K^\ell$, given by choosing $\ell$ of the $k$ coordinates, are all subspaces of dimension strictly less than $\ell$. If $k = \ell$, then the rows of $P(x)$ are $K$-linearly dependent, hence there is a relation between the rows of $P(x)$, and $P$ is not free (as $\ell$ is the minimum number of elements needed to generate it). Otherwise, there are ${k \choose \ell} > \ell$ such projections,  
hence $L$ must be of dimension strictly less than $\ell$, so that again, by definition of $\ell$, $P$ is not free.
\hspace*{\fill}$\square$\\

The following corollary is immediate.
\begin{corollary} Suppose that $\F$ is an injective cogenerator, and suppose that the minimum number of elements required to generate the submodule $P \subset A^k$ is equal to $k$ (i.e. suppose $\ell = k$ in our notation). If $\ker_\F(P)$ is controllable, then $P$ is not free. 
\end{corollary}
\noindent Proof: The cancellation ideal $\frak{i}_k(P)$ is a principal ideal of $A$.  If $P$ is free, then $\frak{i}_k(P) \neq 0$ by the above proposition. Hence, the cancellation variety has dimension $n-1$, and $\ker_\F(P)$ is not controllable by Theorem 3.2 (ii).
\hspace*{\fill}$\square$\\
 
\noindent Example 3.3: Neither the above corollary, nor Theorem 3.2, provide answers when $P$ is not free. For instance, let $A = \mathbb{C}[D_1, D_2, D_3]$, and let $P_1$, $P_2$ the submodules of $A^3$ generated by the rows of the matrices
\[
P_1(D) = \left(\begin{array}{ccc}0 & -D_3 & D_2 \\D_3 & 0 & -D_1 \\-D_2 & D_1 & 0\end{array}\right)
\hspace{3mm}  {\rm and} \hspace {4mm} P_2(D) = \left(\begin{array}{ccc}0 & -D_3 & D_2 \\D_3 & 0 & -D_1 \\-D_1 D_2 & D_1^2 & 0\end{array}\right) ~.
\]
The three rows of each matrix is a minimal set of generators for $P_1$ and $P_2$, and the cancellation ideals  $\frak{i}_3(P_1), \frak{i}_3(P_2) $, are both equal to 0. While $A^3/P_1$ is torsion free, $A^3/P_2$ has torsion; hence, $\ker_\F(P_1)$ is controllable, whereas $\ker_\F(P_2)$ is not.
\vspace{4mm}

We now impose a `causal structure' on a system in the sense that we identify some of its attributes to be the `controller`, and the others to be the ones to be `controlled'. This proceedure rests on an {\it elimination problem for PDE}.
\vspace{1mm}

Let $k = k_1 +k_2$. Consider the split exact sequence
\[ 0 \rightarrow A^{k_1}  \begin{array}{l} \stackrel{i_1}{\longrightarrow}\\ \stackrel{\pi_1}{\longleftarrow} \end{array} ~ A^k  ~ \begin{array}{l} \stackrel{\pi_2}{\longrightarrow}\\ \stackrel{i_2}{\longleftarrow} \end{array} A^{k_2} \rightarrow 0 ~ ,
\]
and the corresponding split exact sequence
\[ 0 \rightarrow \F^{k_1}  \begin{array}{l} \stackrel{\pi_1}{\longleftarrow}\\ \stackrel{i_1}{\longrightarrow} \end{array} ~ \F^k  ~ \begin{array}{l} \stackrel{i_2}{\longleftarrow}\\ \stackrel{\pi_2}{\longrightarrow} \end{array} \F^{k_2} \rightarrow 0 ~.
\]
%Let  $P$ be an $A$-submodule of $A^k$, and $\ker_\F(P) \subset \F^k$ the corresponding differential kernel. 

\begin{proposition} Let $P \subset A^k$ and $\F$  an $A$-submodule of $\Dr$. Then
\begin{equation} {\ker}_\F(\pi_2(P)) =  i_2^{-1}({\ker}_\F(P))  \subset \pi_2({\ker}_\F(P)) \subset {\ker}_\F(i_2^{-1}(P)). \end{equation}  %(and similarly for the corresponding submodules of $\F^{k_2}$).
If $\F$ is injective, then $\pi_2({\ker}_\F(P)) = {\ker}_\F(i_2^{-1}(P))$,
%(and similarly $\pi_2({\ker}_\F(P)) = {\ker}_\F(i_2^{-1}(P))$), 
so that a projection of a distributed system is also a distributed system.
\end{proposition}
\noindent Proof: %The first chain of equality and inclusions in the statement is elementary. For instance, suppose that $p_2(D) \in i_2^{-1}(P)$. Then $(0,p_2(D)) \in P$, so that $(0,p_2(D))(f) = 0$ for all $f \in {\ker}_\F(P)$. It follows that $p_2(D)(\pi_2(f)) = 0$, hence $\pi_2(\ker_\F(P)) \subset \ker_\F(i_2^{-1}(P))$. 
The first split sequence above implies that the sequence
\[ 0 \rightarrow A^
{k_2}/i_2^{-1}(P) \stackrel{i_2}{\longrightarrow} A^k/P \stackrel{\pi_1}{\longrightarrow} A^{k_1}/\pi_1(P) \rightarrow 0\] is exact, and hence so is $0 \rightarrow {\ker}_\F(\pi_1(P)) \stackrel{i_1}{\longrightarrow} {\ker}_\F(P) \stackrel{\pi_2}{\longrightarrow} {\ker}_\F(i_2^{-1}(P))$ exact. This proves the last inclusion in (9). The other inclusion and equality in it are also elementary.

If $\F$ is injective, then 
${\ker}_\F(P) \stackrel{\pi_2}{\longrightarrow} {\ker}_\F(i_2^{-1}(P)) \rightarrow 0$ is also exact. This completes the proof of the propositon.
\hspace*{\fill}$\square$\\

\noindent Example 3.4: If $\F$ is not injective, then the projection of a  system may not be a system. For instance, let $A=\mathbb{C}[\dbydt]$, and let $\F$ be $\D$. Let $\pi_2: \D^2 \rightarrow \D$ be the projection onto the second factor. Let $P \subset A^2$ be the cyclic submodule generated by $(\dbydt, -1)$. Then $\ker_\D(P) = \{(f, {\dbydt}f)~|~ f \in \D \}$,  but $\pi_2(\ker_\D(P))  = \{{\dbydt}f ~|~ f \in \D \}$ is not a differential kernel in $\D$, as ${\dbydt}: \D \rightarrow \D$ is not surjective.

Thus the inclusion $\pi_2(\ker_\F(P)) \subset \ker_\F(i_2^{-1}(P))$ in Proposition 3.5 can be strict when $\F$ is an not an injective $A$-module. In general, the obstruction to the projection $\pi_2(\ker_\F(P))$ being a differential kernel, lies in \index{Ext}${\ext}^1_A(A^{k_1}/\pi_1(P), ~\F)$. \\

\noindent Example 3.5: Let $A = \mathbb{C}[\del_x, \del_y]$. Consider the Cauchy-Riemann equation
\[ 
\begin{array}{lccc} 
 \left( \begin{array}{lc} \del_x & -\del_y \\ \del_y & \del_x  \end{array} \right):  & (\Cinf)^2 &\longrightarrow & (\Cinf)^2 \\& (u,v) & \mapsto & (\del_xu - \del_yv, ~\del_yu + \del_xv) \end{array}
\]
An elment $(u,v)$ in its kernel defines the holomorphic function $u + i v$. 

If $P \subset A^2$ is the submodule generated by the rows of the C-R matrix above, then $i_2^{-1}(P)$ is the principal ideal $(\del_x^2 + \del_y^2)$ generated by the Laplacian. As $\Cinf$ is an injective $A$-module, the above proposition states that the projection $\pi_2: (u,v) \mapsto v$    
maps the kernel of the C-R equation onto the kernel of the Laplacian, namely the set of harmonic functions. \\

\noindent Remark 3.4: In modelling the attributes $f \in \F^{k_1}$ of a phenomenon, it is sometimes convenient to include additional variables $g \in \F^{k_2}$. For instance, if we wish to model the voltages that appear across various branches of a circuit, it might be convenient to also include the currents through these branches.  
These additional variables are the \index{latent variables}`latent' variables, and the attributes of the phenomenon are the \index{manifest variables}`manifest' variables of the model \cite{w}. The modelling process then leads to equations of the form 
\[ P(D)f = Q(D)g~. \]
If $\F$ is an injective $A$-module, then the collection of all $f$ that appear in solutions $(f, g)$ of the above equation is also a differential kernel by the above proposition. In other words, we can eliminate the latent variables $g$ from the model. \\

\noindent Remark 3.5: Sometimes, it might be sufficient only for some projection of $\ker_\F(P)$, say $\pi_2(\ker_\F(P)) \subset \F^{k_2}$ in the above notation, to be controllable. In the case that $\F$ is an injective cogenerator, Proposition 3.5 gives necessary and sufficient condtions for this to be so, namely that $A^{k_2}/i_2^{-1}(P)$ be torsion free. \\

As we have observed, an engineering problem might dicatate that a subset of the attributes of a system be controlled by its other attributes. We illustrate this problem with a most important example, the Maxwell Equations, in the appendix below.

\vspace{1mm}

\subsection{Appendix: Control of the Electromagnetic Field}
 
\vspace{.2cm} 
In the above notation, we now assume that the attributes defined by the $k_1$ coordinates are the controller, and those defined by the $k_2$ coordinates are the ones that are to be controlled. We assume that $\F$ is injective, hence $B = \ker_\F(P)$ defines two systems, the controller $B_1 = \pi_1(B)$ and the system $B_2 = \pi_2(B)$ which is to be controlled. The control problem we now address is the following:
\vspace{1.5mm}

{\em Given the system $B \subset \F^{k}$, the problem is to restrict $B_2$ to a desired subsystem by restricting the behaviour of the controller $B_1$. This is to be achieved by augmenting the laws the controller must satisfy.
What are the subsystems of $B_2$ that can be thus achieved?}
\vspace{2mm}

We start with a few elementary observations.

\begin{lemma}  The subsystem $i_2^{-1}(B)$ of $\pi_2(B)$ is unchanged by additions to the controller laws.  
\end{lemma}
\noindent Proof: Imposing additional controller laws translates to specifying an $A$-submodule of $A^{k_1}$ strictly containing the submodule $i_1^{-1}(P)$. These laws correspond to laws of the form $(p_1, 0) \in A^k$ that are not in $P$. The submodule $P'$ generated by $P$ and these new laws satisfies $\pi_2(P') = \pi_2(P)$. As $\ker_\F(\pi_2(P)) = i_2^{-1}(B)$, this subsystem of $B_2$ remains unchanged when $P$ is enlarged to $P'$.  \hspace*{\fill}$\square$\\

Thus, the possible subsystems of $\F^{k_2}$ that can be achieved by restricting $\pi_1(B)$ with additional controller laws, all contain $i_2^{-1}(B)$, and are contained in $\pi_2(B)$. 

Similarly, every subsystem of $\F^{k_1}$ that can restrict $\pi_2(B)$ is contained in $\pi_1(B)$ and can be assumed to contain $i_1^{-1}(B)$.
\vspace{2mm}

These observations show that the attributes that are to be controlled, and the attributes with which control is to be accomplished, satisfy identical restrictions. The above control problem is thus symmetric with respect to this designation. 
\vspace{2mm}

\noindent Definition: An $A$-submodule of $A^{k_1}$  containing $i_1^{-1}(P)$, and which is contained in $\pi_1(P)$, is said to be admissible with respect to $P$ (and similarly for submodules  of $A^{k_2}$ containing $i_2^{-1}(P)$ and contained in $\pi_2(P)$).

\begin{lemma} Let $\Phi$ assign an $A$-submodule $N$ of $A^{k_2}$ to the submodule $\Phi(N) = i_1^{-1}(i_2(N) + P)$ of $A^{k_1}$. Then $\Phi$ is a bijection between the admissible submodules of $A^{k_2}$ and the admissible submodules of $A^{k_1}$, with respect to $P$.
\end{lemma}
\noindent Proof:  If $N = i_2^{-1}(P)$, then $(i_2(N) + P) = P$, and so $\Phi(i_2^{-1}(P)) = i_1^{-1}(P)$.  
If $N = \pi_2(P)$, it is an elementary calculation that $\Phi(\pi_2(P)) = \pi_c(P)$. 
As the assignment $\Phi$ is inclusion preserving, it follows that it maps an admissible submodule of $A^{k_2}$ to an admissible submodule of $A^{k_1}$, with respect to $P$.

Similarly, for an $A$-submodule $M$ of $A^{k_1}$, define $\Psi (M) = i_2^{-1}(i_1(M)+ P)$. It assigns admissible submodules of $A^{k_1}$ to admissible submodules of $A^{k_2}$, with respect to $P$. It is clear $\Phi$ and $\Psi$ are inverses of one another, hence they are both bijections.    \hspace*{\fill}$\square$\\

\begin{proposition}Let $P \subset A^k$, $\F$ an injective cogenerator and $B = \ker_\F(P)$. Then every subsystem of $\pi_2(B)$ containing $i_2^{-1}(B)$ can be achieved by a unique controller contained in $\pi_1(B)$ and containing $i_1^{-1}(B)$.
\end{proposition}
\noindent Proof: Let $B'$ be a subsystem of $\pi_2(B)$ containing $i_2^{-1}(B)$. It equals $\ker_\F(N)$, for a unique submodule $N$ of $A^{k_2}$. This submodule $N$ is admissible with respect to $P$, and by Lemma 3.2 $\Phi(N)$ equals an admissible submodule of $A^{k_1}$, say $M$.

Suppose that the laws $i_1^{-1}(P)$ of $\pi_1(B)$ are augmented to this submodule $M$. Then the laws of $B$, namely the submodule $P \subset A^k$, are augmented to the submodule $i_1(M) + P$. The projection of the resultant system, $\pi_2(\ker_\F(i_1(M) + P))$, to $\F^{k_2}$ is a system whose laws are given uniquely by $i_2^{-1}(i_1(M) + P)$, namely $\Psi(M)$ of the above lemma. As $\Psi$ is inverse to $\Phi$, $\Psi(M) = N$.

Thus $B'$ is achieved by the unique subsystem of $\pi_1(B)$ defined by the submodule $M$. \hspace*{\fill}$\square$\\

We apply these considerations to the important example of the Maxwell Equations.

Let $A$ be the ring $\C[\del_x, \del_y, \del_z, \del_t]$ of differential operators on space and time. The Maxwell equations (in Gaussian units) are 
\[\nabla \cdot E - 4\pi \rho = 0, \hspace{.3cm}\nabla \cdot B = 0, \]
\[\nabla \times E + \frac{1}{\tt{c}} \del_t B = 0, \hspace{.2cm} \nabla \times B - \frac{1}{\tt{c}}(4\pi J +  \del_t E) = 0 ~, \]
where $E, B$ are the electric and magnetic fields, $\rho, J$, the electric charge and electric current densities, and $\tt{c}$, the speed of light.
\vspace{1mm}
 
The partial differential operator defined by  these equations is \[P(\del): \F^{10} \rightarrow \F^8 ~,\] where

\[
P(\del)  = \left(
\begin{array}{lccccccccc}
 \phantom{-} \del_x & \del_y & \partial_z & 0 & 0 & 0 & -4 \pi \phantom{-} & 0 & 0 & 0 \\
\phantom{-} 0 & 0 & 0 & \del_x & \del_y & \del_z & 0 & 0 & 0 & 0 \\
\phantom{-} 0 & -\del_z \phantom{-} & \del_y & \frac{1}{\tt{c}}\del_t \phantom{-}& 0 & 0 & 0 & 0 & 0 & 0 \\
\phantom{-} \del_z & 0 & -\del_x \phantom{-} & 0 & \frac{1}{\tt{c}}\del_t \phantom{-}& 0 & 0 & 0 & 0 & 0 \\
-\del_y  & \del_x & 0 & 0 & 0 & \frac{1}{\tt{c}}\del_t \phantom{-}& 0 & 0 & 0 & 0 \\
\frac{1}{\tt{c}} \del_t  & 0 & 0 & 0 & \del_z  & -\del_y \phantom{-} & 0 & \frac{4\pi}{\tt{c}} & 0 & 0 \\
\phantom{-} 0 & \frac{1}{\tt{c}} \del_t \phantom{-}& 0 & -\del_z \phantom{-} & 0 &   \del_x  & 0 & 0 & \frac{4\pi}{\tt{c}} & 0 \\
\phantom{-} 0 & 0 & \frac{1}{\tt{c}}\del_t \phantom{-} &  \del_y & - \del_x \phantom{-} & 0 & 0 & 0 & 0 & \frac{4\pi}{\tt{c}}
\end{array}
\right),
\]
\vspace{1mm}

\noindent and where $\F$ is either $\Dr(\R^4)$ or $\Cinf(\R^4)$. The eight rows of this matrix correspond to the two equations involving divergence, and the six equations involving curl. It operates on 
$(E_1, E_2, E_3, B_1, B_2, B_3, \rho, J_1, J_2, J_3) \in \F^{10}$, the components of which are the components of $E, B, \rho$ and $J$. The electro-magnetic system is the kernel $\ker_\F(P(\del)) $ of this operator.

Let $P \subset A^{10}$, be the submodule generated by the rows of $P(\del)$. Then $A^{10}/P$ is torsion free, hence $\ker_\F(P)$ is controllable in the sense of Willems. The problem now is to control the electric and magnetic fields by suitable choices of the control variables $\rho$ and $J$. %In our notation, $(E_1, E_2, E_3, B_1, B_2, B_3) \in A^{k_2}$, and $(\rho, J_1, J_2, J_3) \in A^{k_1}$. 

The submodule $P$ determines the submodules $\pi_2(P)$ and $i_2^{-1}(P)$ of $A^6$, and the submodules $\pi_1(P), ~i_1^{-1}(P)$ of $A^4$. An elementary calculation shows that $i_1^{-1}(P)$ is the submodule of $A^4$ generated by the continuity equation $\del_t \rho + \nabla \cdot J = 0$, whereas $\pi_1(P) = A^4$. 

Similarly, the submodule $i_2^{-1}(P) \subset A^6$ is generated by the laws $\nabla \cdot B = 0$ and $\nabla \times E + \frac{1}{\tt{c}}\del_t B = 0$; these are the homogeneous Maxwell equations given by the submatrix defined by rows 2 to 5, and columns 1 to 6 of $P(\del)$. Finally, the submodule $\pi_2(P)$ is generated by the rows of the $8 \times 6$ submatrix of $P(\del)$ defined by its first 6 columns. They are the Maxwell equations in vacuum, namely the homogeneous equations together with $\nabla \cdot E = 0$ and $\nabla \times B - \frac{1}{\tt{c}} \del_t E = 0$. 
 
By the remark following Lemma 3.1, the laws governing the control variables $\rho$ and $J$ can be any $A$-submodule $M$ of $A^4$ containing the continuity equation. Then, $\rho$ and $J$ would be restricted to lie in $\ker_\F(M)$. In other words, the control variables can be restricted by any system of differential equations containing the continuity equation. %Thus, charge and current could be considered to play the classical role of `inputs'.

The restriction of $\rho$ and $J$ by the laws in $M$, results in a unique subsystem of the system $\ker_\F(i_2^{-1}(P))$ of homogeneous solutions, and containing the vacuum solutions $\ker_\F(\pi_2(P))$. This is the system determined by the admissible submodule $\Psi(M) \subset A^6$ (in the notation of Lemma 3.2). 

\vspace{1mm}
In other words, every achievable subsystem of the electro-magnetic field lies between two systems, one, the solutions of the homogeneous Maxwell equations, and the other, the solutions of the vacuum equations. They are
obtained by imposing additional differential constraints on the current and charge densities, and these constraints translate to laws that the electric and magnetic fields must satisfy, in addition to the homogeneous equations.  As every controller must satisfy the continuity equation, it folows that there is only one controller that accomplishes a given restriction.
\vspace{1.5mm}

As an example, suppose that the electric charge density $\rho$ is set to 0 (so that $\nabla \cdot J = 0$), by the imposition of the law defined by the cyclic submodule of $A^4$ generated by $(1, 0, 0, 0)$. Let $M$ be the submodule of $A^4$ generated by this law together with the continuity equation. Then $\Phi(M) \subset A^6$ is the submodule  generated by the homogeneous equations together with $\nabla \cdot E = 0$, and the system $B_2$ is restricted to $\ker_\F(\Phi(M))$ by this control action. 
\vspace{1.5mm}

Conversely, by Proposition 3.6, every electro-magnetic system contained between these two extremes is achievable by suitably restricting  electric charge and electric current, in addition to satisfying the continuity equation. 
\vspace{2mm} 
 
This is precisely the physics of the electromagnetic field.

\newpage

\section{The Controllable-Uncontrollable Decomposition}

\vspace{2mm}

We continue our study of the influence of the cancellation and characteristic ideals, $\frak{i}_\ell(P)$, $\frak{i}_k(P)$, on the behaviour of the distributed system $\ker_\F(P)$ defined by $P \subset A^k$.
\vspace{1.5mm}

We first obseve that the condition $\frak{i}_\ell \neq 0$ of Theorem 3.2 is superfluous when $n = 1, 2$:

\vspace{1mm}
\noindent $n = 1$: Now every submodule $P$ of $\mathbb{C}({[\dbydt]})^k$ is free (being finitely generated and torsion free). Moreover, an algebraic subset of $\mathbb{C}$ of dimension $n-2$ is empty. Hence, Theorem 3.2 specializes without qualification to the classical Hautus test, as well as to Proposition 3.2.

\vspace{1mm}
\noindent $n = 2$: Now  $A = \mathbb{C}[D_1, D_2]$. Let $P \subset A^k$ be a submodule such that $\ker_\F(P)$ is controllable, where $\F$ is an injective cogenerator. Then, $A^k/P$ is torsion free, hence it injects into a free module $A^{k_1}$, for some $k_1$ (by 
Proposition 3.1). Let $P(D)$ be an $\ell \times k$ matrix whose rows is a minimum set of generators for $P$. The following sequence is then a minimal resolution
\[
\rightarrow A^{\ell_1} \longrightarrow A^\ell \stackrel{P^t(D)}{\longrightarrow} A^k \longrightarrow A^{k_1} \longrightarrow A^{k_1}/(A^k/P) \rightarrow 0 ~,
\]
and hence $A^{\ell_1} = 0$ as the global dimension of the ring $A$ equals 2. The morphism $P^t(D)$ is injective, and this implies that its image $P$ is free. Thus, a controllable two dimensional system is defined by a free submodule, its cancellation ideal $\frak{i}_\ell(P)$ is thus always nonzero. Theorem 3.2(ii)  then implies the following result of Wood et al. (\cite{woo}):
\begin{corollary} Let $n = 2$, and $\F$ an injective cogenerator. Then  $\ker_\F(P)$ is controllable if and only if $P(x)$ drops rank at most at a finite number of points in $\mathbb{C}^2$ (where the $\ell$ rows of $P(D)$ is a minimum set of generators for $P$).
\end{corollary} 
\noindent Proof: By the above discussion it is unnecessary to assume that $\frak{i}_\ell (P)\neq 0$, so that by Theorem 3.2, $\ker_\F(P)$ is controllable if and only if the dimension of the cancellation variety of $\frak{i}_\ell(P)$ is equal to 0. The finite set of points of this variety is precisely where $P(x)$ drops rank.
\hspace*{\fill}$\square$\\

The above discussion together with Corollary 3.4 also implies that when $n = 2$ and $\F$ an injective cogenerator, there cannot be a phenomenon such as the controllable system $\ker_\F(\curl)$. 
\begin{corollary} Let $n = 2$, and $\F$ an injective cogenerator. Then any nonzero controllable system is strictly underdetermined.
\end{corollary}

We return to the case of general $n$. The following example explains the nomenclature  `cancellation ideal'  for $\frak{i}_\ell$:
\vspace{1mm}

\noindent Example 4.1: Consider the special case of a `scalar system' given by $p_1(D)f_1 = p_2(D)f_2$, i.e. the system defined by the kernel  of the map 
\[
\begin{array}{lccc} 
P(D) = (p_1(D), -p_2(D)): & \F^2 & \longrightarrow & \F \\ 
& (f_1,f_2) & \mapsto & (p_1(D),-p_2(D))\left [ \begin{array}{l} f_1\\f_2 \end{array}   \right ] ,
\end{array}
\]
where $\F$ is an injective cogenerator. The cancellation ideal is the ideal $\frak{i}_1(P) = (p_1,p_2)$ generated by $p_1$ and $p_2$. Theorem 3.2 asserts that this kernel is controllable if and only if there is no non-constant $p$ which divides both $p_1$ and $p_2$. This is the classical \index{pole-zero cancellation}
{\it pole-zero cancellation} criterion. 

Next we study the structure of a distributed system that is not necessarily controllable. We first observe that it contains a maximal controllable sub-system, and that the quotient is  isomorphic to an \index{uncontrollable system}uncontrollable system in the following sense.

\begin{definition} A system is uncontrollable or autonomous if none of its nonzero sub-systems is controllable. 
\end{definition}
Thus, a sub-system of an uncontrollable system is also uncontrollable.  As the zero system is (vacuously) uncontrollable, and also controllable, we use these adjectives only for nonzero systems.

We have the following characterisation.

\begin{proposition} Let $\F$ be an injective cogenerator. Then $\ker_\F(P)$ is uncontrollable if and only if $A^k/P$ is a torsion module.
\end{proposition} 
\noindent Proof: Sub-systems of $\ker_\F(P)$ are in bijective correspondence with submodules $P' \subset A^k$ that contain $P$. A nonzero sub-system, say corresponding to $P' \nsubseteq A^k$, is controllable if and only if $A^k/P'$ is torsion free. However, no submodule of $A^k/P$ is torsion free if and only if it is torsion.
\hspace*{\fill}$\square$\\

By Proposition 2.2, a controllable system is the closure of its compactly supported elements. In contrast, there is no compactly supported element in an uncontrollable system.

\begin{proposition} Let $\ker_\F(P)$ be an uncontrollable system as in the above proposition. Then there is no nonzero compactly supported element in it. 
\end{proposition}
\noindent Proof: Suppose $f$ is a nonzero compactly supported element in $\ker_\F(P)$. Let $q(D) \in A^k \setminus P$ be such that $q(D)f \neq 0$, and let $a(D)$ be a nonzero element in $A$ such that $a(D)q(D) \in P$  (such an element exists because $A^k/P$ is torsion). Then $a(D)(q(D)f) = 0$, but this contradicts Paley-Wiener as $q(D)f$ has compact support.  
\hspace*{\fill}$\square$\\

\noindent Remark 4.1: An uncontrollable system does not also contain rapidly decreasing elements in it. The proof is identical to the proof above, and it only requires us to observe that because the Fourier transform maps $\Sc$ to itself, no nonzero element of $A$ has a nonzero solution in $\Sc$.
\vspace{2mm}

Recollect that a prime ideal $p$ of a commutative ring $R$ is an \index{associated prime}associated prime of an $R$-module $M$ if it is equal to the annihilator $\ann(x)$ of some $x \in M$; the set of its associated primes is denoted $\ass(M)$. An element $r$ of $R$ is a \index{zero divisor}zero divisor on $M$ if there is a nonzero $x$ in $M$ such that $rx = 0$. The maximal elements of the family $\{\ann(x) ~| ~0 \neq x \in \M\}$ are associated primes of $M$, hence the union of all the associated primes is the set of all the zero divisors on $M$ (for instance \cite{ma}). It then follows that $M$ is torsion free if and only if 0 is its only associated prime, and is torsion if and only if 0 is not an associated prime.
\vspace{1.5mm}

Now let $P$ be an $A$-submodule of $A^k$. Let $P_0$ be the submodule $\{x \in A^k ~| ~\exists ~a \neq 0 ~with ~ax \in P\}$. Then $P_0$ contains $P$, and the quotient $P_0/P$ is the submodule of $A^k/P$ consisting of its torsion elements. The following sequence 
\begin{equation}
0 \rightarrow P_0/P \longrightarrow A^k/P \longrightarrow A^k/P_0 \rightarrow 0
\end{equation}
is exact, where $A^k/P_0$ is torsion free. In general, given a short exact sequence of $A$-modules, the associated primes of the middle term is contained in the union of the associated primes of the other two modules. However, here it is clear that  we have equality.

\begin{lemma} $\ass (A^k/P) = \ass (P_0/P) \sqcup \ass (A^k/P_0)$ (disjoint union). Hence, if $P \subsetneq P_0 \subsetneq A^k$, then $\ass (A^k/P_0) = \{0\}$ and $\ass (P_0/P)$ is the set of all the nonzero associated primes of $A^k/P$.
\hspace*{\fill}$\square$
\end{lemma}

Applying the functor $\homo_A(-, ~\F)$, where $\F$ is an injective cogenerator, to the above sequence gives the exact sequence
\[
0 \rightarrow {\ker}_\F(P_0) \longrightarrow {\ker}_\F(P) \longrightarrow {\homo}_A (P_0/P,~\F) \rightarrow 0 
\]
As $A^k/P_0$ is torsion free, $\ker_\F(P_0)$ is a controllable sub-system of $\ker_\F(P)$. If $P_1$ is any $A$-submodule of $A^k$ such that $P \subset P_1 \nsupseteq P_0$, then $A^k/P_1$  has torsion elements and $\ker_\F(P_1)$ is not controllable. Hence $\ker_\F(P_0)$ is the largest controllable sub-system of $\ker_\F(P)$  in the sense that any other controllable sub-system is contained in it. We call it the controllable part of the system.

Suppose $P_0/P$ can be generated by $r$ elements; then $P_0/P \simeq A^r/R$, for some submodule $R \subset A^r$, hence $\homo_A(P_0/P, ~\F) \simeq \ker_\F(R)$.  As $A^r/R$ is a torsion module, $\ker_\F(R)$ is uncontrollable. This system is a quotient, and not a sub-system of $\ker_\F(P)$ (unless the above short exact sequence splits). We now show that there is a sub-system  of $\ker_\F(P)$, not canonically determined, that is  isomorphic to this quotient. 
\vspace{2mm}

Recollect that a submodule $N$ of an $R$-module $M$ is a \index{primary submodule}primary submodule of $M$ if every zero divisor $r \in R$ on $M/N$
satisfies $r \in \sqrt{\ann(M/N)}$. It then follows that $\ass(M/N)$ contains only one elment, and if $\ass(M/N) = \{p\}$, then $N$ is a $p$-primary submodule of $M$. We also note that as $A$ is a noetherian domain, $N$ is 0-primary if and only if $M/N$ is torsion free.
\vspace{1.5mm}

The result we need is the following (for instance \cite{ma}).
\vspace{1mm}

\noindent Theorem: Let $R$ be a noetherian ring, and $M$ a finitely generated $R$-module.  Then every proper submodule $N$ of $M$ admits an \index{irredundant primary decomposition}irredundant primary decompostion, which is to say that $N = N_1 \cap  \cdots \cap N_r$, where $N_i$ is a $p_i$-primary submodule of $M$, $p_i \neq p_j$ for $i \neq j$, and where none of the  $N_i$ can be omitted (i.e. $N \neq N_1 \cap \cdots \cap N_{i-1} \cap N_{i+1} \cap \cdots \cap N_r$, for any $i$).  $\ass(M/N)$ then equals $\{p_1, \ldots, p_r\}$. If $p_i$ is a minimal element in $\ass(M/N)$, then the $p_i$-primary component $N_i$ is uniquely determined.
\vspace{2mm}

We first improve the second part of Lemma 2.1(ii) when $\F$ is injective, and when the collection of submodules of $A^k$ is finite.

\begin{lemma} Let $\F$ be an injective $A$-module, then $\sum _{i=1}^r \ker_\F(P_i) = \ker_\F(\bigcap_{i=1}^r P_i)$.
\end{lemma}
\noindent Proof: Consider the following exact `Mayer-Vietoris' sequence
\[
0 \rightarrow A^k/(P_1\cap P_2)
\stackrel{i}{\longrightarrow}
A^k/P_1 \oplus A^k/P_2
\stackrel{\pi}{\longrightarrow} 
A^k/(P_1+P_2)
\rightarrow 0 ~,
\]
where $i([x]) = ([x],-[x])$, $\pi ([x],[y]) = [x+y]$ and where 
$[x]$ indicates the class of $x$ in various quotients. As $\F$ is an injective $A$-module, it follows that
\[
0 \rightarrow 
{\homo}_A(A^k/(P_1+P_2), ~\F)
\longrightarrow
{\homo}_A(A^k/P_1, ~\F) \oplus
{\homo}_A(A^k/P_2, ~\F)
\longrightarrow 
\]\[
{\homo}_A(A^k/(P_1\cap P_2), ~\F) \rightarrow 0 
\] is also exact. But ${\homo}_A(A^k/(P_1+P_2), ~\F) \simeq \ker_\F(P_1+P_2) = \ker_\F(P_1) \cap \ker_\F(P_2)$ by the first part of Lemma 2.1(ii), hence the above exact sequence implies that $\ker_\F(P_1)+\ker_\F(P_2) = \ker_\F(P_1\cap P_2)$.

The lemma now follows by induction. \hspace*{\fill}$\square$\\

\begin{theorem} Let $\F$ be an injective cogenerator, and $P$ a proper submodule of $A^k$. Then $\ker_\F(P)$ can be written as the sum of its controllable part and an uncontrollable sub-system.
\end{theorem}
\noindent Proof: If $A^k/P$ is either torsion free or a torsion module, then there is nothing to be done. Otherwise, let $P = P_0 \cap P_1 \cap \cdots \cap P_r$, $r \geqslant 1$, be an irredundant primary decomposition of $P$ in $A^k$, where $P_0$ is a 0-primary submodule of $A^k$, and the $P_i$,  $i \geqslant 1$, are $p_i$-primary submodules for nonzero primes $p_i$.  Then $A^k/P_0$ is torsion free, and $A^k/(P_1 \cap \cdots \cap P_r)$ is torsion; hence $\ker_\F(P) = \ker_\F(P_0 \cap P_1 \cap \cdots \cap P_r) = \ker_\F(P_0) + \ker_\F(P_1 \cap \cdots \cap P_r)$, where $\ker_\F(P_0)$ is controllable, and $\ker_\F(P_1 \cap \cdots \cap P_r) = \ker_\F(P_1) + \cdots + \ker_\F(P_r)$ is uncontrollable. As 0 is a minimal associated prime of $A^k/P$, $P_0$ is uniquely determined, and in fact equals the $P_0$ of Lemma 4.1. Hence, $\ker_\F(P_0)$ is the controllable part of $\ker_\F(P)$.
\hspace*{\fill}$\square$\\

Propositions 2.2 and 4.2 imply the following corollary.

\begin{corollary} Let $\F$ be either $\Dr$ or $\Cinf$ with their standard topologies, and $P$ an $A$-submodule of $A^k$. Then $\ker_\F(P)$ can be written as the sum of two sub-systems: one, which is the closure of the set of  compactly suported elements, and another, which contains no nonzero compactly supported, nor rapidly decreasing, element.
\end{corollary}

We now describe the nonzero associated primes of $A^k/P$ that determine the uncontrollable part of $\ker_\F(P)$ in the above theorem, in the case when $P$ is a free proper submodule of $A^k$.
 
\begin{proposition} Let $P$  be a free submodule of $A^k$. Then the nonzero associated primes of $A^k/P$ are the principal ideals generated by the irreducible $p(D)$ that divide every element of $\frak{i}_\ell(P)$.
\end{proposition}
Proof: By Proposition 3.4, the cancellation ideal $\frak{i}_\ell(P)$ of $P$ is nonzero. The first half of Theorem 3.2 demonstrates that an irreducible $p(D)$ which divides every generator of $\frak{i}_\ell(P)$, and hence every one of its elements, is a zero divisor on $A^k/P$. Conversely, suppose that $a(D)$ is a zero divisor on $A^k/P$, and that $x$ is an element of $A^k \setminus P$ such that $a(D)x \in P$. Because $A$ is a UFD, $a(D)$ is a unique product of irreducible factors, say $p_1(D) \cdots p_s(D)$ (where the $p_i$ are not necessarily distinct). Then there is an $i, ~1\leqslant i \leqslant r$, such that $p_1(D) \cdots p_i(D)x \in P$, but $p_1(D) \cdots p_{i-1}(D)x \not\in P$.
Thus it follows that the maximal elements of the family of principal ideals generated by zero divisors on $A^k/P$ are the principal ideals generated by irreducible zero divisors. 
Furthermore, by the second half of the proof of Theorem 3.2, every irreducible zero divisor divides every generator of $\frak{i}_\ell(P)$.  Together with the above observation, this shows that the set of irreducible zero divisors on $A^k/P$ is precisely the set of irreducible common factors of the $k \choose \ell$ generators of the cancellation ideal $\frak{i}_\ell(P)$. As the latter set is finite, it follows that the set of irreducible zero divisors on $A^k/P$ is a finite set.

This in turn implies that the nonzero associated primes of $A^k/P$ are precisely the principal ideals generated by the irreducible zero divisors described above. For suppose that some associated prime $p$ is not principal. Let $p_1(D)$ be an irreducible element of $p$  such that its degree is minimum amongst all the (nonzero) elements in $p$. Let $p_2(D)$ be another irreducible element in $p$, not a multiple of $p_1(D)$, again of minimum degree amongst all elements in $p \setminus p_1(D)$. Then $p_1(D) + \alpha p_2(D)$ is irreducible, for all 
$\alpha \in \mathbb{C}$. These infinite number of ireducible elements are all in $p$, and hence are all zero divisors on $A^k/P$, contradicting the assertion above. 
\hspace*{\fill}$\square$\\

By Theorem 3.2(ii), the system defined by a free submodule of $A^k$ is not controllable if and only if its cancellation variety is of codimension 1. In this case, the above proposition, together with the proof of Theorem 4.1, implies the following corollary.

\begin{corollary} Let $P$ be a free submodule of $A^k$, and $\F$ an injective cogenerator. Then, the codimension 1 irreducible components of the cancellation variety $\mathcal{V}(\frak{i}_\ell(P))$ determine the uncontrollable part of $\ker_\F(P)$ in Theorem 4.1.
\end{corollary}

These results, which follow from Theorem 3.2, are valid only when the cancellation ideal is nonzero, which is a specific (though generic) case of underdetermined systems. They are not valid for underdetermined systems whose cancellation ideal equals 0.
\vspace{4mm}

We now study a distributed system from the point of view of its characteristic ideal.

For any submodule $P \subset A^k$, its characteristic ideal $\frak{i}_k(P)$  and characteristic variety $\mathcal{V}(\frak{i}_k(P)) \subset \mathbb{C}^n$ are fundamental objects that determine $\ker_\F(P)$. Indeed, the Integral Representation Theorem of Palamadov \cite{p}, a vast generalization of Malgrange's Theorem quoted in Remark 2.4, states that an element in $\ker_{\Dr}(P)$ is an absolutely convergent integral of its exponential solutions, and these exponential solutions are determined by the points on its characteristic variety. 
\vspace{1mm}

Our study is based on the following elementary lemma.

\begin{lemma} 
$\frak{i}_k(P) \subset \ann(A^k/P) \subset \sqrt{\frak{i}_k(P)}$, and thus $\V(\frak{i}_k(P)) = \V(\ann(A^k/P))$ .
\end{lemma}
Proof: Let $d$ be the determinant of a $k\times k$ matrix $M$ all whose rows belong to the submodule $P \subset A^k$. 
Let $M'$ be the matrix adjoint to $M$, so that $M'M$ is the $k\times k$ diagonal matrix whose diagonal elements are all equal to $d$. But the rows of $M'M$ are also in $P$, and as these rows are 
$d e_1, \ldots , d e_k$  ($\{e_i, ~1 \leqslant i \leqslant k\}$ is the standard basis for $A^k$),  
it follows that $d$ is in $\ann(A^k/P)$. This implies that $\frak{i}_k(P) \subset \ann(A^k/P)$.

Conversely, let $a$ be in $\ann(A^k/P)$ so that $a e_1,\ldots ,a e_k$ are all in $P$. The determinant 
of the $k\times k$ matrix whose rows are these $a e_i$ equals $a^k$. Thus $a^k$ is in $\frak{i}_k(P)$ which implies that $a$ is in 
$\sqrt{\frak{i}_k(P)}$.
\hspace*{\fill}$\square$ 
\vspace{2mm}

Thus, $\ker_\F(P)$ is overdetermined  if  $\ann(A^k/P) \neq 0$. 

\begin{corollary} Let $\F$ be an injective cogenerator, and $P$ a proper submodule of $A^k$. Then $\ker_\F(P)$ is uncontrollable if and only if $\frak{i}_k(P) \neq 0$.  
\end{corollary}
\noindent Proof: By the above lemma, $\frak{i}_k(P) \neq 0$ if and only if $\ann(A^k/P) \neq 0$. As $A^k/P$ is finitely generated, $\ann(A^k/P) \neq 0$ is equivalent to saying that $A^k/P$ is a torsion module. The corollary now follows from Proposition 4.1. 
\hspace*{\fill}$\square$ 
\vspace{3mm}

\begin{corollary} Suppose $P$ is a $p$-primary submodule of $A^k$, $\F$ an injective cogenerator. Then $\ker_\F(P)$ is uncontrollable if $p$ is a nonzero prime, and controllable otherwise. 
\end{corollary}
\noindent Proof: Suppose $p$ is nonzero. As $p$ equals the set of zero divisors on $A^k/P$, it follows that $\sqrt{\ann(A^k/P)} = p$. Thus $\ann(A^k/P)$ is nonzero, and by the above lemma, so is $\frak{i}_k(P)$. 

Conversely, $P$ is 0-primary if and only if $A^k/P$ is torsion free. Then, and only then, is $\ker_\F(P)$ controllable. \hspace*{\fill}$\square$ \\

\begin{corollary} The controllable part of a strictly underdetermined system is nonzero ($\F$ as above).
\end{corollary}
\noindent Proof: The characteristic ideal of a strictly underdetermined system (i.e. $\ell < k$) is equal to 0. \hspace*{\fill}$\square$
\vspace{2mm}

The condition $\ann(A^k/P) = 0$ is weaker than the condition guaranteeing controllability of $\ker_\F(P)$, namely that $A^k/P$ be torsion free. We now explain one implication of this condition on  the dynamics of the system.

 \begin{proposition} Let $P$ be a submodule of $A^k$, and $\F$ an injective $A$-module.  Let $\pi_j: \F^k \rightarrow \F, ~(f_1, \ldots, f_k) \mapsto f_j$ be the projection onto the $j$-th coordinate. Then $\ann(A^k/P) = 0$ if and only if there is a $j, ~1\leqslant j \leqslant k$, such that the restriction $\pi_j: \ker_\F(P) \rightarrow \F$ is surjective.
\end{proposition}

\noindent Proof: Let $i_j: A \rightarrow A^k, ~a \mapsto (0, \ldots, a, \ldots, 0)$ be the inclusion into the $j$-th coordinate.  Now, $\ann(A^k/P) = 0$ is equivalent to saying that there is a $j, ~1\leqslant j \leqslant k$, such that $i_j^{-1}(P) = 0$. By Proposition 3.5, $\pi_j(\ker_\F(P)) = \ker_\F(i_j^{-1}(P))$, as $\F$ is injective. This is in turn equivalent to $\pi_j(\ker_\F(P)) = \F$ by Corollary 3.1.
\hspace*{\fill}$\square$\\

As the only controllable sub-systems of $\F$ are 0 and $\F$ itself, the surjectivity statement above is termed \index{coordinate controllability}{\it coordinate controllability}. Thus a system $\ker_\F(P)$ is coordinate controllable, where $\F$ is an injective $A$-module, if and only if $\ann(A^k/P) = 0$. Lemma 4.3 then implies the following supplement to the previous corollary.

\begin{corollary} A strictly underdetermined system is coordinate controllable.
\hspace*{\fill}$\square$
\end{corollary}

To formulate Willems' behavioural controllability in these terms, consider a general homothety: given $r(D) = (a_1(D), \ldots ,a_k(D))$ in $A^k$, let $i_r: A \rightarrow A^k$ be the morphism $a(D) \mapsto a(D)r(D)$ (which is an injection when $r(D) \neq 0$). Applying $\homo_A(-, ~\F)$ to this morphism gives
$r(D): \F^k \rightarrow \F$, mapping $f=(f_1,\ldots ,f_k)$ to $a_1(D)f_1 + \cdots + a_k(D)f_k$ (and which is a surjection when $r(D) \neq 0$); thus, the injections $i_j$ and the projections $\pi_j$ in the above proposition correspond to $r(D) = (0, \ldots, 1, \ldots, 0)$. Behavioural controllability is then an assertion about the restrictions of the maps $r(D)$ above to $r(D):\ker_\F(P) \rightarrow \F$, for all $r(D)$ in $A^k$.

\begin{proposition} Let $P$ be a submodule of $A^k$, and $\F$ an injective cogenerator. The system $\ker_\F(P)$ is controllable if and only if every image $r(D)(\ker_\F(P))$ is controllable (in other words if and only if every $r(D)(\ker_\F(P))$ is either 0 or all of $\F$).
\end{proposition}
\noindent Proof: Consider the map 
\[\begin{array}{lccc}
m_r: & A & \longrightarrow & A^k/P  \\
 &1 & \mapsto & [r(D)] 
 \end{array}
\]
defined by $r(D) \in A^k$. Applying $\homo_A(-, ~\F)$ to it gives the restriction $r(D): \ker_\F(P) \rightarrow \F$  described above. The map $m_r$ is the 0 map if $r(D) \in P$, hence the image $r(D)(\ker_\F(P))$ equals 0, and is controllable. 

Suppose now that $\ker_\F(P)$ is controllable so that $A^k/P$ is torsion free. Then $m_r$ is an injective map for every $r(D) \notin P$, hence $r(D)(\ker_\F(P))$ equals $\F$ (as $\F$ is an injective module). Conversely, suppose $\ker_\F(P)$ is not controllable. Let $r(D) \notin P$ be a torsion element in $A^k/P$; then $m_r$ is not injective, neither is it the 0 map. Hence $r(D)(\ker_\F(P))$ is not all of $\F$, nor does it equal 0. Thus it is not controllable by Corollary 3.2. 
\hspace*{\fill}$\square$
\vspace{3mm}

 The above discussion raises the question: can the projection of a controllable system to several coordinates also be surjective? The example below shows that in general this need not be the case. 
\vspace{3mm}

\noindent Example 4.2: Let $A = \mathbb{C}[D_1,D_2]$. The system given by the kernel of the map $(D_1, D_2):\F^2 \rightarrow \F$ is controllable as its cancellation variety, the origin in $\mathbb{C}^2$, is of codimension 2. By Proposition 4.4, each $\pi_j:\ker_\F((D_1,D_2))  \rightarrow \F, ~j = 1,2$, ~is surjective, but $\ker_\F((D_1,D_2)) \rightarrow \F^2$ is clearly not. 
\vspace{3mm}

A strongly controllable system does however surject onto several of its coordinates (Remark 3.2). Indeed, if the rank of the free module $A^k/P$ equals $k'$, then $\ker_\F(P) \simeq \F^{k'}$.

\vspace{2mm}
We now have identified three graded notions of controllability of  the system $\ker_\F(P)$, 
each corresponding to increasingly weaker requirements on $A^k/P$. %The relationship of these notions to each other is best seen in terms of projections of the system to various coordinates. 

\vspace{1.5mm}
\FloatBarrier
\begin{table}[h]
\caption{Notions of controllability of $\ker_\F(P)$}
\begin{center}
\begin{tabular}{cc}   Strong controllability \phantom{X} ~$\iff  A^k/P$ is free \\
\phantom{XX-} Behavioural controllability  $\iff  A^k/P$ is torsion free\\
\phantom{-}Coordinate controllability ~$\iff \ann (A^k/P) = 0$
\end{tabular}
\end{center}
\vspace{2mm}
(In the above table, the first equivalence is valid for an arbitrary $A$-module $\F$, the second for $\F$ that are injective cogenerators, and the last for injective $\F$.)
\end{table}
\FloatBarrier

A notion of `inputs' in \cite{w} is a set of cordinates $\{i_1, \cdots i_{k'}\}$ such that the projection $\pi: \ker_\F(P) \rightarrow \F^{k'}$ to these coordinates is surjective.
It reflects an intuitive understanding of an input as a signal which is not restricted in any way, and which can therefore be any element in $\F$ (electric charge $\rho$ and current density $J$ in the example of Maxwell equations in the previous section do not satisfy this condition). 
By Proposition 4.4, uncontrollable systems are precisely those which do not admit inputs (they are called autonomous systems in \cite{ps,w}). 
\vspace{1mm}

We conclude this section with a few additional comments on distributed systems whose behaviours resemble overdetermined lumped systems (\cite{ps,lz}).

\begin{definition} A system $\ker_\F(P)$ is \index{strongly autonomous system}strongly autonomous (or strongly uncontrollable) if $\ker_\F(P)$ is a finite dimensional $\C$-vector space.
\end{definition}

We first observe that if $m$ is a maximal ideal of $A = \mathbb{C}[D_1, \ldots ,D_n]$, say $m = (D_1 - \xi_1, \ldots, D_n - \xi_n)$, then $\ker_{\Dr}(m)$ is the 1-dimensional $\mathbb{C}$-subspace of $\Dr$ spanned by $e^{\imath <\xi, x>}$, where $\xi = (\xi_1, \ldots , \xi_n)$. As $\ker_\F(m) = \ker_{\Dr}(m) \cap \F$, its dimension equals 1 if $e^{\imath <\xi, x>}$ is in $\F$, and 0 otherwise. 

\begin{lemma}  Let $j = (p_1(D) , \ldots , p_i(D) , \ldots , p_r(D) )$ be an ideal of $A$ such that 
$\ker_\F(j)$ is a finite dimensional $\mathbb{C}$-vector space. Let $j' = (p_1(D) , \ldots , p_i^2(D) , \ldots , p_r(D))$. Then $\ker_\F(j')$ 
is also finite dimensional.
\end{lemma}
\noindent Proof:  Let $j_i$ be the ideal $(p_1(D) , \ldots , \widehat{p_i(D)} , \ldots , p_r(D))$, where $\widehat{p_i(D)}$ means that $p_i(D)$ has been
omitted. Consider the following $\mathbb{C}$-linear map:
\[\begin{array}{lccc}
P_i: & \ker_\F(j_i) & \longrightarrow & \F  \\
 & f & \mapsto & p_i(D)f 
 \end{array}
\]
Then $\ker_\F(j)$ equals the kernel of $P_i$, which is finite dimensional by assumption. Hence $\ker_\F(j') = P_i^{-1}(\ker_\F(j))$ is also finite dimensional.
\hspace*{\fill}$\square$\\

\begin{corollary} Let $j$ be an ideal of $\mathbb{C}[D_1, \ldots ,D_n]$. Then $\ker_\F(j)$ is a finite dimensional $\mathbb{C}$-vector space if and only if $\ker_\F(\sqrt{j})$ is.
\end{corollary}
\noindent Proof:  If $\ker_\F(\sqrt{j})$ is finite dimensional, then so is the system defined by every power of $\sqrt{j}$ by the above lemma. As some power of $\sqrt{j}$ is contained in $j$, $\ker_\F(j)$ is also finite dimensional. The other implication is obvious.
\hspace*{\fill}$\square$\\

\begin{proposition} Let $\F$ be either $\Dr$ or $\Cinf$, and $P$ a submodule of $A^k$. Then $\ker_\F(P)$  is strongly 
autonomous if and only if its characteristic variety has dimension 0.
\end{proposition}
\noindent Proof: Suppose that the characteristic variety $\mathcal{V}(\frak{i}_k(P))$ of $\ker_\F(P)$ is a finite set of points. Then the ideal
of this variety, which is the radical of
of $\frak{i}_k(P)$, is a finite intersection of maximal ideals.
The exponential solutions of $\sqrt{\frak{i}_k(P)}$ are the exponential solutions
corresponding to each of these maximal ideals, and these span a finite dimensional
$\mathbb{C}$-vector space, namely
$\ker_\F(\sqrt{\frak{i}_k(P)})$ (by the theorem of Malgrange quoted in Remark 2.4). By the above corollary, $\ker_\F(\frak{i}_k(P))$ is also finite dimensional. 

Suppose that $\frak{i}_k(P) = (a_1(x), \ldots ,a_r(x))$. Then every component of an element $f = (f_1 , \ldots , f_k)$
in $\ker_\F(P)$ is also a homogenous solution of $a_i(D), i = 1, \ldots ,r$, and hence belongs to $\ker_\F(\frak{i}_k(P))$.
It follows that $\ker_\F(P)$ is finite dimensional as well.

Conversely, suppose that the dimension of $\mathcal{V}(\frak{i}_k(P))$ is bigger than or equal to 1. Let $P(D)$ be any matrix whose rows 
generate the submodule $P$. At each point $\xi$ in $\mathcal{V}(\frak{i}_k(P))$, the columns $c_1, \ldots ,c_k$ of $P(x)$ 
are $\mathbb{C}$-linearly dependent, say $\sum\alpha_i c_i = 0$. Then $f = (\alpha_1e^{\imath <\xi, x>}, \ldots ,\alpha_ke^{\imath <\xi, x>})$ 
is in $\ker_\F(P)$.  These infinitely many exponential solutions, corresponding to the points in the characteristic variety, are
linearly independent. 
\hspace*{\fill}$\square$\\

It follows from the above proof that every solution of a strongly autonomous system is entire.
Thus no nonzero element  in it can vanish on any open subset of $\R^n$. 

We compare this with the behaviour of a lumped system which is either finite dimensional, or otherwise contains elements that are not entire. 
In the case of distributed systems, a behaviour could consist entirely of  entire functions, yet be
infinite dimensional (as in the case of an elliptic operator). 

\newpage

\section{Genericity}

\vspace{2mm}

We now ask how large is the class of controllable systems in the set of all systems? More precisely, we topologise the set of systems, and ask whether the subset of controllable systems  is open dense. The results below are with respect to a very coarse topology,  the \index{Zariski topology}Zariski topology. 
\vspace{1.5mm}

This question is motivated by the following considerations. A system's behaviour is modelled by differential equations which involve parameters that are only approximately known, for they are determined by measurements only upto a certain accuracy. Hence, the conclusions and predictions of a theory must be {\it stable} with respect to perturbations of the differential equations appearing in the model in order for it to be an effective theory. 

These ideas go back to the notion of structural stability of autonomous systems due to Andronov and Pontryagin \cite{ap}. 
In contrast, control systems are usually underdetermined, they admit inputs, and the system trajectories are complicated. However, our needs are modest compared to the foundational work in \cite{ap}, where a topological congugacy is sought to be established between the flow of a vector field and the flow of the perturbed vector field.  Instead, all we ask for here is the persistence of controllability under perturbations. In other words, if a system is controllable, is a perturbation of it also controllable? As we work with linear phenomena, we are able to answer this question with respect to perturbations given by a very coarse topology, even though the dynamics of the system are described by partial differential equations.
\vspace{1mm}

In this section, we restrict our attention to the case when the $A$-module $\F$ is an injective cogenerator. By Proposition 2.6, the set of distributed systems in $\F^k$ is in (inclusion reversing) bijection with the set of $A$-submodules of $A^k$. As we wish to study the variation of a system's properties with respect to variations in the differential equations that describe it, we need to topologise the set of all submodules of $A^k$, and towards this we first topologise the set of all matrices with $k$ columns  and entries in $A$.
\vspace{1mm}

Let \index{$A(d)$}$A(d), ~ d\geqslant 0$, be the subset of $A = \mathbb{C}[D_1, \ldots ,D_n]$ consisting of differential operators of degree at most $d$. There are $N(d):={n+d \choose n}$ monomials of degree at most $d$ in the partial derivatives $D_1, \ldots, D_n$, hence we identify $A(d)$ with the affine space $\mathbb{C}^{N(d)}$ equipped with the Zariski topology. There is a natural inclusion of $A(d)$ into $A(d+1)$ as a Zariski closed subset. Thus $\{A(d), d \geqslant 0\}$ is a directed system, and its (strict) direct limit, $\varinjlim{A(d)}$, is the topological \index{$\mathcal{A}$, space of differential operators}space of all differential operators. We denote it by  $\A$.

An element of the ring $A$ is a sum of monomial terms $D ^d := D_1^{d_1} \cdots D_n^{d_n}$ with complex coefficients; in the topological space $\A$ it corresponds to the point whose `coordinates' are these coefficients (if a monomial does not appear in the sum, then its coefficient is 0). Thus, there is a `coordinate axis' in the space $\A$ corresponding to each monomial in the ring $A$. We denote the indeterminate, and the coordinate axis, corresponding to the monomial $D ^d$ by $X_d$ or $X_{ d_1 \cdots d_n}$. The points of the space $\A$ corresponding to the units in the ring $A$ are the points on the axis $X_{0 \cdots 0}$ except for the origin. Its closure is this axis, and is a proper Zariski closed subset of  $\A$.

The  space $\A$ is not Noetherian, for instance the descending sequence of Zariski closed subspaces $\A \supsetneq \{X_{1 \cdots 1} = 0\} \supsetneq \{X_{1 \cdots 1}= 0, X_{2 \cdots 2} = 0\} \supsetneq ...$ ~does not stabilise. If $\{X_d\}$ denotes the union of all the indeterminates $X_{d_1 \cdots d_n}$, then the coordinate ring of the space $\A$ is $\mathbb{C}[\{X_d\}]$. 
\vspace{1mm}

Now let $k$ be a fixed positive integer. Let \index{$M_{\ell,k}$}$M_{\ell,k}$ be the set of differential operators defined by matrices with $k$ columns and $\ell$ rows, and with entries from the ring $A$. Let $M_{\ell,k}(d)$ be the set of those matrices in $M_{\ell,k}$ whose entries are all bounded in degree by $d$. We identify it with the affine space $\mathbb{C}^{\ell k N(d)}$ equipped with the Zariski topology, and just as above, we have the directed system 
$\{M_{\ell,k} (d), ~d\geqslant 0\}$, where $M_{\ell,k} (d) \hookrightarrow M_{\ell,k} (d+1)$ continuously as a Zariski closed subset. The direct limit, $\varinjlim{M_{\ell,k}(d)}$, is  $M_{\ell,k}$ with the Zariski topology, and we denote it by \index{$\mathcal{M}_{\ell,k}$}$\Mellk$ (the space $\A$, in this notation, is $\mathcal{M}_{1,1}$).

Let $S_{\ell,k}$ be the set of submodules of $A^k$ that can be generated by $\ell$ elements, and let the set of distributed systems they define in $\F^k$ be denoted by $B_{\ell,k}$. Define $\Pi_{\ell,k}:M_{\ell,k} \rightarrow S_{\ell,k}$ by mapping a matrix to the submodule generated by its rows, and equip $S_{\ell,k}$ with the quotient topology so that it is continuous. Denote this topological space by \index{$\mathcal{S}_{\ell,k}$}$\Sellk$.  We carry over this topology to $B_{\ell,k}$ via the bijection $P \mapsto \ker_\F(P)$, and denote this space by \index{$\mathcal{B}_{\ell,k}$}$\Bellk$. 

For $k = 1$, we denote $\mathcal{M}_{\ell,1}$ by $\A^\ell$, so that a point in it is a column with $\ell$ entries from $A$. We denote by $\mathcal{I}_\ell$ the space $\mathcal{S}_{\ell,1}$ of ideals of the ring $A$ that can be generated by  $\ell$ elements. In this notation, $\Pi_{\ell,1}: \A^\ell \rightarrow \mathcal{I}_\ell$ maps an element of $\A^\ell$ to the ideal generated by its $\ell$ entries.
\vspace{1.5mm}

A subset of differential operators in $\Mellk$, or of submodules in $\Sellk$ or systems in $\Bellk$, is said to be \index{generic}{\it generic} if it contains an open dense subset in this topology.

\vspace{1mm}
For each $d \geqslant 0$, $M_{\ell,k}(d)$ is irreducible (being isomorphic to affine space), hence so is the direct limit $\Mellk$ irreducible. It follows that $\Sellk$ and $\Bellk$ are also irreducible. Thus every nonempty open subset in any of these spaces is also dense. Nonetheless, we continue to use the phrase `open dense' for emphasis.
\vspace{1.5mm}

We collect a few elementary properties of these spaces.
\begin{lemma} The element $0 \in \mathcal{I}_r$, namely the $0$ ideal of $A$, is closed in $\mathcal{I}_r$, for every $r > 0$.
\end{lemma}
\noindent Proof: In the above notation, $\Pi^{-1}_{r,1}(0)$ equals $0 \in \A^r$. As $\{0\}$ is closed in $A^r(d)$ for every $d$, it is closed in $\A^r$. It follows that $\{0\}$ is a closed point of $\mathcal{I}_r$.  \hspace*{\fill}$\square$
\begin{lemma} Let $I = \{i_1, \ldots , i_{\ell'}\}$ be a set of $\ell'$ indices between 1 and $\ell$, and $J = \{j_1, \ldots ,  j_{k'}\}$ a set of $k'$ indices between 1 and $k$. Let $s: \mathcal{M}_{\ell,k} \rightarrow \mathcal{M}_{\ell',k'}$ map an $\ell \times k$
matrix with entries in $A$ to its $\ell' \times k'$ submatrix determined by the indices $I$ and $J$. Then $s$ is continuous and open.  
\end{lemma}
\noindent Proof: The map $s(d): M_{\ell,k}(d) \rightarrow M_{\ell',k'}(d)$ induced by $s$, is continuous and open for every $d$ as it is a projection $\mathbb{C}^{\ell kN(d)} \rightarrow \mathbb{C}^{\ell' k'N(d)}$ of an affine space onto the affine space given by a subset of its coordinates. Hence it follows that $s$ is continuous and open.
\hspace*{\fill}$\square$
\begin{lemma} The map $\det: \mathcal{M}_{r,r} \rightarrow \A$, mapping a square matrix of size $r$ with entries in $A$ to its determinant, is continuous. Hence, the subset of matrices in $\mathcal{M}_{r,r}$ whose determinant is nonzero, is open dense.
\end{lemma}
\noindent Proof: For each $d \geqslant 0$, the map $\det$ restricts to a map $M_{r,r}(d) \rightarrow A(rd)$ of affine spaces, which is continuous with respect to the Zariski topology as it is given by algebraic operations. Thus the set of matrices with nonzero determinant is open dense in $ M_{r,r}(d)$ for each $d$, and hence in $\mathcal{M}_{r,r}$.
\hspace*{\fill}$\square$
\vspace{1.5mm}

Let $\ell \leqslant k$. Define $\mathfrak{m}_{\ell,k}: \Mellk \rightarrow \A^{k \choose \ell}$ by mapping a matrix to its $\ell \times \ell$ minors (written as a column in some fixed order). This map is continuous by the above lemmas. If the rows of two matrices in $\Mellk$ generate the same submodule of $A^k$, then their images in $\mathcal{I}_{{k \choose \ell}}$, under the composition $\Pi_{{k \choose \ell},1}\circ \frak{m}_{\ell,k}$, are equal. Thus the map $\frak{m}_{\ell,k}$ descends to a map $\frak{i}_{\ell,k}: \Sellk \rightarrow \mathcal{I}_{{k \choose \ell}}$, mapping a submodule $P$ to the ideal generated by the maximal minors of any matrix in $\Mellk$ whose rows generate $P$ . In other words, the following diagram commutes:
\begin{equation}
\begin{tikzcd}
\mathcal{M}_{\ell,k} \arrow{r}{\Pi_{\ell,k}} \arrow[swap]{d}{\mathfrak{m}_{\ell,k}} & \mathcal{S}_{\ell,k} \arrow{d}{\mathfrak{i}_{\ell,k}} \\
\phantom{xx} \A^{k \choose \ell}
\phantom{xx} \arrow{r}{\Pi_{{k \choose \ell},1}} & ~\phantom{x} \mathcal{I}_{k \choose \ell}
\end{tikzcd}
\end{equation}

\noindent Remark 5.1: If $P \in \Sellk$ can be generated by fewer than $\ell$ elements, then clearly  $\frak{i}_{\ell,k}(P) = 0$. On the other hand, if the minimum number of generators for $P$ is $\ell$,  then $\frak{i}_{\ell,k}(P)$ is the cancellation ideal $\mathfrak{i}_\ell(P)$ of $P$.
\begin{lemma} 
The map $\mathfrak{i}_{\ell,k} : \Sellk \rightarrow \mathcal{I}_{k \choose \ell}$ is continuous.
\end{lemma}
\noindent Proof: The maps $\mathfrak{m}_{\ell,k}$, $\Pi_{\ell,k}$ and 
$\Pi_{{k \choose \ell},1}$ are continuous, and as the space $\Sellk$ is a quotient of $\Mellk$, so is $\mathfrak{i}_{\ell,k}$ also continuous. 
\hspace*{\fill}$\square$
\begin{proposition} Let $\ell \leqslant k$. Then the set of submodules in $\Sellk$ that are free of rank $\ell$ is open dense in $\Sellk$. 
\end{proposition}
\noindent Proof: By Proposition 3.4, a submodule in $\Sellk$ is free of rank $\ell$, if and only if its cancellation ideal is nonzero. By Lemmas 5.1 and 5.4, the complement of this set of submodules, viz. $\frak{i}_{\ell,k}^{-1}(0)$, is closed in $\Sellk$. \hspace*{\fill}$\square$
\begin{corollary} The hypothesis in Theorem 3.2, that the cancellation ideal be nonzero, holds for a Zariski open dense set of submodules of $\Sellk$. 
\end{corollary}

Similar results hold, when $\ell \geqslant k$, for characteristic ideals of submodules of $A^k$ in $\Sellk$. Now let $\mathfrak{m}^t_{\ell,k}: \Mellk \rightarrow \A^{\ell \choose k}$ map a matrix to its $k \times k$ minors (written as a column in some fixed order). This map is continuous.  It again descends to a continuous map $\frak{i}^t_{\ell,k}: \Sellk \rightarrow \mathcal{I}_{\ell \choose k}$, mapping a submodule $P$  to the ideal generated by the $k \times k$ minors of any matrix in $\Mellk$ whose rows generate $P$. Its image $\frak{i}^t_{\ell,k}(P)$ is  the characteristic ideal $\frak{i}_k(P)$ of the submodule $P$.

Analogous to Proposition 5.1, we have the following result.
\begin{proposition} Let $\ell \geqslant k$. Then the set of submodules in $\Sellk$ whose characteristic ideals are nonzero, is open dense.
\end{proposition}
We immediately conclude a genericity result.
\begin{theorem} Let $\ell \geqslant k$. Then the set of uncontrollable systems in $\Bellk$ is open dense. In other words, uncontrollability is generic in $\Bellk$ when $\ell \geqslant k$.
\end{theorem}
\noindent Proof: By Corollary 4.5, a system $\ker_\F(P)$ is uncontrollable if and only if its characteristic ideal $\frak{i}_k(P)$ is nonzero. \hspace*{\fill}$\square$
\vspace{2mm}

The above construction of the topological space $\Mellk$, and of the spaces $\Sellk$ and $\Bellk$, assumes that a perturbation of a system governed by $\ell$ laws results again in a system governed by the same number $\ell$ of laws. We could drop this assumption as follows. Consider a differential operator in $\Mellk$ as an element in $\mathcal{M}_{\ell',k}$, where $\ell' > \ell$, by appending to its $\ell$ rows the $(\ell' - \ell) \times k$ matrix $0_{\ell'-\ell,k}$, whose every entry is the 0 operator in $A$. Then $\Mellk \hookrightarrow \mathcal{M}_{\ell',k}$ continuously, with image a proper Zariski closed subspace. This construction descends to the level of submodules of $A^k$, via the projections $\{\Pi_{\ell,k}, \ell \geqslant 1\}$, and hence to distributed systems in $\F^k$. Thus we have the chain of inclusions 
\[
\mathcal{M}_{1,k} \hookrightarrow \mathcal{M}_{2,k} \hookrightarrow \cdots \hookrightarrow \mathcal{M}_{\ell-1,k} \hookrightarrow \Mellk \hookrightarrow \cdots
\]
and the corresponding chains 
\[
\mathcal{S}_{1,k} \hookrightarrow \mathcal{S}_{2,k} \hookrightarrow \cdots \hookrightarrow \mathcal{S}_{\ell-1,k} \hookrightarrow \Sellk \hookrightarrow \cdots
\] 
\[
\mathcal{B}_{1,k} \hookrightarrow \mathcal{B}_{2,k} \hookrightarrow \cdots \hookrightarrow \mathcal{B}_{\ell-1,k} \hookrightarrow \Bellk \hookrightarrow \cdots
\]
of submodules of $A^k$ and systems in $\F^k$.
We can now consider the directed systems given by $\{\Mellk\}, \{\Sellk\}, \{\Bellk\}, ~\ell = 1, 2 \ldots$ Their direct limits, denoted $\mathcal{M}(k)$, $\mathcal{S}(k)$ and $\mathcal{B}(k)$, are the spaces of matrices with $k$ columns and entries in $A$, of submodules of $A^k$, and of distributed systems in $\F^k$, respectively, all with the Zariski topolgy.
Just as above, a subset of differential operators in $\mathcal{M}(k)$, or submodules in $\mathcal{S}(k)$ or systems in $\mathcal{B}(k)$, is said to be generic if it contains an open dense subset in this topology.
\vspace{1mm}

We can therefore either study perturbations of a system defined by $\ell$ laws 
within the space of systems all defined by $\ell$ laws, i.e. within $\mathcal{B}_{\ell,k}$, or study perturbations 
which might result in different numbers of laws by studying the inclusion $\Bellk \hookrightarrow \mathcal{B}(k)$.

However, as $\mathcal{B}_{\ell-1,k}$ embeds in $\mathcal{B}_{\ell,k}$ as a proper Zariski closed subset, the topology we have defined implies that, generically, the number of laws governing a system cannot decrease, but can only increase. 
%(This is akin to the statement that the number 0 can be easily perturbed to become nonzero, but it is unlikely that a nonzero number, when perturbed, will become 0.) 
In other words, in the space of systems that are described by $\ell$  laws, 
those systems that could be described by a fewer number of laws is a Zariski
closed subset, and those that need to be described 
by all the $\ell$ laws is open dense.

Here we are concerned only with properties of systems belonging to an open dense set. Moreover, answers to questions might depend upon whether the system is underdetermined or overdetermined, for instance Theorem 5.1 is applicable to overdetermined systems. Hence
we confine ourselves to the case where the number of laws defining a 
system remains constant under perturbations. In other words, we study genericity questions within $\Bellk$.
\vspace{2mm}

For the rest of this section, we consider strictly underdetermined systems, i.e. the case when $\ell < k$, and show that controllability is generic in $\Bellk$.

The proof rests on Theorem 3.2(ii). The cancellation ideal $\frak{i}_\ell(P)$ of a submodule $P$ in $\Sellk$ is generically nonzero by Proposition 5.1, and is generated by the  $k \choose \ell$ many maximal minors of any $P(D) \in \Mellk$, whose $\ell$ rows is a minimum set of generators for  $P$. A theorem of Macaulay (for instance Appendix to $\S 13$ in \cite{ma}) states that the codimension of such a \index{determinantal ideal}determinantal  ideal, if  proper, is bounded above by $k-\ell +1$. We show that when $(k-\ell+1) \leqslant n$,  there is a Zariski  open subset of $\Mellk$ where either the cancellation ideal attains this  codimension, or is equal to $A$. On the other hand, when $(k-\ell+1) > n$, we show that there is a  Zariski open subset of $\Mellk$ where the cancellation ideal equals $A$. Together, this will imply that controllability is generic in $\Bellk$, for all $\ell < k$, because $k-\ell+1$ is at least equal to 2. 
\vspace{1mm}

We begin the proof with a few elementary observations.
\begin{lemma}
Let $\mu: \A^r \times \A^r \rightarrow \A$ be the map $((a_i), (b_i)) \mapsto \sum_1^r a_ib_i$. Then $\mu ^{-1} (0)$ is a proper Zariski closed subset of the space $\A^r \times \A^r$.
\end{lemma}
 \noindent Proof: The map $\mu$ restricts to $\mu(d): A(d)^r \times A(d)^r \rightarrow A(2d)$. It is given by adding and multiplying the coefficients of  $a_i$ and $b_i$, and hence is continuous. The point 0 is closed in $A(2d)$, hence $\mu(d) ^{-1} (0)$ is closed in $A(d)^r \times A(d)^r$. The direct limit of $\{\mu(d)^{-1}(0)\}_{d=0,1,\ldots}$ equals $\mu^{-1}(0)$, hence it is closed in $\A^r \times \A^r$.
 \hspace*{\fill}$\square$ 
\begin{lemma}
Let $I$ be a proper ideal of $A$. Then $I$ is a proper  Zariski closed subset of the space $\A$. Hence the set of elements of $\A$ which are not in $I$ is an open dense subset of  $\A$. 
\end{lemma}
\noindent Proof: Let $I$ be generated by $\{a_1, \ldots ,a_r\}$, and let $a$ denote the point $(a_1, \ldots, a_r) \in A^r$. Define the map $\mu_a: \A^r \rightarrow \A$ by $\mu_a(b_1, \ldots ,b_r) = \sum a_ib_i$. This is a $\mathbb{C}$-linear map as the coefficients of the operator $\sum a_ib_i$ are $\mathbb{C}$-linear combinations of the coefficients of the $b_i$. Its image is precisely the ideal $I$. 

Let the maximum of the degrees of the $a_i$ be $s$, then $s \geqslant 1$ as $I$ is proper. For each $d$, the map $\mu_a$ restricts to a map $\mu_a(d): A(d)^r \rightarrow A(d+s)$. Its image, say $I_a(d + s)$, is contained in $I \cap A(d + s)$.
As $I$ is proper, $I_a(d + s)$ is a proper linear subspace of $A(d+s) (\simeq \mathbb{C}^{N(d+s)}$), for all $d \gg 0$. It is a proper Zariski closed subset of $A(d+s)$, and its vanishing ideal is generated by linear forms (in the indeterminates $\{X_{d_1 \cdots d_n}\}$, and whose coefficients are polynomial functions of the coefficients of  the operators $a_1, \ldots, a_r$). Its complement is then open dense. 
For $d_1 < d_2$, the map $\mu_a(d_2)$ restricts to $\mu_a(d_1)$, hence $I_a(d_1 + s)  \subset I_a(d_2 + s)$. The direct limit of the closed subspaces $\{I_a(d + s)\}_{d = 0, 1, \ldots}$ is $I$, hence $I$ is a proper Zariski closed subset of $\A$, and its complement is open dense.  
\hspace*{\fill}$\square$
\vspace{2mm}

\noindent Remark 5.2:  More generally, every $\mathbb{C}$-linear subspace of $\A$ is Zariski closed, given by the common zeros of linear forms in the indeterminates $\{X_d\}$. So is therefore every affine linear subset of $\A$. 
\vspace{1mm}

\begin{proposition}
For $a = (a_1, \ldots , a_r) \in \A^r$, let $I_a$ be the Zariski closed subset of $\A$ consisting of the elements in the ideal $(a_1, \ldots ,a_r)$ generated by the coordinates $a_1, \ldots , a_r$ of $a$ $($namely, the above lemma$)$. Then the set $G = \{(a, a_{r+1}) ~| ~ a \in \A^r, ~ a_{r+1} \notin I_a\}$ contains a Zariski open subset of  $\A^r \times \A$.
\end{proposition}
\noindent Proof: Fix $s \geqslant 0$. Each $x \in A(s)^r$ defines, for every $d \geqslant 0$, a $\C$-linear map $\mu_x(d): A(d)^r \rightarrow A(d+s)$. If we represent this map by a matrix, say $M_x(d)$, the image of $\mu_x(d)$ is the column span  of $M_x(d)$, and we denote it by $I_x(d+s)$ (as in the above lemma). There is a Zariski open dense subset of $A(s)^r$, say $U(s)$, such that the matrix $M_x(d)$, corresponding to every $x \in U(s)$, is of full column rank.   The column span of $M_x(d)$ is a linear subspace of $A(d+s)$, and is a variety defined by linear forms. These linear forms belong to the kernel of the transpose $M_x^{\t}(d)$ (as the annihilator of an image of a linear map equals the kernel of its transpose). By `Cramer's rule', this kernel is spanned by multi-linear forms in the coefficients of $M_x(d)$.  Hence, as $x$ varies in $U(s)$, the subset  $\{(x, x_{r+1}) ~|~ x \in U(s), x_{r+1} \in I_x(d+s) \}$ of $A(s)^r \times A(d+s)$ is \index{locally closed}locally closed (i.e. the intersection of a closed set with an open set). The subset $\{(x, x_{r+1}) ~|~ x \in U(s), x_{r+1} \notin I_x(d+s)\}$ is  then open in $U(s) \times A(d+s)$, for all $d$, and thus its direct limit (with respect to $d$), $G(s) = \{(x, x_{r+1}) ~| ~ x \in U(s), x_{r+1} \notin I_x\}$,  is open in $U(s) \times \A$.

For $s < s'$, $G(s) \subset G(s')$, and the direct limit of $\{G(s)\}_{s=0,1, \ldots}$ equals $G$. This proves the proposition.  \hspace*{\fill}$\square$
\begin{lemma}
Let $I$ be a proper ideal of $A$. Then the set of zero divisors on $A/I$ is closed in the space $\A$, and hence the set of nonzero divisors on 
$A/I$ is open dense.
\end{lemma}
\noindent Proof: The set of zero divisors on $A/I$ is the union of its finite number of associated primes, and this finite union is closed in $A$ (by Lemma 5.6).
\hspace*{\fill}$\square$
\vspace{1mm}

Suppose $I$ is a maximal ideal of $A = \C[D_1, \ldots, D_n]$, it can then be generated by $n$ elements, say $I = (a_1, \ldots ,a_n)$.  The ideal $I + (a)$ is  equal to $A$ if and only if the element $a$ does not belong to $I$. Thus $I + (a) = A$, for $a$ in an  open dense subset of $\A$ (by Lemma 5.6). %This is analogous to saying that the set of units in $\C$ is Zariski open, because $A/I \simeq\mathbb{C}$.
Furthermore, let $M$ be the set of points $a = (a_1, \ldots ,a_n) \in \A^n$ such that the ideal $(a_1, \ldots ,a_n)$ generated by its coordinates is maximal (we show below in Proposition 5.7 that $M$ contains an open dense subset of  $\A^n$). Then by the above proposition, the set of points $(a_1, \ldots , a_n, a_{n+1}) \in A^{n+1}$ such that $(a_1, \ldots ,a_n) \in M$ and the ideal generated by these $n+1$ coordinates is proper, is closed in the subspace $M \times \A$ of $\A^n \times \A$. 

The opposite is however the case for proper ideals generated by fewer than $n$ elements.% namely the proposition below. and is suggested by the following heuristic:Let $I$ be a proper ideal generated by $r < n$ elements, then the codimension of $I$ is at most $r$ (by Krull's Theorem), and the dimension of $A/I$ is at least $n-r$. By Noether normalization, $A/I$ is isomorphic to an integral extension of a polynomial ring with number of indeterminates at least $n-r$. Therefore, the set of units in $A/I$ is contained in a proper Zariski closed set, and the set of elements $a \in A$ such that the sum $I + (a)$ is not equal to $A$ contains an open dense subset of the space $\A$. Indeed, we have the following proposition.
\begin{proposition}
Let $r \leq n$. and let $P_r$ be the set of elements $(a_1, \ldots ,a_r) \in \A^r$ such that the ideal $(a_1, \ldots ,a_r)$ generated by its coordinates is a proper ideal of $A$. Then $P_r$ contains an open dense subset of $\A^r$.
\end{proposition}

We use the following result of Brownawell \cite{b}:

\noindent Theorem (Brownawell): {\em Suppose the ideal generated by the polynomials $p_1, \ldots , p_m$ equals $\mathbb{C}[x_1, \ldots, x_n]$, where the degree of the $p_i$ is less than or equal to $d$. Then there are polynomials $q_1, \ldots ,q_m$ such that $\sum p_i q_i = 1$, where the degree of the $q_i$ is less than or equal to $D = n^2 d^n + nd$.} 
\vspace{1mm}

\noindent Proof of proposition:  It suffices to prove the statement for $r = n$, for by Lemma 5.2, if $P_n$ contains an open dense subset of $\A^n$, then its projection to $\A^r$ also contains an open dense subset. The coordinates of these points generate proper ideals of $A$. Thus this projection of $P_n$ is contained in $P_r$.

We first show that $P_n \cap A(d)^n$ contains a nonempty euclidean open subset of $A(d)^n$. 

To say that the ideal $I_a$ generated by the coordinates of a point $a = (a_1, \ldots ,a_n) \in A(d)^n$ is a proper ideal of $A$ is to say that the map $\sigma_a:\mathbb{C}^n \rightarrow \mathbb{C}^n$ defined by $x \mapsto (a_1(x), \ldots, a_n(x))$ has nonempty inverse image $\sigma_a^{-1}(0)$.
Consider a point $a_o$ in $A(d)^n$ such that the corresponding map $\sigma_{a_o}$ has a regular point $x$ in $\sigma_{a_o}^{-1}(0)$. The rank of $\sigma_{a_o}$ equals $n$ at $x$, hence all smooth maps sufficiently close to $\sigma_{a_o}$ in the \index{compact-open topology}compact-open topology also include 0 in their images (Inverse Function Theorem). Restricting to maps given by elements in $A(d)^n$ as above, this means that there is an open neighbourhood of  $a_o$ in the {\em euclidean topology} on $A(d)^n$ such that the variety of the ideal generated by the coordinates of every point in it, is nonempty. All these ideals are thus proper ideals of $A$.

We now show that $P_n \cap A(d)^n$ contains a nonempty Zariski open subset.

Suppose that $I_a = A$ for some $a = (a_1, \ldots, a_n)$ in $A(d)^n$. 
Then by the theorem of Brownawell, there are elements  $b_1, \ldots ,b_n$ of $A$, of degree bounded by $\Delta = n^2 d^n + nd$, such that $\sum a_i b_i = 1$. This is a Zariski closed condition on the coefficients of the $a_i$ and $b_i$ (as in Lemma 5.5), and it defines a proper affine variety in the affine space $A(d)^n \times A(\Delta)^n$. Its projection $C$ to the first $n$ coordinates is a \index{constructible set}constructible set in $A(d)^n$ (by Chevalley's Theorem 1.8.4 in EGA IV, the image of a variety is a constructible set, i.e., a finite union of locally closed sets). The coordinates of a point in $C$ generates the unit ideal $A$, and the coordinates of points in the complement $C' = A(d)^n \setminus C$, which is also constructible, generate proper ideals of $A$. 

We have shown at the outset that $C'$ contains a euclidean open subset; as it is constructible, it must therefore contain a nonempty Zariski open subset of $A(d)^n$. This is true for every $d$, hence the set of points in $A^n$ which generate a proper ideal of $A$, contains an open dense subset of $\A^n$.
 \hspace*{\fill}$\square$
\begin{proposition}
Let $r < n$, and let $U_{r+1}$ be the set of points $(a_1, ..., a_{r+1}) \in P_{r+1}$ such that $a_1 \neq 0$ in $A$, $a_2 \neq 0$ in $A/(a_1)$, ... and $a_{r+1} \neq 0$ in $A/(a_1, \ldots ,a_r)$. Then $U_{r+1}$ is open dense in $P_{r+1}$, and hence contains an open dense subset of $\A^{r+1}$.
\end{proposition}
\noindent Proof:  The point 0 is closed in $\A$, hence the statement is true for $r = 0$. Assume by induction that $U_r$ is open dense in $P_r$. The subset $\{(a,b_{r+1}) ~| ~a \in U_r, ~b_{r+1} \in I_a\}$ is contained in the complement of the  set $G$ of Proposition 5.3, it is thus contained in a closed subset of $U_r \times \A$. Its complement is then open dense in $U_r \times \A$; hence $U_{r+1}$, which is the intersection of this complement with $P_{r+1}$, is open dense in $P_{r+1}$, and hence contains a dense open subset of $\A^{r+1}$.   
\hspace*{\fill}$\square$
\begin{proposition}
Let $r < n$, and let $N_{r+1}$ be the set of points $(a_1, \ldots , a_{r+1}) \in U_{r+1}$ such that $a_1$ is a nonzero divisor (nzd) on $A$, $a_2$ is a nzd on $A/(a_1)$, ... and $a_{r+1}$ is a nzd on $A/(a_1, \ldots ,a_r)$. Then $N_{r+1}$ is open dense in $U_{r+1}$, and hence contains an open dense subset of $\A^{r+1}$. 
\end{proposition}
\noindent Proof: The statement is true for $r = 0$ as now $N_1$ equals $U_1$ ($A$ is an integral domain). Assume by induction that  $N_r$ is open dense in $U_r$. For $a = (a_1, \ldots ,a_r)$ in $N_r$, let $Z_a$ be the set of zero divisors on $A/(a_1, \ldots ,a_r)$; it is a closed subset of the space $\A$ by Lemma 5.7. We need to show that the set $Z = \{(a, b_{r+1}\} ~| ~a \in N_r, b_{r+1} \in Z_a\}$ is closed in $N_r \times \A$, for its complement in $U_{r+1}$ is precisely $N_{r+1}$. 
To show this, it suffices by Proposition 5.3 to show that every element in the vanishing ideal of $Z_a$ is a polynomial function of elements in the vanishing ideal of the closed set $I_a$ of the space $\A$.

Let $\mu: \A \times \A \rightarrow \A$ be the multiplication map, mapping $(a, b)$ to $a b$ (the map of Lemma 5.5 for $r = 1$). Denote also by $\mu$ its restriction $\mu : \A \times (\A \setminus I_a) \rightarrow \A$. Then $V = \mu ^{-1}(I_a)$ is locally closed in $\A \times \A$, as it is is a Zariski closed subset of the space $\A \times (\A \setminus I_a)$. Let its vanishing ideal be $J$; its elements are polynomial functions of the elements of the vanishing ideal of $I_a$, and therefore polynomial functions of the coefficients of the components $a_1, \ldots , a_r$ of $a$. 

The set $Z_a$ is the projection of $V$ to the first factor $\A$. In general a projection is not Zariski closed as the space $\A$ is not complete, but here $Z_a$ is indeed closed. Hence its vanishing ideal equals $i ^{-1} (J)$, where $i$ is the inclusion of $\A$ in $\A \times (\A \setminus I_a)$.  

As $i$ is an algebraic map, it follows that elements of  $i ^{-1} (J)$ are polynomial functions of the coefficients of $a_1, \ldots ,a_r$. This completes the proof of the proposition. 
\hspace*{\fill}$\square$
\vspace{2mm}

Recall the definition of a \index{Cohen-Macaulay}Cohen-Macaulay ring:  a sequence $a_1, \ldots ,a_t$ in a Noetherian ring $R$ is regular if (i) the ideal $(a_1, \ldots , a_t) \subsetneq R$, and (ii) $a_1$ is a nonzero divisor in $R$, and for each $i$, $2 \le i \le t$, $a_i$ is a nonzero divisor on $R/(a_1, \dots , a_{i-1})$. Thus, the above proposition  asserts that for $r < n$, the set of points in $\A^{r+1}$ corresponding to regular sequences contains an open dense subset. The \index{depth}depth of an ideal $I$ is the length of any maximal regular sequence in $I$. The ring $R$ is Cohen-Macaulay if for every ideal $I$ of $R$, $\dep(I) = \codim(I)$. 
\vspace{2mm}

\begin{proposition}
For $r \leqslant n$, the set of points $(a_1, \ldots , a_r) \in \A^r$ such that the ideal generated by its coordinates has codimension $r$, contains an open dense subset. In particular, the set $M$ of points $(a_1, \ldots ,a_n)$ in $\A^n$ such that the ideal $(a_1, \ldots ,a_n)$ is maximal, contains an open dense subset of $\A^n$.
\end{proposition}
\noindent Proof: The set of points in $\A^r$ corresponding to regular sequences contains an open dense subset of $\A^r$. As the ring $A$ is Cohen-Macaulay, the codimension of an ideal generated by such a sequence equals $r$. The statement on the set $M$ of the proposition now follows as the dimension of $A$ equals $n$.
\hspace*{\fill}$\square$
\begin{corollary}
Let $r > n$. Then the set of points $(a_1, \ldots ,a_r)$ in $\A^r$ such that the ideal $(a_1, \ldots ,a_r)$ equals $A$, contains a Zariski open subset of $\A^r$.   
\end{corollary}
\noindent Proof: It suffices to prove the statement for $r = n+1$. By the remarks preceding Proposition 5.4, the set of points $a \in M \times \A$ such that its coordinates generate a proper ideal is a proper closed subset, hence its complement is open in $M \times \A$, and so contains an open dense subset of $\A^{n+1}$.  \hspace*{\fill}$\square$
\vspace{3mm}

We can now prove the main theorem of the section.

\begin{theorem} 
(i) Let $\ell < k$ be such that $k - \ell + 1 \leqslant n$. Then the set of differential operators in $\Mellk$ which define controllable systems contains an open dense subset of $\Mellk$. 

(ii) Let $k - \ell + 1 > n$. Then an open dense set of differential operators in $\mathcal{M}_{\ell,k}$ define strongly controllable systems.
\end{theorem}

\noindent Proof: (i) Set $r ={k \choose \ell}$ and $s = k - \ell  + 1$.   Let $\Mellk^*$ be the subset of those operators in $\Mellk$ whose cancellation ideals are nonzero. By Lemmas 5.2 and 5.3, $\Mellk^*$ is an open dense subset of $\Mellk$.

Suppose $P(D) \in \Mellk^*$ is such that its cancellation ideal $\frak{i}_\ell(P(D))$ is equal to $(1)$. Then by Proposition 3.3, $\ker_\F(P)$ is strongly controllable. 
 
Assume now that the cancellation ideal of $P(D)$ is a nonzero proper ideal of $A$.  By a theorem of Macaulay quoted earlier, its codimension is bounded by $s$. We show that this bound is attained by an open  subset of $\mathcal{M}^*_{\ell,k}$.

Consider the following matrix in $\mathcal{M}_{\ell,k}$

\[ P(D)  = \begin{pmatrix} a_1 & a_2 & \cdots & a_s & 0 & \cdots & 0 \\
0  & 0 & \cdots & 0 & & &\\ 
\vdots &  \vdots & & \vdots & & I_{\ell-1}\\
0 &  0 & \cdots & 0 &&  &
\end{pmatrix} 
\]
where $I_{\ell-1}$ is the $(\ell-1) \times (\ell-1)$ identity matrix, and the $a_i$ are arbitrary nonzero elements of $A$. Its nonzero maximal minors are $\{a_1, a_2, \ldots , a_s \}$, hence its cancellation ideal equals $(a_1, a_2, \ldots, a_s)$. 

As $s \leq n$, the image of the composition  $\mathcal{M}^*_{\ell,k} \stackrel {\frak{m}_{\ell,k}} {\longrightarrow} \A^r \stackrel{\pi}{\longrightarrow} \A^s$ contains the open  set of points in $\A^s$ whose coordinates are regular sequences of length $s$ (Proposition 5.7);  here the second map is the projection of $\A^r$ to  $\A^s$, where the indices of $s$ are determined by the nonzero maximal minors described above, in the ordered set of all the $r$ maximal minors. Hence there is an open subset of operators in $\mathcal{M}^*_{\ell,k}$ with the property that the ideal generated by the $s$ minors described above is of codimension $s$. The cancellation ideals of all the operators in this open subset of $\mathcal{M}^*_{\ell,k}$ then have codimension at least $s$, and hence codimension equal to $s$. 
As $s \geqslant 2$, these operators define controllable systems by Theorem 3.2(ii).

(ii) If $s =  k - \ell + 1 > n$, then by Corollary 5.2 there is an open dense subset of points in $\A^s$ whose coordinates generate the unit ideal in $A$. Let $(a_1, \ldots , a_s)$ be such a point; there is a matrix just as in (i) above, whose maximal minors are these $a_i$.  Hence the set $U$ of matrices in
 $\mathcal{M}^*_{\ell,k}$ such that for every $P(D) \in U$, the ideal generated by the corresponding $s$ minors
equals $A$, is nonempty open, and so open dense. Thus $\frak{i}_{\ell,k}(P(D)) = A$, and $\ker_\F(P)$  is strongly controllable. \hspace*{\fill}$\square$
\vspace{2mm}

This genericity result in the space of operators descends to the space of distributed systems.

\begin{theorem}
(i) (i) Let $\ell < k$ be such that $k - \ell + 1 \leqslant n$. Then the set of distributed systems in $\Bellk$ which define controllable systems contains an open dense subset of $\Bellk$. 

(ii) Let $k - \ell + 1 > n$. Then an open dense set of systems in $\Bellk$ define strongly controllable systems.
\end{theorem}
\noindent Proof: (i) We show that the image of the set of operators defining controllable systems in $\Mellk$ under the map $\Pi_{\ell,k}$ contains an open subset of $\Sellk$. We refer to the commutative diagram (11) in the argument below.

Let $s=k-\ell +1$ as before, and let $X \subset \mathcal{I}_s$ be the subset of those ideals whose codimensions equal $s$.
By the above theorem, it suffices to show that $X$ contains an open subset of $\mathcal{I}_s$. As $\mathcal{I}_s$ is equipped
with the quotient topology given by the surjection $\Pi_{s,1}: \A^s \rightarrow \mathcal{I}_s$, we
need to show that $U := \Pi_{s,1}^{-1}(X)$ contains an open subset of $\A^s$.

Let $I \in X$, and let $(a_1, \ldots ,a_s) \in \Pi_{s,1}^{-1}(I)$. If $a_1, \ldots ,a_s$ is a
regular sequence, then there is nothing to be done by Proposition 5.6.  Suppose that it is not.
As $A$ is Cohen-Macaulay, $\dep(I) = s$, hence there is a regular sequence $a'_1, \ldots ,a'_s$
in $I$.
Each $a'_i$ is an $A$-linear combination of $a_1, \ldots ,a_s$, hence it follows that there is an
$A$-linear map $L: \A^s \rightarrow \A^s$ such that $L(a_1, \ldots ,a_s) = (a'_1, \ldots ,a'_s)$.
Again by Proposition 5.6, there is a neighbourhood $W$ of $(a'_1, \ldots ,a'_s)$ such that the coordinates of each point in it is a regular sequence. The map $L$ is continuous, hence
$L^{-1}(W)$ contains an open neighbourhood of $(a_1, \ldots ,a_s)$;  let $(b_1, \ldots ,b_s)$ be
a point in it. Then $L(b_1, \ldots ,b_s) = (b'_1, \ldots ,b'_s)$ is in $W$, and so the sequence
$b'_1, \ldots ,b'_s$ is regular. Each $b'_i$ belongs to the ideal $J = (b_1, \ldots , b_s)$,
hence $\dep(J) = s$, and so also the codimension of $J$ equals $s$. Thus  $L^{-1}(W) \subset U$, and $U$ is open.

(ii) The proof follows from Theorem 5.2(ii); it  suffices to observe that by Corollary 5.2, the ideal $(1)$ is open in $\mathcal{I}_s$, as $s > n$.
\hspace*{\fill}$\square$\\

We summarise the principal results of this section (Theorems 5.1 and 5.3) in the following statement.

\begin{theorem} Controllability is generic in $\Bellk$ when $\ell < k$. On the other hand, uncontrollability is generic in $\Bellk$ when $\ell \geqslant k$.
\end{theorem}

\hspace*{\fill}$\square$\\
%\end{corollary}

\newpage

\section{Pathologies}
\vspace{.2cm}

We have already pointed out several instances where some fact which is the case for a signal space $\F$ that is an injective cogenerator, is not true in general (for instance Proposition 2.6 and Example 2.4). We now point out in more detail, in the context of the Sobolev spaces, the problems that can arise.
\vspace{2mm}

Let $k$ be a field and $R$ a commutative $k$-algebra. Let $F$ be an $R$-module and $\d(F)$ = $\homo_k(F,~k)$ its algebraic dual. The functor $\d$, from the category of $R$-modules to itself, is contravariant  and exact, . Let 
\[
0 \rightarrow M_1 \longrightarrow M_2 \longrightarrow M_3 \rightarrow 0
\]
be an exact sequence of $R$-modules. Suppose that $F$ is a flat $R$-module. Then 
\[
0 \rightarrow  F \otimes _R M_1 \longrightarrow
F \otimes _R M_2 \longrightarrow 
F \otimes _R M_3  \rightarrow 0
\]
is exact, so that  
\[
0 \rightarrow \d(F \otimes _R M_3) \longrightarrow
\d(F \otimes _R M_2) \longrightarrow 
\d(F \otimes _R M_1) \rightarrow 0
\]
is also exact. But $\d(F \otimes_R M_i) =$
$\homo_k(F \otimes_R M_i,~ k) \simeq$
$\homo_R (M_i, ~\homo_k(F, ~k))$ by the adjointness of the pair $(\homo,~\otimes)$. Thus
\[
 0 \rightarrow {\homo}_R(M_3,~\d(F)) \longrightarrow
{\homo}_R (M_2,~\d(F)) \longrightarrow
{\homo}_R(M_1,~\d(F)) \rightarrow 0 
\]
is exact, which is to say that $\d(F)$ is an injective 
$R$-module. Reversing the above argument shows that $F$ is flat if and only if its algebraic dual $\d(F)$ is injective.

Suppose further that $F$ is faithfully flat. Then if $M$ is any nonzero module, $F\otimes_R M$ is nonzero, so that $\d(F\otimes_R M) \simeq \homo_R (M, ~\d(F))$ is also nonzero, which is to say that $\d(F)$ is an injective cogenerator. Reversing this shows that $F$ is faithfully flat if and only if $\d(F)$ is an injective cogenerator.

Thus one could say that in the category of $R$-modules, flatness and injectivity are adjoint properties, as are also the properties of faithful flatness and injective cogeneration.

Suppose now that $k=\C$ and $R = A=\mathbb{C}[D_1, \ldots ,D_n]$. We recollect the $A$-module structure of the classical spaces of distributions from Section 2, on which this entire theory of control rests: 
\vspace{2mm}

\noindent 1. The locally convex topological vector space $\mathcal{D}'$ of distributions, as well as the space $\mathcal{C}^\infty$ of smooth functions, are injective cogenerators as $A$-modules; 

\vspace{1mm}
\noindent 2. Their {\it topological} duals $\mathcal{D}$ and $\mathcal{E}'$, respectively, of compactly supported smooth functions and distributions are faithfully flat $A$-modules; 

\vspace{1mm}
\noindent 3. The space $\mathcal{S}'$ of tempered distributions is an injective $A$-module that is not a cogenerator. Its topological dual $\mathcal{S}$, the Schwartz space of rapidly decreasing smooth functions, is a flat module that is not faithfully flat.
\vspace{2mm}

As all the topological $A$-modules listed above are reflexive topological $k$-vector spaces, the question arises as to whether flatness and injectivity are also adjoint properties in the category of topological $A$-modules. 

We show that the answer to this question is negative by exhibiting an elementary counter-example of a flat topological module whose topological dual is not injective, but which is again flat. \\

For every $s$ in $\mathbb{R}$, the Sobolev space 
\index{$\mathcal{H}^s$}$\mathcal{H}^s$ on $\mathbb{R}^n$ of order $s$ is the space of temperate
distributions $f$ whose Fourier transform $\hat{f}$ is a measurable 
function such that 
\[
\|f\|_s = \big(\frac{1}{(2\pi)^n} \int_{\mathbb{R}^n} |\hat{f}(\xi)|^2 (1+|\xi|^2)^s d\xi \big)^\frac{1}{2} < \infty
\]
$\mathcal{H}^s$ is a Hilbert space with norm $\| \phantom{.} \|_s$. When
$s>t$, $\mathcal{H}^s \hookrightarrow \mathcal{H}^t$ is a continuous inclusion. If $p(D)$ is an element of $A$ of order $r$, then it maps $\mathcal{H}^s$ into $\mathcal{H}^{s-r}$.
If the family $\{\mathcal{H}^s, s \in \mathbb{R}\}$ is considered an increasing family of vector spaces indexed by the directed set $\mathbb{R}$, then its 
direct limit $\overrightarrow{\mathcal{H}} := \varinjlim {\mathcal{H}^s}$
is the union $\bigcup_{s \in \mathbb{R}}\mathcal{H}^s$ of all the Sobolev spaces and is an $A$-module that is strictly contained in 
$\mathcal{S}'$. Instead, if the family $\{\mathcal{H}^s, s \in \mathbb{R}\}$ 
is considered a decreasing family of vector spaces, then its inverse limit 
$\overleftarrow{\mathcal{H}} := \varprojlim {\mathcal{H}^s}$ is the 
intersection $\bigcap_{s \in \mathbb{R}}\mathcal{H}^s$ of the Sobolev spaces, and is again an $A$-module. This intersection contains the 
Schwartz space $\mathcal{S}$ but is strictly larger than it. As $\mathbb{Z}$
is cofinal in $\mathbb{R}$ in either of the above two situations, these limits are also the direct and inverse limits of the countable family 
$\{\mathcal{H}^s, s \in \mathbb{Z}\}$.

Further, if each $\mathcal{H}^s$ is given its Hilbert space topology, then
the union $\overrightarrow{\mathcal{H}}$ will be equipped with the direct limit topology, i.e. the strongest topology so that each 
$\mathcal{H}^s \hookrightarrow \overrightarrow{\mathcal{H}}$ is continuous. 
(As the topology that $\mathcal{H}^s$ inherits from $\mathcal{H}^t$ for
$s>t$ is strictly weaker, this direct limit is {\em not} a strict direct limit.) With this topology $\overrightarrow{\mathcal{H}}$ is a 
locally convex space which is bornological and barrelled. On the other hand, the intersection $\overleftarrow{\mathcal{H}}$ with the inverse limit topology, i.e. the weakest topology such that each inclusion 
$\overleftarrow{\mathcal{H}} \hookrightarrow \mathcal{H}^s$ is continuous, is a Fr\'echet space. As the dual of 
$\mathcal{H}^s$ is $\mathcal{H}^{-s}$, $\overrightarrow{\mathcal{H}}$ and 
$\overleftarrow{\mathcal{H}}$
are duals of each other. Thus $(\overrightarrow{\mathcal{H}})' = \overleftarrow{\mathcal{H}}$ and
$(\overleftarrow{\mathcal{H}})' = \overrightarrow{\mathcal{H}}$, so that they are both reflexive.

\begin{lemma} $\overrightarrow{\mathcal{H}}$ and $\overleftarrow{\mathcal{H}}$ are both torsion free $A$-modules.
\end{lemma}
\noindent Proof: Suppose $f$ is any element in $\overrightarrow{\mathcal{H}}$
or $\overleftarrow{\mathcal{H}}$ such that $p(D)f = 0$ for some nonzero $p(D)$. Fourier transform implies that $p(\xi)\hat{f}(\xi) = 0$, which implies that the support of the measurable function $\hat{f}$ is contained in the real variety of the polynomial $p(\xi)$, a set of measure zero. Hence $\hat{f} = 0$, and so is therefore $f=0$.
\hspace*{\fill}$\square$
\vspace{2mm}

\begin{corollary} Let $A = \mathbb{C}[{\dbydt}]$, the $\mathbb{C}$-algebra of ordinary differential operators. Then the Sobolev limits 
$\overrightarrow{\mathcal{H}}(\mathbb{R})$ and 
$\overleftarrow{\mathcal{H}}(\mathbb{R})$ are both flat $A$-modules, but not faithfully flat.
\end{corollary}
\noindent Proof: As $\mathbb{C}[{\dbydt}]$ is a principal ideal domain, torsion free implies flat \cite{ma}. Thus it remains to show that the Sobolev limits are not faithfully flat. 

Let $p({\dbydt}) = 1+{\dbydts}$. Then for any element $f$ in either of the two Sobolev limits, the Fourier inverse of $(1+\xi^2)^{-1} \hat{f}(\xi)$ is also in the corresponding Sobolev limit. This implies that $p(\dbydt)$ defines a surjective morphism on the Sobolev limits, so that $\frak{m}\overrightarrow{\mathcal{H}}(\mathbb{R}) = \overrightarrow{\mathcal{H}}(\mathbb{R})$ and 
$\frak{m}\overleftarrow{\mathcal{H}}(\mathbb{R}) = \overleftarrow{\mathcal{H}}(\mathbb{R})$ 
for the maximal ideals $\frak{m}$ of $\mathbb{C}[{\dbydt}]$ that contain $p({\dbydt})$ . This shows that the two limits are not faithfully flat.
\hspace*{\fill}$\square$
\vspace{2mm}

\begin{lemma} $\overrightarrow{\mathcal{H}}$ and 
$\overleftarrow{\mathcal{H}}$ are not divisible (hence not injective)   $A$-modules.
\end{lemma}
\noindent Proof: Let $f$ in $\mathcal{S}'$ be such that $\hat{f}(\xi)$ equals 
$\sqrt{\xi}$ in a neighbourhood of 0 and is rapidly decreasing at infinity. Such an $f$ is in every Sobolev space $\mathcal{H}^s$, and so is in both $\overrightarrow{\mathcal{H}}$ as well as  
$\overleftarrow{\mathcal{H}}$. However, the Fourier inverse of  $\xi^{-1}\hat{f}(\xi)$ is not in any Sobolev space and is therefore not in either of the Sobolev limits. This implies that ${\dbydt}$ does not define a surjective morphism on the Sobolev limits, so that they are not divisible (hence not injective) modules.
\hspace*{\fill}$\square$
\vspace{2mm}

The answer to the question posed at the start is thus negative:

\begin{theorem} Flatness and injectivity are not adjoint properties in the category of topological $A$-modules.
\hspace*{\fill}$\square$
\end{theorem}

While the Sobolev limits on $\mathbb{R}$ are flat modules over the ring of differential operators, the corresponding fact is false in $\mathbb{R}^n, n\geqslant 2$.
\vspace{1.5mm}

\noindent Example 5.1 (The \index{deRham complex}Sobolev-deRham complex on $\mathbb{R}^3$):  Recollect the deRham complex (in the classical spaces) from Example 3.2. We now show that the Sobolev-deRham complex 
\[ \overrightarrow{\mathcal{H}} \stackrel{\grad}{\longrightarrow} (\overrightarrow{\mathcal{H}})^3 \stackrel{\curl}{\longrightarrow} (\overrightarrow{\mathcal{H}})^3 \stackrel{\div}{\longrightarrow} \overrightarrow{\mathcal{H}} \]
is {\it not} exact, and similarly for  
$\overleftarrow{\mathcal{H}}$.
Thus we exhibit an element in the kernel of 
$\curl:(\overrightarrow{\mathcal{H}})^3 \rightarrow (\overrightarrow{\mathcal{H}})^3$ which is not in the image of 
$\grad:\overrightarrow{\mathcal{H}} \rightarrow (\overrightarrow{\mathcal{H}})^3$, and similarly for 
$\overleftarrow{\mathcal{H}}$. By the equational criterion for flatness \cite{ma}, this will show that the Sobolev limits are not flat $A$-modules.

Let $h$ be an element in the kernel of $\curl:(\overrightarrow{\mathcal{H}})^3 \rightarrow (\overrightarrow{\mathcal{H}})^3$. 
As the deRham complex is exact in the space $\mathcal{S}'$ of temperate distributions, there is certainly an $f$ in $\mathcal{S}'$ (and unique upto an additive constant) such that 
$\grad (f) = h$. Thus $h =(\partial_xf,\partial_yf,\partial_zf)$. We claim that there is an $f$ in $\mathcal{S}'$ which is not in any Sobolev space but such that $\partial_xf$, $\partial_yf$ and $\partial_zf$ are all three in every one of them.

The Sobolev spaces are subspaces of $\mathcal{S}'$ on which is defined the Fourier transform.
Hence by Fourier transformation, the operators $\grad$ and $\curl$ become matrix operators with polynomial entries -
\[
\widehat{\grad} = \imath\begin{pmatrix}
     \xi_x     \\
     \xi_y     \\
     \xi_z     \\ 
\end{pmatrix} ;\widehat{\curl} =\imath\begin{pmatrix}
   0  & -\xi_z & \xi_y    \\
   \xi_z  &  0 & -\xi_x   \\
   -\xi_y &  \xi_x &  0   \\
\end{pmatrix}
\]
It now suffices to find an $\hat{f}$ (here $\hat{\phantom{oo}}$ denotes the Fourier transform) which is not in $\mathcal{L}^2(\mathbb{R}^3)$ but such that $\xi_x\hat{f}$, $\xi_y\hat{f}$ and $\xi_z\hat{f}$ are all in $\mathcal{L}^2(\mathbb{R}^3)$, with respect to the measure $(1+r^2)^sd\vol$, for every $s$, where $r=\sqrt{\xi_x^2+\xi_y^2+\xi_z^2}$~. But $\xi_x\hat{f}$, $\xi_y\hat{f}$ and $\xi_z\hat{f}$ are all in $\mathcal{L}^2$ if and only if $r\hat{f}$ is in $\mathcal{L}^2$. Thus it suffices to find an $\hat{f}$ which is not in $\mathcal{L}^2$ but such that $r\hat{f}$ is in $\mathcal{L}^2$.

This is elementary, for let $\hat{f}$ be any function which is rapidly decreasing at infinity, and which at 0 is $O(r^\alpha)$, $-\frac{5}{2} < \alpha \leqslant -\frac{3}{2}$.
Then 
\[
\int_{\mathbb{R}^3}|\hat{f}|^2 (1+r^2)^s d\vol = \int_{\mathbb{R}^3}|\hat{f}|^2(1+r^2)^s r^2drd\theta d\phi 
\]
This integral, in some neighbourhood $\Omega$ of 0, is therefore of the order of 
\[
\int_{\Omega}r^{2\alpha +2}dr d\theta d\phi
\]
which is not finite. 

On the other hand 
\[
\int_{\mathbb{R}^3}r^2|\hat{f}|^2(1+r^2)^sd\vol = \int_{\mathbb{R}^3}|\hat{f}|^2r^4(1+r^2)^sdrd\theta d\phi
\]
is finite. 
Thus this $f$ is not in any 
$\mathcal{H}^s$ whereas $\grad(f)$ is in every $\mathcal{H}^s$. The image of $\grad$ is therefore strictly contained in the kernel of $\curl$. Similarly the image of $\curl$ is also strictly contained in the kernel of $\div$. The Sobolev-deRham complex is therefore not exact at either place.

Thus, in the case of the Sobolev limits, there is an obstruction to a kernel being an image, other than the existence of torsion elements ($A^3/C$ is torsion free, where $C$ is the submodule of $A^3$ generated by the rows of $\curl$, as in Examples 2.1 and 2.2). \hspace*{\fill}$\square$
\vspace{2mm}

Some results require only that $\F$ be an  injective $A$-module, and are thus true for $\Sr$, as well as for $\Dr$ and $\Cinf$. Proposition 3.5 on the elimination problem is one such result. In $\D$, elimination is not always possible (Example 3.4). 

%\noindent Example 5.2: If $\F$ is not injective, then the projection of a  system may not be a system. For instance, let $A=\mathbb{C}[\dbydt]$, and let $\F = \D$. Let $\pi_2: \D^2 \rightarrow \D$ be the projection onto the second factor. Let $P \subset A^2$ be the cyclic submodule generated by $(\dbydt, -1)$. Then $\ker_\D(P) = \{(f, {\dbydt}f)~|~ f \in \D \}$,  but $\pi_2(\ker_\D(P))  = \{{\dbydt}f ~|~ f \in \D \}$ is not a differential kernel in $\D$, as ${\dbydt}: \D \rightarrow \D$ is not surjective.

%Thus the inclusion $\pi_2(\ker_\F(P)) \subset \ker_\F(i_2^{-1}(P))$ in Proposition 2.5 can be strict when $\F$ is an not an injective $A$-module. In general, the obstruction to the projection $\pi_2(\ker_\F(P))$ being a differential kernel, lies in \index{Ext}${\ext}^1_A(A^p/\pi_1(P), ~\F)$. \\

We point out a couple of other pathologies. For instance,  Lemma 4.2 improves containment in Lemma 2.1(ii) to equality, when $\F$ is injective. In general,  containment is strict. 
\vspace{1.5mm}

\noindent Example 6.3: Let $A = \mathbb{C}[\dbydt]$ and
$\F = \D$. Let
$P_1$ and $P_2$ be cyclic submodules of
$A^2$ generated by $(1,0)$ and $(1,~-\dbydt)$, respectively. 
Then $P_1 \cap P_2$ is the 0 submodule,
so that $\ker_{\D}(P_1 \cap P_2)$ is
all of $\D^2$. On the other hand,
$\ker_{\D}(P_1)=\{(0, f) ~|~ f \in \D \}$ and
$\ker_{\D}(P_2)=\{({\dbydt}g,~g)~|~ g\in \D\}$.
Thus an element $(u,v)$ in $\D^2$ 
is in $\ker_{\D}(P_1)+ ~\ker_{\D}(P_2)$
only if $u={\dbydt}(g), ~v=f+g$, 
where $f$ and $g$ are arbitrary elements in $\D$. Let now $u$ be any
(nonzero) non-negative compactly supported smooth function.
Then $(u,0)$, which is in
$\ker_{\D}(P_1 \cap P_2)$, is, however, not in
$\ker_{\D}(P_1)+\ker_{\D}(P_2)$,
as $g(t) = \int_{-\infty}^{t}dg
= \int_{-\infty}^{t}u~dt$ is not compactly supported. 

Thus $\ker_\D(P_1) + \ker_\D(P_2) \subsetneq \ker_\D(P_2)$. In fact, it turns out that that $\ker_\D(P_1) + \ker_\D(P_2)$ is not the kernel of any differential operator (viz., Corollary 7.3 below). 

In general there is an obstruction to a sum of two systems in $\F$ being a system, and this obstruction can be located in \index{Ext}${\ext}^1_ A (A^k /(P_1 + P_2 ), \F)$ (viz. Lemma 4.2), which, for an injective $\F$, equals zero.   \\

The sum of two systems, if it be a system, is the one obtained by connecting the two \index{parallel connection} in `parallel'. Their intersection is the system corresponding to a \index{series connection}`series' connection. By Lemma 2.1(ii), the series connection of infinitely many systems is again a system (defined by the sum of the submodules defining these systems), and Lemma 4.2 states that the parallel connection of finitely many systems is a system, when the signal space is an injective $A$-module. However, the example below shows that the parallel connection of infinitely many systems need not be one,  even when $\F$ is an injective cogenerator. 
\vspace{2mm}

\noindent Example 6.4: Let $A=\C[{\dbydt}]$, and let $\F = \Cinf(\mathbb{R})$. Let $I$ be the principal ideal $({\dbydt})$, and let $P_i = I^i$, $i \geqslant 0$. Clearly, $\bigcap_{i \geqslant 0}P_i = 0$. It follows that $\ker_{\Cinf}(\bigcap_{i \geqslant 0}P_i) = \Cinf$. However, $\ker_{\Cinf}(P_i)$ equals the $\C$-vector space of polynomials of degree less than $i$, and hence the sum $\sum_{i \geqslant 0}\ker _{\Cinf}(P_i)$ equals the space of all polynomials, which is strictly contained in $\ker_{\Cinf}(\bigcap_{i \geqslant 0}P_i)$. It is clear that this infinite sum of systems is not a system at all.

\newpage

\section{The Nullstellensatz for Systems of PDE}

\vspace{.2cm}

We have studied in some detail controllable systems defined in a signal space $\F$ which is an injective cogenerator, specifically $\Dr$ and $\Cinf$. On the other hand, if $\F$ is a flat $A$-module, for instance $\D$, $\Er$, or $\Sc$, then every differential kernel is an image by the equational crterion of flatness \cite{ma}; thus every system in such a signal space admits a vector potential and is controllable. We now focus on the space $\Sr$ of tempered distributions. It is an injective $A$-module, but not a cogenerator (namely Remark 2.5). The system $\ker_{\Sr}(P)$, defined by $P \subset A^k$, is controllable if $A^k/P$ is torsion free, but this is no longer a necessary condition (Theorem 3.1(i) and the remark following it).  
\vspace{2mm}

The first problem we encounter when studying systems in a space $\F$ which is not an injective cogenerator is that Proposition 2.6 is not valid, and there is  no longer a bijective correspondence between systems in $\F^k$ and submodules of $A^k$ (Example 2.4). This prompted the definition of $\M(\ker_\F(P))$, for $P$ a submodule of $A^k$, in Remark 2.7, and we recall it now.

\begin{definition} Let $P \subset A^k$. Then $\M(\ker_\F(P))$ is the submodule of all $p(D) \in A^k$ such that the kernel of the map $p(D): \F^k \rightarrow \F$ contains $\ker_\F(P)$. It is the \index{closure, Willems}Willems closure (or simply, the closure) of $P$ with respect to $\F$. We denote it by $\bar{P}_\F$ (or simply by $\bar{P}$ if $\F$ is clear from the context). $P$ is said to be closed with respect to $\F$ if $P = \bar{P}_\F$. 
\end{definition}
Thus for instance, Example 2.4 states that $(\overline{D + \imath}) = (1)$ with respect to $\Sr$. We have also observed that every submodule $P$  is closed when $\F$ is an injective cogenerator.
\vspace{1mm}

We give another useful description of $\bar{P}_\F$. As $\homo_A(A^k/P, ~\F) \simeq \ker_\F(P) = \ker_\F(\bar{P}_\F) \simeq \homo_A(A^k/\bar{P}_\F, ~\F)$, it follows that every $\phi: A^k \rightarrow \F$ which vanishes on $P$ also vanishes on $\bar{P}_\F$. Conversely, if for $p \in A^k$, there is some $\phi: A^k \rightarrow \F$ which vanishes on $P$, but $\phi(p) \neq 0$, then this $p \notin \bar{P}_\F$.
\vspace{1.5mm}

If $P_1$ and $P_2$ are submodules of $A^k$, then clearly $\ker_\F(P_1) = \ker_\F(P_2)$ if and only if $\bar{P_1} = \bar{P_2}$. We collect a few other elementary observations below.

\begin{lemma} Let $P \subsetneq A^k$, and  $\F$,  an $A$-submodule  of $\Dr$. Then $\bar{P}_\F$ is contained in the 0-primary component $P_0$ of $P$ (in the notation of (9) of Section 4).
\end{lemma}
\noindent Proof: If $A^k/P$ is a torsion module, then $P_0 = A^k$, and there is nothing to be done. Otherwise, let $p \in A^k \setminus P_0$, and let $[p]$ denote its class in  $A^k/P_0$ or in $A^k/P$. It suffices to
show that there is a map $\phi: A^k/P \rightarrow \F$ which does not vanish at $[p]$.

As $A^k/P_0$ is finitely generated and torsion free, it embeds into a free $A$-module, say $A^r$.  Let the image of $[p]$ in $A^r$ be $q(D) = (q_1(D), \ldots, q_r(D))$, where we assume, without loss of generality, that $q_1(D)$ is nonzero. As $\F$ is a  faithful $A$-module, namely Corollary 3.1, there is an $f \in \F$ such that $q_1(D)f \neq 0$. Let $\psi: A \rightarrow \F$ be defined by $\psi(1) = f$. Then the composition $A^k/P_0 \hookrightarrow A^r \stackrel{\pi_1}{\longrightarrow} A \stackrel{\psi}{\longrightarrow} \F$, where $\pi_1$ is the projection to the first factor, is a map $A^k/P_0 \rightarrow \F$ which does not vanish at $[p]$. As $P \subset P_0$, this  map induces a map $\phi: A^k/P \rightarrow \F$ which also does not vanish at $[p]$.
\hspace*{\fill}$\square$
 \begin{lemma} Let $\{B_i\}$ be a collection of systems in $\F^k$, then  $\sum_i \M(B_i) \subset \M(\bigcap_iB_i)$, and $\M(\sum_i B_i) =  \bigcap_i \M(B_i)$.
\end{lemma}
\noindent Proof:  This is analogous to Lemma 2.1(ii). \hspace*{\fill}$\square$
\begin{lemma} For any signal space $\F$, $\ker_\F \circ \M$ is the identity map on systems. 
\end{lemma}
\noindent Proof: Let $B = \ker_\F(P)$; then $\M(B) = \bar{P}$, and $\ker_\F(\bar{P}) = B$. \hspace*{\fill}$\square$
\begin{corollary} For any $\F$, the closure of $P$ with respect to $\F$ is closed with respect to $\F$.
\end{corollary}
\noindent Proof: $\bar{\bar{P}} = \M\ker_\F(\M\ker_\F(P)) = \M(\ker_\F \circ \M(\ker_\F(P)) = \M\ker_\F(P) = \bar{P}$. \hspace*{\fill}$\square$
\vspace{2mm}

While the inclusion in Lemma 7.2 might be strict, there is equality of closures.

\begin{lemma} Suppose $\{B_i\}$ is a collection of systems in $\F^k$, then  the closure of $\sum_i \M(B_i)$ with respect  to $\F$ equals $\M(\bigcap_iB_i)$.
\end{lemma}
\noindent Proof: By Lemmas 2.1(ii) and 7.3, $\ker_\F(\sum_i\M(B_i)) = \bigcap_i\ker_\F(\M(B_i)) = \bigcap_iB_i = \ker_\F\M(\bigcap_iB_i)$. Hence the closure of  $\sum_i\M(B_i)$ equals the closure of $\M(\bigcap_iB_i)$, and this equals $\M(\bigcap_iB_i)$ by the above corollary. \hspace*{\fill}$\square$\\

By Lemma 7.3, there is now a bijective correspondence between systems in $\F^k$ and submodules of $A^k$ closed with respect to $\F$. {\it This is completely analogous to the familiar Galois correspondence between affine varieties in $\C^n$ and radical ideals of $\C[x_1, \ldots, x_n]$.} The notion of closure here is the analogue of the radical of an ideal, and its calculation is the analogue of the Hilbert Nullstellensatz.\\

Consider two signal spaces $\F_1 \subset \F_2$, and let $P$ be a submodule of $A^k$. Suppose $p \in A^k \setminus \bar{P}_{\F_1}$. Then there is a map $\phi: A^k \rightarrow \F_1$ vanishing on $P$ such that $\phi(p) \neq 0$, and hence $i \circ \phi: A^k \rightarrow \F_2$ does not also vanish at $p$, where the map $i$ is the inclusion of $\F_1$ in $\F_2$. Thus always $\bar{P}_{\F_2} \subset \bar{P}_{\F_1}$.

\begin{lemma} Let $\F_1 \subset \F_2$. Assume that for every $f \in \F_2$, $f \neq 0$, there exists a homomorphism $\phi_f: \F_2 \rightarrow \F_1$ such that $\phi_f(f) \neq 0$. Then $\bar{P}_{\F_1} = \bar{P}_{\F_2}$ for any submodule $P \subset A^k$. 
\end{lemma} 
\noindent Proof: Let $p \in A^k \setminus \bar{P}_{\F_2}$, and let $\phi:A^k \rightarrow \F_2$ be a map vanishing on $P$ such that $\phi(p) = f \neq 0$. Then $\phi_f \circ \phi:A^k \rightarrow \F_1$ vanishes on $P$ and $\phi_f \circ \phi(p) \neq 0$, and hence also $p \in A^k \setminus \bar{P}_{\F_1}$.
\hspace*{\fill}$\square$
\begin{lemma} Let ${\F}_1 \subset {\F}_2$, and $P$ a submodule of $A^k$. If $P$ is closed with respect to ${\F}_1$, then it is also closed with respect to $\F_2$. 
\end{lemma}
\noindent Proof: As always $\bar{P}_{\F_2} \subset \bar{P}_{\F_1}$, it follows that if $P = \bar{P}_{\F_1}$, then also $P = \bar{P}_{\F_2}$. \hspace*{\fill}$\square$ 
\begin{lemma} If $P \subset A^k$ is not closed with respect to a signal space $\F$, then $\ker_\F(P)$ is not dense in $\ker_{\Dr}(P)$.
\end{lemma}
\noindent Proof: By assumption, $P$ is strictly contained in its closure $\bar{P}_\F$. As $\ker_\F(P)$ equals $\ker_\F(\bar{P}_\F) = \ker_{\Dr}(\bar{P}_\F) \cap \F^k$, its closure in $(\Dr)^k$ 
is contained in $\ker_{\Dr}(\bar{P}_\F)$, a closed subspace of $({\Dr})^k$ that is strictly contained in $\ker_{\Dr}(P)$ (as every submodule of $A^k$ is closed with respect to $\Dr$).\hspace*{\fill}$\square$
\begin{lemma} Let $\{P_i\}$ be a collection of submodules of $A^k$, each closed with respect to $\F$. Then $\bigcap_i P_i$ is also closed with respect to $\F$.
\end{lemma}
\noindent Proof: As each $P_i$ is closed, it follows that
$\bigcap_iP_i \subset \M\ker_\F(\bigcap_iP_i) \subset \M(\sum_i\ker_\F(P_i)) = \bigcap_i\M\ker_\F(P_i) = \bigcap_iP_i$,
where the second inclusion follows from Lemma 2.1(ii). This implies equality everywhere, and thus that $\bigcap_iP_i = \M\ker_\F(\bigcap_iP_i)$.\hspace*{\fill}$\square$
\vspace{2.5mm}

A question more general than the above lemma is whether $\M\ker_\F(\bigcap_iP_i)$  equals $\bigcap_i\M\ker_\F(P_i)$ for an arbitrary collection of submodules $\{P_i\}$ of $A^k$. In other words, is the closure of an intersection equal to the intersection of the closures? (In our notation, is $\overline{\bigcap_iP_i} = \bigcap_i\bar{P_i}$~?) It is trivially true when $\F$ is an injective cogenerator, but in general it is only the containment  $\M\ker_\F(\bigcap_iP_i) \subset \bigcap_i\M\ker_\F(P_i)$ that holds. The other inclusion is not always true, as the following example demonstrates.
\vspace{2mm}

\noindent Example 7.1: Let $I$ be any nonzero proper ideal of $A$. Let $P_i = I^i, i \geqslant 0$, and let $\F = \D(\mathbb{R}^n)$, the space of smooth compactly supported functions. As each $P_i$ is a nonzero ideal, $\ker_\D(P_i) = 0$ by Paley-Wiener. Thus each $\bar{P_i} = A$, and hence $\bigcap_i\bar{P_i} = A$. On the other hand, $\bigcap_iP_i = 0$ (for instance, $A$ is an integral domain, hence its $I$-adic completion is Hausdorff by a theorem of Krull). Hence $\overline{\bigcap_iP_i} = 0$. 
\vspace{2mm}

The question therefore is whether the we have equality when the intersection is finite, i.e. does $\overline{\bigcap_{i=1}^mP_i} = \bigcap_{i=1}^m\bar{P_i}$ hold? This would be analogous to the fact that the radical of a finite intersection of ideals equals the intersection of the individual radicals. The only way to answer such a question seems to be (to this author) to first characterise the closures of submodules with respect to the  space $\F$, and we turn to this problem next.

\begin{proposition} Let $P$ be a proper submodule of $A^k$, and let $\F$ be $\D$, $\Er$, or $\Sc$. ~Let $P = P_0 \cap P_1 \cap \cdots \cap P_r$ be an irredundant primary decomposition of $P$ in $A^k$, where $P_0$ is a 0-primary submodule of $A^k$, and $P_i$,  $i \geqslant 1$, is a $p_i$-primary submodule, $p_i$ a nonzero prime. Then the closure $\bar{P}_\F$ of $P$ with respect to $\F$ equals $P_0$. Thus $P$ is closed with respect to $\F$ if and only if $A^k/P$ is torsion free. 
\end{proposition}
\noindent Proof: Suppose $P$ does not have a 0-primary component. Then $A^k/P$ is a torsion module, $\ker_{\Dr}(P)$ is uncontrollable, and there is no nonzero compactly supported or rapidly decaying decaying element in it (Proposition 4.2 and the remark following it).  Thus $\ker_\F(P) = 0$, and $\bar{P}_\F = A^k$.

Now assume that $P$ has a 0-primary component $P_0$. We have already noted that it is uniquely determined by $P$. We first show that $\ker_\F(P) = \ker_\F(P_0)$, and hence that $\bar{P}_\F = (\bar{P_0})_\F$.

As $\F$ is a flat $A$-module, $\ker_\F(P)$ is an image, say of the map $R(D): \F^{k_1} \rightarrow \F^k$. Consider the extension of this map to $\Dr$, namely $R(D): ({\Dr})^{k_1} \rightarrow  ({\Dr})^k$. Its image is equal to the controllable part of $\ker_{\Dr}(P)$, which equals $\ker_{\Dr}(P_0)$ by Theorem 3.1. Thus, the image $R(D)(\F^{k_1})$ is contained in $\ker_\F(P_0)$, and is hence equal to it. 

We now claim that $P_0$ is closed with respect to $\F$, so that $(\bar{P_0})_\F = P_0$. For if not, then by Lemma 7.7, $\ker_\F(P_0)$ is not dense in $\ker_{\Dr}(P_0)$. But this contradicts Proposition 2.2 because $\ker_{\Dr}(P_0)$ is controllable as $A^k/P_0$ is torsion free. \hspace*{\fill}$\square$\\

\noindent Remark 7.1: We can also describe the above closure as follows. Let $\pi: A^k \rightarrow A^k/P$ be the canonical surjection. Then the closure $\bar{P}_\F$~equals $\pi^{-1}(T(A^k/P)$, where $T(A^k/P)$ is the submodule of torsion elements of $A^k/P$ (similar descriptions for the closure in all the classical spaces in \cite{sell}). 
\vspace{1.5mm}

We turn now to the space $\Sr$ of tempered distributions. We use the following statement on primary decompositions.

Suppose $Q = \bigcap_{i = 1}^rP_i$ is a primary decomposition of $Q$ in $A^k$, where $P_i$ is $p_i$-primary. Suppose $I$ is an ideal such that $I \subset p_i$ for $i = q+1, \ldots, r$, and that $I$ is not contained in the other $p_i$. Then the ascending chain of submodules $(Q:I) \subset (Q:I^2) \subset \cdots $ stabilizes to the submodule $\bigcap_{i=1}^qP_i$, and hence this submodule is independent of the primary decomposition of $Q$.

\begin{proposition} Let $P$ be a submodule of  $A^k$, and let $\F = \Sr$. Let $P = P_1\cap \cdots \cap P_r$ be an irredundant primary decomposition of $P$ in $A^k$, where $P_i$ is $p_i$-primary. Let $\mathcal{V}(p_i)$, the variety of the ideal $p_i$ in $\C^n$,
contain real points for $i=1,\dots ,q$, and not for $i=q+1,\dots ,r$. Then the closure $\bar{P}_{\Sr}$ of $P$ with respect to $\Sr$ equals $P_1\cap \cdots \cap P_q$. Thus $P$ is closed with respect to $\Sr$ if and only if every $\mathcal{V}(p_i)$ contains real points.
\end{proposition}
\noindent Proof:   If we set $I = p_{q+1}\cap \cdots \cap p_r$ in the result quoted just above, then it follows that $P':=  P_1\cap \cdots \cap P_q$ is independent of the primary decomposition of $P$. 

We first claim that $\ker_{\Sr}(P')$ equals $\ker_{\Sr}(P)$. As $P \subset P'$, it suffices to show that $\ker_{\Sr}(P) \subset \ker_{\Sr}(P')$. If this is not true, then there is some $f$ in $\ker_{\Sr}(P)$, and a $p(D) \in P' \setminus P$, such that $p(D)f \neq 0$. However, for every $a(D)$ in the ideal $(P:p(D))$, $a(D)(p(D)f) = 0$. Taking Fourier transforms gives $a(\xi)(\widehat{p(D)f})(\xi) = 0$, hence $ (\supp(\widehat{p(D)f}))$ is contained in $\mathcal{V}(a) \cap \R^n$ (where $\supp$ denotes support). As $P=\bigcap_{i=1}^rP_i$, $(P:p(D)) = \bigcap_{i=1}^r(P_i:p(D))$. The submodule $P_i$ is $p_i$-primary in $A^k$, and so it follows that $\sqrt{(P_i:p(D))}$ either equals $A$ (if $p(D) \in P_i$), or it equals $p_i$ (if $p(D) \notin P_i$). Hence $\sqrt{(P:p(D))}$ is equal to the intersection of a subset of $\{p_{q+1}, \ldots, p_r\}$. It now follows that $ (\supp(\widehat{p(D)f}))$ is contained in $(\bigcup_{i=q+1}^r\V(p_i)) \cap \R^n$. By assumption, none of the $\V(p_i)$, $i = q+1, \ldots, r,$ intersects $\R^n$, hence the support of $\widehat{p(D)f}$ is empty. Thus $p(D)f = 0$, contrary to assumption. 

We now show that $P' = \bar{P}_{\Sr}$ . So let $p(D)$ be any element of $A^k \setminus P'$, and consider the exact sequence
\[ 0 \rightarrow A/(P':p(D)) \stackrel{p}{\longrightarrow} A^k/P' \stackrel{\pi}{\longrightarrow} A^k/P'+(p(D)) \rightarrow 0 \]
where the morphism $p$ above maps $[a(D)]$ to  $[a(D)p(D)]$, and $\pi$ is the canonical surjection. Applying the functor ${\homo}_A(- ~ ,\Sr)$ gives the sequence
\[
0\rightarrow {\homo}_A(A^k/P'+(p(D)),\Sr) \longrightarrow {\homo}_A(A^k/P',\Sr) \stackrel{p(D)}{\longrightarrow} {\homo}_A(A/(P':p(D)),\Sr) \rightarrow 0\]
which is exact because $\Sr$ is an injective $A$-module.
Observe now that $V((P':p(D)))$ is the union of some of the varieties $\V(p_1), \ldots, \V(p_q)$, hence by assumption there is a real point, say $\xi_0$ on it. Then the function $e^{\imath <\xi_0,\xi>}$ is tempered, and so belongs to the last term ${\homo}_A(A/(P':p(D)),~\Sr)$ above. It is therefore nonzero, and this implies that $\ker_{\Sr}(P'+p(D))$ is strictly contained in  $\ker_{\Sr}(P')$. The proposition follows.  \hspace*{\fill}$\square$\\

We can now answer the question following Example 7.1, when $\F$ is one of the classical spaces.

\begin{corollary} Let $P_i$, $i= 1, \ldots, m$, be submodules of $A^k$, and let $\F$ be a classical space. Then  $(\overline{\bigcap_{i=1}^mP_i})_\F = \bigcap_{i=1}^m\bar{P_i}_\F$.
\end{corollary}
\noindent Proof: By induction, it suffices to prove the statement for $m=2 $.
\vspace{1mm}

\noindent (i) When $\F = \Dr$ or $\Cinf$, the statement is trivial, because now every submodule of $A^k$ is closed.

As  the inclusion $ (\overline{P_1 \cap P_2})_\F \subset \bar{P_1}_\F \cap \bar{P_2}_\F$ is true for every $\F$, we need to show the reverse inclusion for $\Sr$, $\Sc$, $\Er$ and $\D$.

\vspace{1mm}
\noindent (ii) Let $P_1 = \bigcap_{i=1}^{r_1}Q_i$ and $P_2 = \bigcap_{j=1}^{r_2}Q'_j$ be irredundant primary decompositions of $P_1$ and $P_2$ in $A^k$ respectively, where $Q_i$ and $Q'_j$ are $p_i$-primary and $p'_j$-primary, respectively, $i = 1, \ldots, r_1, ~j = 1, \ldots, r_2$. Suppose that the varieties of $p_1, \ldots, p_{q_1}, p'_1, \ldots, p'_{q_2}$ contain real points, whereas those of $p_{q_1+1}, \ldots, p_{r_1}, p'_{q_2+1} \ldots, p'_{r_2}$ do not. By Proposition 7.2, $\bar{P_1}_{\Sr} \cap \bar{P_2}_{\Sr}  = ( \bigcap_{i=1}^{q_1}Q_i) \cap (\bigcap_{j=1}^{q_2}Q'_j)$. 

On the other hand, $(\bigcap_{i=1}^{r_1}Q_i) \cap (\bigcap_{j=1}^{r_2}Q'_j)$  is a  primary decomposition of $P_1 \cap P_2$, though perhaps not irredundant. An irredundant primary decomposition can however be obtained from it by omitting, if necessary, some of the $Q_i$ or $Q'_{j}$. Thus the set of associated primes of $P_1 \cap P_2$ is a subset of $\{p_1, \ldots, p_{r_1}, p'_1, \ldots, p'_{r_2}\}$. This implies that those associated primes of $P_1 \cap P_2$ whose varieties contain real points is a subset of $\{p_1, \ldots, p_{q_1}, p'_1, \ldots, p'_{q_2}\}$. Clearly then, by Proposition 7.2, $\bar{P_1}_{\Sr} \cap \bar{P_2}_{\Sr}  \subset (\overline{P_1 \cap P_2})_{\Sr}$, which proves the corollary for $\Sr$. 

\vspace{1mm}
\noindent (iii) Now  let $\F = \Sc$, $\Er$ or $\D$. Let $p \in A^k$ be in $\bar{P_1}_\F \cap \bar{P_2}_\F$. Then there exist nonzero $a_1$ and $a_2$ in $A$ such that $a_1p \in P_1$ and $a_2p \in P_2$ (by Remark 7.1). As $A$ is a domain, the product $a_1a_2 \neq 0$, and  $a_1a_2p \in P_1 \cap P_2$, which is to say that $p \in (\overline{P_1 \cap P_2})_\F$ ~.   \hspace*{\fill}$\square$\\

We can now complete the argument in Example 6.3 and exhibit two kernels whose sum is not a kernel.

\begin{corollary} Let $\F$ be a classical space, and $P_i$, $i = 1, \ldots, r$, submodules of $A^k$. Then $\ker_\F(\bigcap_{i=1}^rP_i)$ is the smallest kernel containing $\sum_{i=1}^r\ker_\F(P_i)$.
\end{corollary}
\noindent Proof: By Lemma 4.2, the two expressions are equal if $\F$ is an injective $A$-module. Now let $\F = \Sc$, $\Er$ or $\D$.
Suppose $\ker_\F(P)$ contains both $\ker_\F(P_1)$ and $\ker_\F(P_2)$, for some $P\subset A^k$. Then $\bar{P}_\F$ is contained in both $\bar{P_1}_\F$ and $\bar{P_2}_\F$, and hence in $\bar{P_1}_\F \cap \bar{P_2}_\F$~, which equals $(\overline{P_1 \cap P_2})_\F$ by the above corollary. It  now follows that 
\[ {\ker}_\F(P_1\cap P_2) = {\ker}_\F((\overline{P_1 \cap P_2})_\F) \subset {\ker}_\F(\bar{P}_\F) = {\ker}_\F(P).
\]
Hence, any kernel which contains both $\ker_\F(P_1)$ and $\ker_\F(P_2)$ also contains $\ker_\F(P_1\cap P_2)$. Thus $\ker_\F(P_1 \cap P_2)$ is the smallest kernel that contains $\ker_\F(P_1) + \ker_\F(P_2)$. The lemma follows by induction. \hspace*{\fill}$\square$\\

We improve Theorem 3.1(i) to now provide a necessary condition for a system in $\Sr$ to be controllable. 

\begin{proposition} Let $P$ be a submodule of $A^k$. The system $\ker_{\Sr}(P)$ is controllable if and only if $A^k/\bar{P}_{\Sr}$ is torsion free, and hence if and only if the system admits a vector potential.
\end{proposition}
\noindent Proof: Suppose $P = P_0 \cap P_1\cap \cdots \cap P_r$ is an irredundant primary decomposition of $P$ in $A^k$, where $P_0$ is 0-primary, and $P_i$ is $p_i$-primary, $p_i$ nonzero prime, for $i \geqslant 1$. Let $\mathcal{V}(p_i)$
contain real points for $i=0,\dots ,q$, and not for $i=q+1,\dots ,r$. Then the closure $\bar{P}_{\Sr}$ of $P$ with respect to $\Sr$ equals $P_0 \cap P_1 \cap \cdots \cap P_q$. By definition $\ker_{\Sr}(P) = \ker_{\Sr}(\bar{P}_{\Sr})$, and this in turn is equal to $\sum_{i=0}^q \ker_{\Sr}(P_i)$ by Lemma 4.2. In this sum, only the summand $\ker_{\Sr}(P_0)$ is controllable. Hence $\ker_{\Sr}(P)$ is controllable if and only if the index $q$ is equal to 0, which is to say that $A^k/\bar{P}_{\Sr}$ is torsion free.  As $\Sr$ is injective, $\ker_{\Sr}(\bar{P}_{\Sr})$ admits a vector potential by  Theorem 3.1(i).

Suppose now that there is no 0-primary component in the primary decomposition of $P$, so that $A^k/\bar{P}_{\Sr}$ is either torsion or $\bar{P}_{\Sr} = A^k$. In the first case, $\ker_{\Dr}(P)$ is uncontrollable by Proposition 4.1, and has no nonzero compactly supported elements in it by Proposition 4.2. It follows that $\ker_{\Sr}(P)$ also does not have nonzero compactly supported elements, and therefore cannot admit a vector potential by Proposition 2.2. In the latter case, $\ker_{\Sr}(P) = 0$ is trivially controllable, whose vector potential is the 0 map (as in Remark 3.1). 
\hspace*{\fill}$\square$\\

We put together Theorem 3.1(ii), Proposition 7.1 and Proposition 7.3 in the statement below. {\it It is the central result of these notes.}

\begin{theorem} Let $\F$ be a classical signal space. Let $P$ be a submodule of $A^k$ and let $\bar{P}_\F \subset A^k$ be its closure with respect to $\F$. Then the system $\ker_\F(P)$ is controllable if and only if $A^k/\bar{P}_\F$ is torsion free, and hence if and only if the system admits a vector potential.
\end{theorem}

We conclude our discussion of controllability with an important question. Theorem 7.1 states that in the classical spaces $\Dr$, $\Cinf$, $\Sr$, $\Sc$, $\Er$ and $\D$, the notion of controllability is equivalent to the existence of a vector potential. We ask if this is so over any signal space $\F$? In other words, is Kalman's notion of a controllable system, suitably generalised, nothing more - nor less - than the possibility of describing the dynamics of the system by means of a vector potential? \\

We continue our discussion of the Nullstellensatz for PDE with the calculation of the closure of a submodule of $A^k$ with respect to spaces of periodic functions. These results are from \cite{dps}; they highlight the importance of rational and integral points on algebraic varieties in the study of PDE, just as Proposition 7.2 does for real points.
\vspace{1.5mm}

Consider the lattice $2\pi \mathbb{Z}^n$ and the space of smooth functions on $\R^n$ periodic with respect to it. We may identify these functions with smooth functions on the torus $\t = \R^n/2\pi \mathbb{Z}^n$, and hence we denote this space by $\Cinf(\t)$. An element $f:\R^n \rightarrow \mathbb{C}$ 
in $\mathcal{C}^\infty (\t)$ is represented by its Fourier series $f(x) = 
\sum _{\xi \in \mathbb{Z}^n} c_\xi
e^{\imath<\xi,x>}$, where the coefficients
$c_\xi \in \mathbb{C}$ are required to satisfy the following estimate:
for every integer $k\geq 1$ there exists a constant $C_k>0$ such that
$|c_\xi|\leq \frac{C_k}{(1+\sum _{j=1}^n|\xi_j|)^k}$ for all $\xi$. 
%The space of distributions on $\t$ has a similar description, however with different requirements on the absolute values $|c_a|$.)
The space $\Cinf(T)$ is Fr\'echet, and a topological $A$-module. 

For positive integers $N_1$ dividing $N_2$,
the natural inclusion $\mathcal{C}^\infty (\mathbb{R}^n/2\pi N_1\mathbb{Z}^n) \hookrightarrow \mathcal{C}^\infty (\mathbb{R}^n/ 2\pi N_2\mathbb{Z}^n)$ is an $A$-module morphism and identifies the first 
space with a closed subspace of the second. The corresponding directed set $\{\mathcal{C}^\infty (\mathbb{R}^n/2\pi N\mathbb{Z}^n)~|~N \in \mathbb{N}\}$ of Fr\'echet spaces defines the (strict) direct limit $\underset{\rightarrow }{\lim} \ \mathcal{C}^\infty (\mathbb{R}^n/2\pi N\mathbb{Z}^n)$, and is a locally convex bornological and barrelled topological vector space. It can be identified with functions on the inverse limit $\underset{\leftarrow}{\lim }\ \mathbb{R}^n/2\pi N\mathbb{Z}^n$, called a `protorus' in \cite{dps} and denoted $\pt$. Hence we denote the above direct limit by $\mathcal{C}^\infty (\pt)$. The elements of $A$ act continuously on it so that $\mathcal{C}^\infty(\pt)$ is a topological $A$-module. An element $f:\R^n \rightarrow \C$ in it is
represented by its Fourier series $f(x)=\sum _{\xi \in \mathbb{Q}^n} c_\xi e^{\imath<\xi,x>}$, where the support of $f$, i.e.,
$\{\xi \in \mathbb{Q}^n|\ c_\xi\neq 0\}$, is a subset of
$\frac{1}{N}\mathbb{Z}^n$ for some integer $N\geq 1$ depending on $f$, and where the $|c_\xi|$ satisfies the same requirement of rapid decrease as above. By uniform convergence of the partial sums of a Fourier series, it follows that 
\[ p(D)(\sum _{\xi \in \mathbb{Q}^n} c_\xi
e^{\imath<\xi,x>})=
\sum _{\xi \in \mathbb{Q}^n} c_\xi p(\xi)
e^{\imath<\xi,x>}.\]

The space $\mathcal{C}^\infty(\pt)$ (and hence also $\mathcal{C}^\infty(\t)$) is not an injective $A$-module, indeed it is not even divisible. For instance, the image of $D_1: \mathcal{C}^\infty (\pt)\rightarrow \mathcal{C}^\infty (\pt)$ consists of the 
elements $f$ with support in $\{(\xi_1,\dots ,\xi_n)\in \mathbb{Q}^n|\ 
\xi_1\neq 0\}$, and is therefore not surjective. 
These spaces are also not flat $A$-modules (the kernel of the above map cannot be an image).

We enrich the above spaces of periodic functions by adjoining the coordinate functions $x_1, \ldots, x_n$ to obtain the subalgebras $\Cinf(\t)[x_1, \ldots, x_n]$ and $\Cinf(\pt)[x_1, \ldots, x_n]$ of $\Cinf(\R^n)$. It is easy to see that for $n = 1$, $\Cinf(\t)[x]$ and $\Cinf(\pt)[x]$ are injective $\mathbb{C}[\dbydt]$-modules. However for $n > 1$, $\Cinf(\t)[x_1, \ldots, x_n]$ and $\Cinf(\pt)[x_1, \ldots, x_n]$ are not divisible $A$-modules, \cite{dps}. Neither are they flat $A$-modules. Nonetheless, we determine the PDE Nullstellensatz in these spaces as follows. 

Let $\Cinf(\t)_{\fin}$ be the $A$-submodule of $\Cinf(\t)$ consisting of those elements which are finitely supported, i.e. those elements whose Fourier series expansion is a finite sum. Similarly define the submodule $\Cinf(\pt)_{\fin} \subset \Cinf(\pt)$. We quote the following result (Proposition 2.1) from \cite{dps}.
\begin{proposition} The subalgebras $\Cinf(\t)_{\fin}[x_1, \ldots, x_n]$ and $\Cinf(\pt)_{\fin}[x_1, \ldots, x_n]$ of $\Cinf(\R^n)$ are injective $A$-modules (they are the injective hulls of $\Cinf(\t)_{\fin}$ and $\Cinf(\pt)_{\fin}$ respectively).
\end{proposition}

\begin{lemma} The signal spaces $\mathcal{C}^\infty (\t)_{\fin}[x_1,\dots ,x_n]$ and $\mathcal{C}^\infty(\t) [x_1,\dots ,x_n]$
define the same closure for every submodule of $A^k$. Similarly, the spaces $\mathcal{C}^\infty (\pt)_{\fin}[x_1,\dots ,x_n]$ and $\mathcal{C}^\infty(\pt) [x_1,\dots ,x_n]$ define the same closure.
\end{lemma}
\noindent Proof: For an $a \in \mathbb{Z}^n$, define the homomorphism
\[\pi_a:\mathcal{C}^\infty(\t) [x_1,\dots ,x_n]\rightarrow  \mathcal{C}^\infty(\t)_{\fin}[x_1,\dots ,x_n]\]
 by $\pi_a (\sum _{\xi \in \mathbb{Z}^n}
p_\xi(x)e^{\imath<\xi,x>} )= p_a(x)e^{\imath<a,x>}.$ The first statement now follows from Lemma 7.5. The second statement follows similarly.
\hspace*{\fill}$\square$
\begin{proposition}  (i) Let $P$ be a submodule of $A^k$. Let 
$P = P_1\cap \dots \cap P_r$ be an irredundant primary decomposition of $P$ in $A^k$, where $P_i$ is $p_i$-primary. Let $\mathcal{V}(p_i)$, the variety defined by $p_i$,
contain rational points for $i=1,\ldots,q$ and not for $i=q+1,\ldots,r$. Then the closure $\bar{P}$ with respect to 
$\F = \mathcal{C}^\infty(\pt)[x_1,\dots ,x_n]$ is equal to $P' = P_1\cap \cdots \cap P_q$. Thus $P$ equals its closure with respect to $\F$ if and only if every $\mathcal{V}(p_i)$ contains rational points.

\noindent (ii) Let $P_1\cap \dots \cap P_r$ be an irredundant primary decomposition of $P$ as in (i). Let $\mathcal{V}(p_i)$
contain integral points for $i=1,\ldots,q$ and not for $i=q+1,\ldots,r$. Then the closure $\bar{P}$ with respect to 
$\F = \mathcal{C}^\infty(\t)[x_1,\dots ,x_n]$ is equal to $P_1\cap \cdots \cap P_q$. Thus $P$ equals its closure if and only if every $\mathcal{V}(p_i)$ contains integral points.
\end{proposition}
\noindent Proof: (i) By Lemma 7.9, it suffices to take $\F = \mathcal{C}^\infty (\pt)_{\fin}[x_1,\dots ,x_n]$. 

We first need to show that $\ker_\F(P) = \ker_\F(P')$, and this follows exactly as the first half of the proof of Proposition 7.2 with $\mathbb{Q}^n$ replacing $\R^n$ in it (for if $f=\sum_{\xi\in\mathbb{Q}^n}p_\xi(x)e^{\imath<\xi,x>}$, a finite sum, then it is a tempered distribution and its Fourier transform $\hat{f}$ equals $\sum_{\xi \in \mathbb{Q}^n} p_\xi(D)\delta_\xi$, where $\delta_\xi$ is the Dirac distribution supported at $\xi$, so that the support of $\hat{f}$ is contained in $\mathbb{Q}^n$).
%If this is not true, then there is some $f$ in $\ker_\F(P)$, and a $p(D) \in P' \setminus P$, such that $p(D)f \neq 0$. However, for every $a(D)$ in the ideal $(P:p(D))$, $a(D)(p(D)f) = 0$. Taking Fourier transforms gives $a(\xi)(\widehat{p(D)f})(\xi) = 0$, hence $ (\supp(\widehat{p(D)f}))$ is contained in $\mathcal{V}(a) \cap \mathbb{Q}^n$. As $P=\bigcap_{i=1}^rP_i$, $(P:p(D)) = \bigcap_{i=1}^r(P_i:p(D))$. The submodule $P_i$ is $p_i$-primary in $A^k$, and so it follows that $\sqrt{(P_i:p(D))}$ either equals $A$ (if $p(D) \in P_i$), or it equals $p_i$ (if $p(D)\notin P_i$). Hence $\sqrt{(P:p(D))}$ is equal to the intersection of a subset of $\{p_{q+1}, \ldots, p_r\}$. 
%It then follows, again as in the proof of Proposition 7.2 that $ (\supp(\widehat{p(D)f}))$ is contained in $(\bigcup_{i=q+1}^r\V(p_i)) \cap \imath\mathbb{Q}^n$. By assumption, none of the $\V(p_i)$, $i = q+1, \ldots, r,$ intersects $\mathbb{Q}^n$, hence the support of $\widehat{p(D)f}$ is empty. Thus $p(D)f = 0$, contrary to assumption. 

The proof that $P' = \bar{P}_\F$ also follows just as the second half of the proof of Proposition 7.2 (with a rational point replacing the real point $\xi_0$). The second exact sequence in that proof remains exact because the space $\mathcal{C}^\infty (\pt)_{\fin}[x_1,\dots ,x_n]$ is injective.

\vspace{2mm}

\noindent (ii) It suffices to take $\F = \mathcal{C}^\infty (\t)_{\fin}[x_1,\dots ,x_n]$. Now the proof is identical to the proof of (i) with $\mathbb{Z}$ replacing $\mathbb{Q}$.
\hspace*{\fill}$\square$\\

We collect the results of Propositions 2.6, 7.1, 
7.2 and 7.5 in the statement below. 

\begin{theorem} (Nullstellensatz for PDE) (i) Every submodule of $A^k$ is closed with respect to $\Dr$ and $\Cinf$.

\noindent (ii) Let $P$ be a submodule of  $A^k$, and let $\F = \Sr$. Let $P = P_1\cap \cdots \cap P_r$ be an irredundant primary decomposition of $P$ in $A^k$, where $P_i$ is $p_i$-primary. Let $\mathcal{V}(p_i)$, the variety of the ideal $p_i$ in $\C^n$,
contain real  points for $i=1,\dots ,q$, and not for $i=q+1,\dots ,r$. Then the closure $\bar{P}_{\Sr}$ of $P$ with respect to $\Sr$ equals $P_1\cap \cdots \cap P_q$. Thus $P$ is closed with respect to $\Sr$ if and only if every $\mathcal{V}(p_i)$ contains real points.

\noindent (iii) Let $P \subset A^k$, and let $\F = \mathcal{C}^\infty(\pt)[x_1,\dots ,x_n]$. Let $P = P_1\cap \cdots \cap P_r$ be an irredundant primary decomposition as in (ii). Let $\mathcal{V}(p_i)$, the variety of the ideal $p_i$ in $\C^n$,
contain rational  points for $i=1,\dots ,q$, and not for $i=q+1,\dots ,r$. Then the closure $\bar{P}_\F$ of $P$ with respect to $\F$ equals $P_1\cap \cdots \cap P_q$. Thus $P$ is closed with respect to $\F$ if and only if every $\mathcal{V}(p_i)$ contains rational points.

\noindent (iv) Let $P \subset A^k$, and let $\F = \mathcal{C}^\infty(\t)[x_1,\dots ,x_n]$. Let $P = P_1\cap \cdots \cap P_r$ be an irredundant primary decomposition as above. Let $\mathcal{V}(p_i)$, the variety of the ideal $p_i$ in $\C^n$,
contain integral  points for $i=1,\dots ,q$, and not for $i=q+1,\dots ,r$. Then the closure $\bar{P}_\F$ of $P$ with respect to $\F$ equals $P_1\cap \cdots \cap P_q$. Thus $P$ is closed with respect to $\F$ if and only if every $\mathcal{V}(p_i)$ contains integral points.

\noindent (v) Let $P \subset A^k$, and let $\F$ be $\D$, $\Er$, or $\Sc$. ~Let $P = P_0 \cap P_1 \cap \cdots \cap P_r$ be an irredundant primary decomposition of $P$ in $A^k$, where $P_0$ is a 0-primary submodule of $A^k$, and $P_i$,  $i \geqslant 1$, is a $p_i$-primary submodule, $p_i$ a nonzero prime. Then the closure $\bar{P}_\F$ of $P$ with respect to $\F$ equals $P_0$. Thus $P$ is closed with respect to $\F$ if and only if $A^k/P$ is torsion free. 
\end{theorem}
\vspace{2mm}

As we have already explained, these results are analogues of Hilbert's Nullstellensatz, and the closure that we have defined is analogous to the radical of an ideal. 
\vspace{1mm}

In conclusion, we expand briefly on the above comment. Consider the category $\frak{D}$ of differential kernels whose objects are $\{\homo_A(-, \F)\}$, where the first factor is a finitely generated $A$-module and $\F$ a fixed $A$-submodule of $\Dr$, and whose morphisms are given by composing either with $A$-module maps between these finitely generated modules or with $A$-module endomorphisms of $\F$. This is the category that we have studied in these notes, and the results here suggest that it bears many formal similarities to the category of affine algebraic sets. We now ask: is it possible to `patch' together these objects in $\frak{D}$ to get a larger category just as we get the category of algebraic varieties from the category of affine algebraic sets? If so, it would then be psossible to extend the theory of control explained here for $\frak{D}$ to this larger category of systems.

\newpage 
\section{Acknowledgement}
\vspace{.5cm}

These notes are a modification of \cite{stek}. It now includes a discussion on the `achievable subspaces' of the space of solutions of the Maxwell equations, extra material on the PDE Nullstellensatz, and many typos and infelicities weeded out. %I am grateful to Alex Martsinkovsky for including them in his volume on Functor Categories. 
I am very grateful to Lev Lokutsievsky for his invitation to lecture at the Steklov Institute. I had earlier lectured on some of this material at IIT Bombay and Padova University, and I thank Madhu Belur and Maria Elena Valcher for their invitations.

I am grateful to Diego Napp, Harish Pillai, Marius van der Put and Paula Rocha with whom I collaborated on some of the work reported here. I am also grateful to many colleagues for many discussions, especially Arul Shankar for his quick proof of Theorem 3.2.

%\newpage

%\printindex{}

\end{document}